\documentclass[10pt]{article}

\usepackage{xr}
\usepackage[doi=false,isbn=false,url=false,
	hyperref=auto,
	sorting=nyt,
	maxnames=10,
	maxcitenames=3,
	backend=biber,
	texencoding=auto,
	giveninits=true,
	block=space,
	style=numeric,
]{biblatex}

\DefineBibliographyStrings{english}{%
	andothers = {et al}
} 
 
\DeclareFieldFormat{eprint:arxiv}{%
	\href{https://arxiv.org/abs/#1}{\emph{arXiv:\allowbreak #1}}}

\DeclareFieldFormat{eprint:mrnumber}{%
	\href{http://www.ams.org/mathscinet-getitem?mr=MR#1}{MR\allowbreak #1}}

\DeclareFieldFormat{eprint:jstor}{%
	\href{http://www.jstor.org/stable/#1}{JSTOR\allowbreak #1}}

\DeclareFieldFormat{eprint:customeprint}{%
	\em Available at \upshape \ttfamily \href{https://#1}{#1}}

\DeclareFieldFormat{eprint:online}{%
	Available \href{https://#1}{online}}

\DeclareFieldFormat{eprint:inprep}{%
	\emph{In preparation}}

\DeclareFieldFormat{eprint:manual}{%
	\emph{#1}}

\DeclareFieldFormat{eprint:onarxiv}{%
	\emph{#1}}

\DeclareFieldFormat{eprint:toappear}{%
	\mbox{\emph{To appear}}}

\DeclareFieldFormat{eprint:accepted}{%
	\emph{Accepted at {#1}; to appear}}

\DeclareSourcemap{
	\maps[datatype=bibtex]{
		\map{
			\step[fieldsource=mrnumber,		fieldtarget=eprint, final]
			\step[fieldset=eprinttype,		fieldvalue=mrnumber]
		}
		\map{
			\step[fieldsource=arxiv,		fieldtarget=eprint, final]
			\step[fieldset=eprinttype,		fieldvalue=arxiv]
		}
		\map{
			\step[fieldsource=jstor,		fieldtarget=eprint, final]
			\step[fieldset=eprinttype,		fieldvalue=jstor]
		}
		\map{
			\step[fieldsource=customeprint,	fieldtarget=eprint, final]
			\step[fieldset=eprinttype,		fieldvalue=customeprint]
		}
		\map{
			\step[fieldsource=online,		fieldtarget=eprint, final]
			\step[fieldset=eprinttype,		fieldvalue=online]
		}
		\map{
			\step[fieldsource=inprep,		fieldtarget=eprint, final]
			\step[fieldset=eprinttype,		fieldvalue=inprep]
		}
		\map{
			\step[fieldsource=manual,		fieldtarget=eprint, final]
			\step[fieldset=eprinttype,		fieldvalue=manual]
		}
		\map{
			\step[fieldsource=onarxiv,		fieldtarget=eprint, final]
			\step[fieldset=eprinttype,		fieldvalue=onarxiv]
		}
		\map{
			\step[fieldsource=toappear,		fieldtarget=eprint, final]
			\step[fieldset=eprinttype,		fieldvalue=toappear]
		}
		\map{
			\step[fieldsource=accepted,		fieldtarget=eprint, final]
			\step[fieldset=eprinttype,		fieldvalue=accepted]
		}
	}
}

\DeclareFieldFormat{doi}{%
	\href{https://doi.org/#1}{DOI}%
}

\DeclareFieldFormat{url}{%
	\href{https://#1}{\texttt{#1}}%
}

\DeclareFieldFormat{comments}{%
	\emph{#1}%
} 
 
\DeclareFieldFormat
	[article,inbook,incollection,inproceedings,patent,thesis,unpublished]
	{title}{{#1\isdot}}

\DeclareFieldFormat
	[article,book,inbook,incollection,inproceedings,patent,thesis,unpublished, online]
	{date}{\mkbibparens{#1}}
\DeclareFieldFormat
	[article,book,inbook,incollection,inproceedings,patent,thesis,unpublished]
	{volume}{\mkbibbold{#1}}
\DeclareFieldFormat
	{pages}{\mkbibparens{#1}}

\DeclareFieldFormat{edition}{%
	\ifinteger{#1}
	{\mkbibordedition{#1}~\bibstring{edition}\addcomma}
	{#1~\bibstring{edition}\addcomma}}

\renewbibmacro*{volume+number+eid}{%
	\printfield{volume}%
\setunit*{\adddot}%
	\printfield{number}%
\setunit{\space}%
	\printfield{eid}%
}

\DeclareBibliographyDriver{article}{%
	\usebibmacro{bibindex}%
	\usebibmacro{begentry}%
	\usebibmacro{author}%
\setunit{\addspace}%
	\usebibmacro{date}%
\newunit\newblock
	\usebibmacro{title}%
\newunit\newblock
	\printfield{journaltitle}\addperiod%
\setunit*{\addspace}%
	\usebibmacro{volume+number+eid}%
\setunit*{\addspace}\newblock
	\printfield{pages}
\setunit*{\addspace}
	\printfield{comments}%
	\usebibmacro{eprint}%
\setunit*{\addspace}%
	\printfield{doi}%
	\usebibmacro{finentry}%
}

\DeclareBibliographyDriver{online}{%
	\usebibmacro{bibindex}%
	\usebibmacro{begentry}%
	\usebibmacro{author}%
\setunit{\addspace}%
	\usebibmacro{date}%
\newunit\newblock
	\usebibmacro{title}%
\newunit\newblock
	\usebibmacro{eprint}%
	\usebibmacro{finentry}%
}

\DeclareBibliographyDriver{book}{%
	\usebibmacro{bibindex}%
	\usebibmacro{begentry}%
	\usebibmacro{author}%
\setunit{\addspace}\newblock
	\usebibmacro{date}%
\newunit\newblock
	\usebibmacro{title}%
\setunit*{\addspace}\newblock
	\printfield{edition}%
\newunit\newblock
	\printlist{publisher}
\setunit*{\addspace}%
	\printfield{volume}%
\setunit*{\addspace}\newblock%
	\printfield{comments}%
	\usebibmacro{eprint}%
\setunit*{\addspace}
	\printfield{doi}
	\usebibmacro{finentry}%
}

\DeclareBibliographyDriver{thesis}{%
	\usebibmacro{bibindex}%
	\usebibmacro{begentry}%
	\usebibmacro{author}%
\setunit{\addspace}\newblock
	\usebibmacro{date}%
\newunit\newblock
	\usebibmacro{title}.%
\setunit*{\addspace}\newblock
	\space Thesis, \printlist{institution}
	\usebibmacro{eprint}%
	\printfield{comments}%
\setunit*{\addspace}
	\printfield{doi}%
	\usebibmacro{finentry}%
}

\DeclareBibliographyDriver{inproceedings}{%
	\usebibmacro{bibindex}%
	\usebibmacro{begentry}%
	\usebibmacro{author}%
\setunit{\addspace}%
	\usebibmacro{date}%
\newunit\newblock
	\usebibmacro{title}%
\newunit\newblock
	\printfield{booktitle}\addcomma%
\setunit*{\addspace}%
	\printfield{series}\addcomma%
\setunit*{\addspace}\newblock
	\usebibmacro{editor}\addcomma
\setunit*{\addspace}\newblock
	\printlist{publisher}
\setunit*{\addspace}%
	\printfield{volume}%
\setunit*{\addspace}%
	\printfield{pages}
\setunit*{\addspace}
	\usebibmacro{eprint}%
	\printfield{comments}%
\setunit*{\addnnspace}
	\printfield{doi}
	\usebibmacro{finentry}%
}

\DeclareBibliographyDriver{inbook}{%
	\usebibmacro{bibindex}%
	\usebibmacro{begentry}%
	\usebibmacro{author}%
\setunit{\addspace}%
	\usebibmacro{date}%
\newunit\newblock
	\usebibmacro{title}%
\newunit\newblock
	\printfield{booktitle}\addcomma%
\setunit*{\addspace}%
	\printfield{series}\addcomma%
\setunit*{\addspace}\newblock
	\usebibmacro{editor}\addcomma
\setunit*{\addspace}\newblock
	\printlist{publisher}
\setunit*{\addspace}%
	\printfield{volume}%
\setunit*{\addspace}%
	\printfield{pages}
\setunit*{\addspace}
	\usebibmacro{eprint}%
	\printfield{comments}%
\setunit*{\addspace}
	\printfield{doi}
	\usebibmacro{finentry}%
}

\DeclareBibliographyDriver{incollection}{%
	\usebibmacro{bibindex}%
	\usebibmacro{begentry}%
	\usebibmacro{author}%
\setunit{\addspace}%
	\usebibmacro{date}%
\newunit\newblock
	\usebibmacro{title}%
\newunit\newblock
	\printfield{booktitle}\addcomma%
\setunit*{\addspace}%
	\printfield{series}\addcomma%
\setunit*{\addspace}\newblock
	\printlist{publisher}
\setunit*{\addspace}%
	\printfield{volume}%
\setunit*{\addspace}%
	\printfield{pages}
\setunit*{\addspace}
	\usebibmacro{eprint}%
	\printfield{comments}%
\setunit*{\addspace}
	\printfield{doi}
	\usebibmacro{finentry}%
}

\AtEveryBibitem{%
	\clearfield{day}%
	\clearfield{month}%
} 
 
\DeclareNameAlias{sortname}{given-family}
\DeclareNameAlias{default}{given-family}

\newcommand{\acknofootnote}{%
	\blfootnote{%
		\hspace*{1em}%
		Jonathan Hermon%
	\hfill%
		Sam Olesker-Taylor%
		\hspace*{1em}%
	\\%
		\href{mailto:jhermon@math.ubc.ca}{jhermon@math.ubc.ca},
		\href{http://www.math.ubc.ca/~jhermon/}{math.ubc.ca/$\sim$jhermon/}%
	\hfill%
		\href{mailto:sam.ot@posteo.co.uk}{sam.ot@posteo.co.uk},
		\href{https://mathematicalsam.wordpress.com}{mathematicalsam.wordpress.com}%
	\\%
		University of British Columbia, Vancouver, Canada%
	\hfill%
		Department of Mathematical Sciences, University of Bath, UK%
	\\%
		Supported by EPSRC EP/L018896/1 and an NSERC Grant%
	\hfill%
		Supported by EPSRC Grants 1885554 and EP/N004566/1%
	\par\smallskip\par
	\centering%
		The vast majority of this work was undertaken while both authors were at the University of Cambridge%
	}
}

\newcommand*{\MM}{\ensuremath{\Gamma}}

\newcommand{\printtoc}[1]{%
	\ifthenelse%
		{\equal{#1}{1}}%
		{\sffamily\boldmath\tableofcontents\unboldmath\normalfont}%
		{\newpage\small\sffamily\boldmath\tableofcontents\unboldmath\normalfont\normalsize}%
	}

\newcommand{\nextresult}{%
	\setcounter{introthm}{\value{introthm}}
	\setcounter{introcor}{\value{introthm}}
	\setcounter{introconj}{\value{introthm}}
	\setcounter{introdefn}{\value{introthm}}
	\setcounter{intrormkT}{\value{introthm}}
}

\usepackage{indentfirst, environ}
\usepackage[utf8]{inputenc}

\usepackage{empheq}

\usepackage{manyfoot}
\SetFootnoteHook{\hspace*{-1.8em}}
\DeclareNewFootnote{bl}[gobble]
\newcommand{\blfootnote}[1]{\footnotebl{\sffamily#1}}

\usepackage{graphicx}
\graphicspath{{../Figures/}}
\usepackage[labelsep=period]{caption}

\usepackage{chngcntr}
\counterwithin{figure}{section}

\usepackage[UKenglish]{babel}
\usepackage{amsmath, amsthm, amssymb, mathrsfs, bm, mathtools}
\usepackage{wasysym}
\usepackage{mathtools}

\allowdisplaybreaks

\usepackage{cases}

\usepackage{titlesec}
\usepackage{titletoc}

\titleformat{\title}
{\sffamily \Huge \bfseries \scshape \boldmath}{}{}{}
\titleformat{\section}
{\sffamily \Large \bfseries \scshape \boldmath}{\thesection}{1em}{}
\titleformat{\subsection}
{\sffamily \large \bfseries \scshape \boldmath}{\thesubsection}{1em}{}
\titleformat{\subsubsection}
{\sffamily \normalsize \bfseries \scshape \boldmath}{\thesubsubsection}{1em}{}
\titleformat{\paragraph}
{\sffamily \normalsize \bfseries \scshape \boldmath}{\theparagraph}{1em}{}
\titleformat{\subparagraph}[runin]
{\sffamily \normalsize \bfseries \scshape \boldmath}{\thesubparagraph}{1em}{}

\def\IfAmpersandUseAlign#1#2&#3\EndIfAmpersandUseAlign
{%
	\if\relax\detokenize{#3}\relax
	\begin{equation*}%
	#1%
	\end{equation*}%
	\else
	\begin{align*}%
	#1%
	\end{align*}%
	\fi
}
\def\[#1\]%
{%
	\IfAmpersandUseAlign{#1}#1&\EndIfAmpersandUseAlign
}

\usepackage{eqparbox}
\newcommand{\eqmathsbox}[3][\mathrel]{%
	#1{\eqmakebox[#2]{$\displaystyle#3$}}%
}

\usepackage{scrextend}

\usepackage{enumitem}

\newcommand{\numberingroman}{%
	\renewcommand{\labelenumi}{(\roman{enumi})}%
	\renewcommand{\theenumi}{(\roman{enumi})}%
}

\setlist[description]{%
	topsep		= 0pt,		
	noitemsep,				
	font		= {\mdseries\itshape},	
}

\newcommand{\ggr}{G}
\newcommand{\gab}{G^\ab}
\newcommand{\gcom}{G^\com}
\newcommand{\gk}{G_k}

\newcommand{\UU}{U}
\newcommand{\HH}{H}

\newcommand{\ugr}{\UU_{p,d}}
\newcommand{\uab}{\UU_{p,d}^\ab}
\newcommand{\ucom}{\UU_{p,d}^\com}
\newcommand{\uk}{(\UU_{p,d})_k}

\newcommand{\hgr}{\HH_{p,d}}
\newcommand{\hab}{\HH_{p,d}^\ab}

\newcommand{\hk}{(\HH_{p,d})_k}

\usepackage[dvipsnames,hyperref]{xcolor}

\newcommand{\nt}{\addtocounter{equation}{1}\tag{\theequation}}

\newcommand{\bcdot}{\ensuremath{\bm{\cdot}}}
\let\mod\relax
\DeclareMathOperator{\mod}{\, mod}
\DeclareMathOperator{\MOD}{mod}

\let\div\relax
\DeclareMathOperator{\div}{div}

\newcommand{\QUAD}[1]{\quad\text{#1}\quad}
\newcommand{\Quad}[1]{
	\mathchoice
	{\quad\text{#1}\quad}
	{\text{ #1 }}
	{\text{ #1 }}
	{\text{ #1 }}
}

\newcommand{\whp}{\text{whp}\xspace}

\newcommand{\forallZ}{\text{for all $Z$}\xspace}

\newcommand{\Qforall}{\Quad{for all}}
\newcommand{\Qfor}{\Quad{for}}
\newcommand{\Qand}{\Quad{and}}
\newcommand{\Qwhere}{\Quad{where}}
\newcommand{\Qwhen}{\Quad{when}}

\newcommand{\id}{\mathsf{id}}

\newcommand{\typ}{\mathsf{typ}}

\newcommand{\cq}{\coloneqq}
\renewcommand{\epsilon}{\varepsilon}

\newcommand{\eps}{\epsilon}

\usepackage{xifthen}

\newcommand{\binomt}[2]{ \textstyle \binom{#1}{#2} \displaystyle }

\newcommand{\maxt}[1]{ \textstyle \max_{#1} \displaystyle }

\newcommand{\mint}[1]{ \textstyle \min_{#1} \displaystyle }

\usepackage{calc}
\newlength{\halfplusheight}
\newcommand{\MAX}[1]{\mathop{\raisebox{\halfplusheight}{\(\displaystyle\max_{#1}\)}}}

\newcommand{\LIMSUP}[1]{\mathop{\raisebox{\halfplusheight}{\(\displaystyle\limsup_{#1}\)}}}
\newcommand{\LIMINF}[1]{\mathop{\raisebox{\halfplusheight}{\(\displaystyle\liminf_{#1}\)}}}

\DeclareMathOperator*{\sumTT}{\textstyle\sum}
\newcommand{\sumT}[2][]{
	\ifthenelse{\isempty{#1}}
	{\sumTT_{#2}}
	{\sumTT_{#2}^{#1}}
}
\newcommand{\sumt}[2][]{
	\ifthenelse{\isempty{#1}}
	{\textstyle \sum_{#2}      \displaystyle}
	{\textstyle \sum_{#2}^{#1} \displaystyle}
}
\newcommand{\sumd}[2][]{
	\ifthenelse{\isempty{#1}}
	{\displaystyle \sum_{#2}}
	{\displaystyle \sum_{#2}^{#1}}
}

\newcommand{\intt}[2][]{
	\ifthenelse{\isempty{#1}}
	{\textstyle \int_{#2}      \displaystyle}
	{\textstyle \int_{#2}^{#1} \displaystyle}
}

\newcommand{\prodt}[2][]{
	\ifthenelse{\isempty{#1}}
	{\textstyle \prod_{#2}      \displaystyle}
	{\textstyle \prod_{#2}^{#1} \displaystyle}
}
\newcommand{\prodd}[2][]{
	\ifthenelse{\isempty{#1}}
	{\prod_{#2}}
	{\prod_{#2}^{#1}}
}

\usepackage{xspace}

\let\originalexp\exp
\let\exp\relax

\DeclareRobustCommand{\exp} [1]{\originalexp(#1)}
\newcommand{\expb} [1]{\originalexp\bigl( #1 \bigr)}

\newcommand{\abs}  [1]{| #1 |}
\newcommand{\absb} [1]{\big| #1 \bigr|}

\newcommand{\normb} [1]{\big\lVert #1 \bigr\rVert}

\newcommand{\rbr} [1]{ ( #1 ) }
\newcommand{\rbb} [1]{\bigl( #1 \bigr)}
\newcommand{\rbB} [1]{\Bigl( #1 \Bigr)}

\newcommand{\sbb} [1]{\bigl[ #1 \bigr]}

\newcommand{\bra} [1]{ \{ #1 \} }
\newcommand{\brb} [1]{\bigl\{ #1 \bigr\}}

\newcommand{\dist}{\mathrm{dist}}

\DeclareMathOperator{\diam}{diam}

\newcommand{\tdiam}{t_\mathrm{diam}}

\newcommand{\mix}{\mathrm{mix}}

\newcommand{\rel}{\mathrm{rel}}

\newcommand{\st}{{ \ \mathrm{st} \ }}

\newcommand{\ab} {\mathrm{ab}}
\newcommand{\com}{\mathrm{com}}

\DeclareMathOperator{\step}{step}

\newcommand{\Unif}{\mathrm{Unif}}
\newcommand{\iid}{\mathrm{iid}}

\newcommand{\Geom}{\mathrm{Geom}}

\newcommand{\Bern}{\mathrm{Bern}}

\newcommand{\tmix}{t_\mix}

\newcommand{\trel}{t_\rel}

\newcommand{\Ninn}{{N\in\mathbb{N}}}

\newcommand{\cups}{\cup \cdots \cup}

\newcommand{\floor}[1]{\lfloor #1 \rfloor}

\newcommand{\ceil}[1]{\lceil #1 \rceil}
\newcommand{\ceilb}[1]{\bigl\lceil #1 \bigr\rceil}

\newcommand{\midb}{\bigm\vert}

\newcommand{\one}  [1]{\bm1( #1 )}
\newcommand{\oneb} [1]{\bm1\bigl( #1 \bigr)}

\newcommand{\tvb}[1]{\bigl\lVert #1 \bigr\rVert_{\mathrm{TV}}}

\newcommand{\logk}[1][]{
	\ifthenelse{\equal{}{#1}}
	{\log k}
	{(\log k)^{#1}}
}

\newcommand{\logn}[1][]{
	\ifthenelse{\equal{}{#1}}
	{\log n}
	{(\log n)^{#1}}
}

\newcommand{\logm}[1][]{
	\ifthenelse{\equal{}{#1}}
	{\log m}
	{(\log m)^{#1}}
}

\newcommand{\loglogn}[1][]{
	\ifthenelse{\equal{}{#1}}
	{\log\log n}
	{(\log\log n)^{#1}}
}

\newcommand{\prt}[2][]{
	\ifthenelse{\equal{}{#1}}
	{\mathbb{P}(#2)}
	{\mathbb{P}_{#1}(#2)}
}
\newcommand{\pr}[2][]{
	\mathchoice
	{\ifthenelse{\isempty{#1}}
		{\mathbb{P}\bigl(#2\bigr)}
		{\mathbb{P}_{#1}\bigl(#2\bigr)}}
	{\ifthenelse{\isempty{#1}}
		{\mathbb{P}(#2)}
		{\mathbb{P}_{#1}(#2)}}
	{\ifthenelse{\isempty{#1}}
		{\mathbb{P}(#2)}
		{\mathbb{P}_{#1}(#2)}}
	{\ifthenelse{\isempty{#1}}
		{\mathbb{P}(#2)}
		{\mathbb{P}_{#1}(#2)}}
}
\newcommand{\prb}[2][]{
	\ifthenelse{\equal{}{#1}}
	{\mathbb{P}\bigl( #2 \bigr)}
	{\mathbb{P}_{#1}\bigl( #2 \bigr)}
}
\newcommand{\prB}[2][]{
	\ifthenelse{\equal{}{#1}}
	{\mathbb{P}\Bigl( #2 \Bigr)}
	{\mathbb{P}_{#1}\Bigl( #2 \Bigr)}
}
\newcommand{\prbb}[2][]{
	\ifthenelse{\equal{}{#1}}
	{\mathbb{P}\biggl( #2 \biggr)}
	{\mathbb{P}_{#1}\biggl( #2 \biggr)}
}
\newcommand{\prBB}[2][]{
	\ifthenelse{\equal{}{#1}}
	{\mathbb{P}\Biggl( #2 \Biggr)}
	{\mathbb{P}_{#1}\Biggl( #2 \Biggr)}
}
\newcommand{\prs}[2][]{
	\ifthenelse{\equal{}{#1}}
	{\mathbb{P}\left( #2 \right)}
	{\mathbb{P}_{#1}\left( #2 \right)}
}

\newcommand{\qr}[2][]{
	\mathchoice
	{\ifthenelse{\isempty{#1}}
		{\mathbb{Q}\bigl(#2\bigr)}
		{\mathbb{Q}_{#1}\bigl(#2\bigr)}}
	{\ifthenelse{\isempty{#1}}
		{\mathbb{Q}(#2)}
		{\mathbb{Q}_{#1}(#2)}}
	{\ifthenelse{\isempty{#1}}
		{\mathbb{Q}(#2)}
		{\mathbb{Q}_{#1}(#2)}}
	{\ifthenelse{\isempty{#1}}
		{\mathbb{Q}(#2)}
		{\mathbb{Q}_{#1}(#2)}}
}
\newcommand{\qrb}[2][]{
	\ifthenelse{\equal{}{#1}}
	{\mathbb{Q}\bigl( #2 \bigr)}
	{\mathbb{Q}_{#1}\bigl( #2 \bigr)}
}
\newcommand{\qrB}[2][]{
	\ifthenelse{\equal{}{#1}}
	{\mathbb{Q}\Bigl( #2 \Bigr)}
	{\mathbb{Q}_{#1}\Bigl( #2 \Bigr)}
}
\newcommand{\qrbb}[2][]{
	\ifthenelse{\equal{}{#1}}
	{\mathbb{Q}\biggl( #2 \biggr)}
	{\mathbb{Q}_{#1}\biggl( #2 \biggr)}
}
\newcommand{\qrBB}[2][]{
	\ifthenelse{\equal{}{#1}}
	{\mathbb{Q}\Biggl( #2 \Biggr)}
	{\mathbb{Q}_{#1}\Biggl( #2 \Biggr)}
}
\newcommand{\qrs}[2][]{
	\ifthenelse{\equal{}{#1}}
	{\mathbb{Q}\left( #2 \right)}
	{\mathbb{Q}_{#1}\left( #2 \right)}
}

\newcommand{\ext}[2][]{
\ifthenelse{\equal{}{#1}}
{\mathbb{E}(#2)}
{\mathbb{E}_{#1}(#2)}
}
\newcommand{\ex}[2][]{
	\mathchoice
	{\ifthenelse{\isempty{#1}}
		{\mathbb{E}\bigl(#2\bigr)}
		{\mathbb{E}_{#1}\bigl(#2\bigr)}}
	{\ifthenelse{\isempty{#1}}
		{\mathbb{E}(#2)}
		{\mathbb{E}_{#1}(#2)}}
	{\ifthenelse{\isempty{#1}}
		{\mathbb{E}(#2)}
		{\mathbb{E}_{#1}(#2)}}
	{\ifthenelse{\isempty{#1}}
		{\mathbb{E}(#2)}
		{\mathbb{E}_{#1}(#2)}}
}
\newcommand{\exb}[2][]{
	\ifthenelse{\equal{}{#1}}
	{\mathbb{E}\bigl( #2 \bigr)}
	{\mathbb{E}_{#1}\bigr( #2 \bigr)}
}
\newcommand{\exB}[2][]{
	\ifthenelse{\equal{}{#1}}
	{\mathbb{E}\Bigl( #2 \Bigr)}
	{\mathbb{E}_{#1}\Bigl( #2 \Bigr)}
}
\newcommand{\exbb}[2][]{
	\ifthenelse{\equal{}{#1}}
	{\mathbb{E}\biggl( #2 \biggr)}
	{\mathbb{E}_{#1}\biggl( #2 \biggr)}
}
\newcommand{\exBB}[2][]{
	\ifthenelse{\equal{}{#1}}
	{\mathbb{E}\Biggl( #2 \Biggr)}
	{\mathbb{E}_{#1}\Biggl( #2 \Biggr)}
}

\newcommand{\fx}[2][]{
	\ifthenelse{\equal{}{#1}}
	{\mathbb{F}(#2)}
	{\mathbb{F}_{#1}(#2)}
}
\newcommand{\fxb}[2][]{
	\ifthenelse{\equal{}{#1}}
	{\mathbb{F}\bigl( #2 \bigr)}
	{\mathbb{F}_{#1}\bigr( #2 \bigr)}
}
\newcommand{\fxB}[2][]{
	\ifthenelse{\equal{}{#1}}
	{\mathbb{F}\Bigl( #2 \Bigr)}
	{\mathbb{F}_{#1}\Bigl( #2 \Bigr)}
}
\newcommand{\fxbb}[2][]{
	\ifthenelse{\equal{}{#1}}
	{\mathbb{F}\biggl( #2 \biggr)}
	{\mathbb{F}_{#1}\biggl( #2 \biggr)}
}
\newcommand{\fxBB}[2][]{
	\ifthenelse{\equal{}{#1}}
	{\mathbb{F}\Biggl( #2 \Biggr)}
	{\mathbb{F}_{#1}\Biggl( #2 \Biggr)}
}

\newcommand{\Var}[1]{\mathbb{V}\mathrm{ar}(#1)}
\newcommand{\Varb}[2][]{
	\ifthenelse{\equal{}{#1}}
	{\mathbb{V}\mathrm{ar} \bigl(#2\bigr)}
	{\mathbb{V}\mathrm{ar}_{#1} \bigl(#2\bigr)}
}
\newcommand{\VAR}[2][]{
	\ifthenelse{\equal{}{#1}}
	{\mathrm{Var}(#2)}
	{\mathrm{Var}_{#1}(#2)}
}

\newcommand{\Oh}  [1]{\mathcal{O}( #1 )}
\newcommand{\Ohb} [1]{\mathcal{O}\bigl( #1 \bigr)}

\newcommand{\oh}  [1]{o( #1 )}

\newcommand{\mbn}{\mathbb{N}}

\newcommand{\mbp}{\mathbb{P}}

\newcommand{\mbr}{\mathbb{R}}

\newcommand{\mbz}{\mathbb{Z}}

\newcommand{\mca}{\mathcal{A}}
\newcommand{\mcb}{\mathcal{B}}
\newcommand{\mcc}{\mathcal{C}}
\newcommand{\mcd}{\mathcal{D}}
\newcommand{\mce}{\mathcal{E}}

\newcommand{\mci}{\mathcal{I}}
\newcommand{\mcj}{\mathcal{J}}

\newcommand{\mcl}{\mathcal{L}}

\newcommand{\mcw}{\mathcal{W}}

\newcommand{\mfgcd}{\mathfrak{g}}

\newcommand{\toinf}[1]{\ensuremath{#1\to\infty}}
\newcommand{\asinf}[1]{\text{as \ensuremath{#1\to\infty}}}

\newcommand{\tozero}[1]{\ensuremath{#1\to0}}

\newcommand{\Ninf}{{N\to\infty}}

\makeatletter
\newenvironment{subtheorem}[1]{%
	\def\subtheoremcounter{#1}%
	\refstepcounter{#1}%
	\protected@edef\theparentnumber{\csname the#1\endcsname}%
	\setcounter{parentnumber}{\value{#1}}%
	\setcounter{#1}{0}%
	\expandafter\def\csname the#1\endcsname{\theparentnumber\alph{#1}}%
	\expandafter\def\csname theH#1\endcsname{thm.\theparentnumber\alph{#1}}%
	\unskip\ignorespaces
}{%
	\setcounter{\subtheoremcounter}{\value{parentnumber}}%
	\ignorespacesafterend
}
\makeatother
\newcounter{parentnumber}

\makeatletter
\newenvironment{subtheorem-num}[1]{%
	\def\subtheoremcounter{#1}%
	\refstepcounter{#1}%
	\protected@edef\theparentnumber{\csname the#1\endcsname}%
	\setcounter{parentnumber}{\value{#1}}%
	\setcounter{#1}{0}%
	\expandafter\def\csname the#1\endcsname{\theparentnumber.\arabic{#1}}%
	\expandafter\def\csname theH#1\endcsname{thm.\theparentnumber.\arabic{#1}}%
	\unskip\ignorespaces
}{%
	\setcounter{\subtheoremcounter}{\value{parentnumber}}%
	\ignorespacesafterend
}
\makeatother

\usepackage[colorlinks, pdfauthor = {Jonathan\ Hermon\ and\ Sam\ Olesker-Taylor}, pdfusetitle]{hyperref}
\hypersetup{
	colorlinks,
	linkcolor	= {blue},
	filecolor	= {red},
	urlcolor	= {blue},
	citecolor	= {ForestGreen}
}

\newcommand{\qedtriangle}{\renewcommand{\qedsymbol}{\ensuremath{\triangle}}}

\usepackage[capitalise]{cleveref}

\crefformat{equation}{\textup{(#2#1#3)}}
\crefmultiformat{equation}{\textup{(#2#1#3}}{\textup{, #2#1#3)}}{\textup{, #2#1#3}}{\textup{, #2#1#3)}}
\crefrangeformat{equation}{\textup{(#3#1#4--#5#2#6)}}

\newenvironment{Proof}[1][\proofname]{%
	\proof[\upshape\bfseries\sffamily\boldmath#1]
}{\endproof}

\usepackage{etoolbox}
\usepackage{needspace}

\newtheoremstyle{sfsl}
{1\baselineskip}		
{1\baselineskip}		
{\slshape}				
{}						
{\bfseries\sffamily}	
{.}						
{0.5em}					
{\thmname{#1}\thmnumber{ #2}\thmnote{ {\mdseries(#3)}}}

\newtheoremstyle{sfup}
{1\baselineskip}		
{1\baselineskip}		
{\upshape}				
{}						
{\bfseries\sffamily}	
{.}						
{0.5em}					
{\thmname{#1}\thmnumber{ #2}\thmnote{ {\mdseries(#3)}}}

\theoremstyle{sfsl}

\newtheorem*{thm*}{Theorem}
\newtheorem{thm} {Theorem}[section]
\crefname{thm}{Theorem}{Theorems}

\newtheorem*{introthm*}{Theorem}
\newtheorem{introthm}{Theorem}

\crefname{introthm}{Theorem}{Theorems}

\newtheorem*{cor*}{Corollary}

\crefname{cor}{Corollary}{Corollaries}

\newtheorem*{introcor*}{Corollary}

\crefname{introcor}{Corollary}{Corollaries}

\newtheorem*{introconj*}{Conjecture}

\crefname{introconj}{Conjecture}{Conjectures}

\newtheorem*{introques*}{Question}

\crefname{introques}{Question}{Questions}

\newtheorem*{lem*}    {Lemma}
\newtheorem{lem} [thm]{Lemma}
\crefname{lem}{Lemma}{Lemmas}

\newtheorem*{introlem*}{Lemma}

\crefname{introlem}{Lemma}{Lemmas}

\newtheorem*{prop*}    {Proposition}
\newtheorem{prop} [thm]{Proposition}
\crefname{prop}{Proposition}{Propositions}

\newtheorem*{clm*}    {Claim}
\newtheorem{clm} [thm]{Claim}
\crefname{clm}{Claim}{Claims}

\newtheorem*{defn*}    {Definition}
\newtheorem{defn} [thm]{Definition}
\crefname{defn}{Definition}{Definitions}

\newtheorem*{introdefn*}{Definition}
\newtheorem{introdefn}{Definition}

\crefname{introdefn}{Definition}{Definitions}

\newtheorem{innercustomgenericsl}{\customgenericnamesl}
\providecommand{\customgenericnamesl}{}
\newcommand{\newcustomtheoremsl}[2]{%
	\newenvironment{#1}[1]
	{%
		\renewcommand\customgenericnamesl{#2}%
		\renewcommand\theinnercustomgenericsl{##1}%
		\innercustomgenericsl
	}
	{\endinnercustomgenericsl}
}

\newcustomtheoremsl{customthm}{Theorem}
\newcustomtheoremsl{customprop}{Proposition}
\newcustomtheoremsl{customlemma}{Lemma}
\newcustomtheoremsl{customrmk}{Remark}
\newcustomtheoremsl{customconj}{Conjecture}
\newcustomtheoremsl{customdefn}{Definition}
\newcustomtheoremsl{customquestion}{Question}
\newcustomtheoremsl{customopenproblem}{Open Problem}
\newcustomtheoremsl{customhyp}{Hypothesis}

\newtheorem*{conj*}   {Conjecture}

\crefname{conj}{Conjecture}{Conjectures}

\newenvironment{conj-ind*}
	{\begin{quote}\textsf{\textbf{Conjecture.}}\slshape}
	{\end{quote}}
\newenvironment{conj-ind}
	{\begin{quote}\vspace{-\glueexpr\baselineskip+\topsep}\begin{customconj}}
	{\end{customconj}\end{quote}}

\newenvironment{question-ind*}
	{\begin{quote}\textsf{\textbf{Question.}}\slshape}
	{\end{quote}}
\newenvironment{question-ind}
	{\begin{quote}\vspace{-\glueexpr\baselineskip+\topsep}\begin{customquestion}}
	{\end{customquestion}\end{quote}}

\newenvironment{openproblem-ind*}
	{\begin{quote}\textsf{\textbf{Open Problem.}}\slshape}
	{\end{quote}}
\newenvironment{openproblem-ind}
	{\begin{quote}\vspace{-\glueexpr\baselineskip+\topsep}\begin{customopenproblem}}
	{\end{customopenproblem}\end{quote}}

\newtheorem*{hypothesis*}{Hypothesis}

\newtheorem*{hyp*}{Hypothesis}
\newtheorem{hyp}{Hypothesis}

\crefname{hyp}{Hypothesis}{Hypotheses}

\newtheorem*{rmk*}{Remark}

\theoremstyle{sfup}

\newtheorem{innercustomgenericup}{\customgenericnameup}
\providecommand{\customgenericnameup}{}
\newcommand{\newcustomtheoremup}[2]{%
	\newenvironment{#1}[1]
	{%
		\renewcommand\customgenericnameup{#2}%
		\renewcommand\theinnercustomgenericup{##1}%
		\innercustomgenericup
	}
	{\endinnercustomgenericup}
}

\crefname{exm} {Example}{Examples}
\crefname{exmT}{Example}{Examples}

\newtheorem{rmkT}[thm]{Remark}
	\newenvironment{rmkt}
	{\pushQED{\qed}\renewcommand{\qedsymbol}{\ensuremath{\triangle}}\rmkT}
	{\popQED\endrmkT}
\crefname{rmk} {Remark}{Remarks}
\crefname{rmkT}{Remark}{Remarks}

\newtheorem*{rmkTT}{Remark}
\newenvironment{rmkt*}
	{\pushQED{\qed}\renewcommand{\qedsymbol}{\ensuremath{\triangle}}\rmkTT}
	{\popQED\endrmkTT}

\newcustomtheoremup{customrmkT}{Remark}
\newenvironment{customrmkt}
	{\pushQED{\qed}\renewcommand{\qedsymbol}{\ensuremath{\triangle}}\customrmkT}
	{\popQED\endcustomrmkT}

\crefname{rmks} {Remarks}{Remarks}
\crefname{rmksT}{Remarks}{Remarks}

\newtheorem*{rmks*} {Remarks}
\newtheorem*{rmksTT}{Remarks}
\newenvironment{rmkst*}
	{\pushQED{\qed}\renewcommand{\qedsymbol}{\ensuremath{\triangle}}\rmksTT}
	{\popQED\endrmksTT}

\newtheorem{intrormkT}{Remark}
	\newenvironment{intrormkt}
	{\pushQED{\qed}\renewcommand{\qedsymbol}{\ensuremath{\triangle}}\intrormkT}
	{\popQED\endintrormkT}

\crefname{intrormk} {Remark}{Remarks}
\crefname{intrormkT}{Remark}{Remarks}

\newtheorem*{intrormk*} {Remark}
\newtheorem*{intrormkTT}{Remark}
\newenvironment{intrormkt*}
	{\pushQED{\qed}\renewcommand{\qedsymbol}{\ensuremath{\triangle}}\intrormkTT}
	{\popQED\endintrormkTT}

\newtheorem*{exm*} {Example}
\newtheorem*{exmTT}{Lemma}
	\newenvironment{exmt*}
	{\pushQED{\qed}\renewcommand{\qedsymbol}{\ensuremath{\triangle}}\exmTT}
	{\popQED\endexmTT}

\newtheorem*{note*} {Note}
\newtheorem*{noteTT}{Note}
	\newenvironment{notet*}
	{\pushQED{\qed}\renewcommand{\qedsymbol}{\ensuremath{\triangle}}\noteTT}
	{\popQED\endnoteTT}

\catcode`\^^A=14

\usepackage{xparse}

\newcounter{mixedsubequations}

\ExplSyntaxOn

\int_new:N \g_mixedsubeq_int

\NewDocumentEnvironment{mixedsubequations}{o}
{
	\IfNoValueTF { #1 }
	{
		\addtocounter{equation}{-\g_mixedsubeq_int}
		\stepcounter{mixedsubequations}
	}
	{
		\int_gset:Nn \g_mixedsubeq_int { \clist_count:n { #1 } }
		\clist_map_inline:nn { #1 }
		{
			\refstepcounter{equation}\label{##1}
		}
		\addtocounter{equation}{-\g_mixedsubeq_int}
		\setcounter{mixedsubequations}{1}
	}
	\domixedsubequations
}
{\ignorespacesafterend}

\NewDocumentCommand{\domixedsubequations}{}
{
	\cs_set:Npx \theequation
	{
		\exp_not:o { \theequation }
		\exp_not:n { \alph{mixedsubequations} }
	}
	\ignorespaces
}
\ExplSyntaxOff

\makeatletter
\let\save@mathaccent\mathaccent
\newcommand*\if@single[3]{%
  \setbox0\hbox{${\mathaccent"0362{#1}}^H$}%
  \setbox2\hbox{${\mathaccent"0362{\kern0pt#1}}^H$}%
  \ifdim\ht0=\ht2 #3\else #2\fi
  }
\newcommand*\rel@kern[1]{\kern#1\dimexpr\macc@kerna}
\newcommand*\widebar[1]{\@ifnextchar^{{\wide@bar{#1}{0}}}{\wide@bar{#1}{1}}}
\newcommand*\wide@bar[2]{\if@single{#1}{\wide@bar@{#1}{#2}{1}}{\wide@bar@{#1}{#2}{2}}}
\newcommand*\wide@bar@[3]{%
  \begingroup
  \def\mathaccent##1##2{%
    \let\mathaccent\save@mathaccent
    \if#32 \let\macc@nucleus\first@char \fi
    \setbox\z@\hbox{$\macc@style{\macc@nucleus}_{}$}%
    \setbox\tw@\hbox{$\macc@style{\macc@nucleus}{}_{}$}%
    \dimen@\wd\tw@
    \advance\dimen@-\wd\z@
    \divide\dimen@ 3
    \@tempdima\wd\tw@
    \advance\@tempdima-\scriptspace
    \divide\@tempdima 10
    \advance\dimen@-\@tempdima
    \ifdim\dimen@>\z@ \dimen@0pt\fi
    \rel@kern{0.6}\kern-\dimen@
    \if#31
      \overline{\rel@kern{-0.6}\kern\dimen@\macc@nucleus\rel@kern{0.4}\kern\dimen@}%
      \advance\dimen@0.4\dimexpr\macc@kerna
      \let\final@kern#2%
      \ifdim\dimen@<\z@ \let\final@kern1\fi
      \if\final@kern1 \kern-\dimen@\fi
    \else
      \overline{\rel@kern{-0.6}\kern\dimen@#1}%
    \fi
  }%
  \macc@depth\@ne
  \let\math@bgroup\@empty \let\math@egroup\macc@set@skewchar
  \mathsurround\z@ \frozen@everymath{\mathgroup\macc@group\relax}%
  \macc@set@skewchar\relax
  \let\mathaccentV\macc@nested@a
  \if#31
    \macc@nested@a\relax111{#1}%
  \else
    \def\gobble@till@marker##1\endmarker{}%
    \futurelet\first@char\gobble@till@marker#1\endmarker
    \ifcat\noexpand\first@char A\else
      \def\first@char{}%
    \fi
    \macc@nested@a\relax111{\first@char}%
  \fi
  \endgroup
}
\makeatother

\numberwithin{equation}{section}

\title{\sffamily Cutoff for Random Walks on Upper Triangular Matrices}

\author{\sffamily Jonathan Hermon\quad Sam Olesker-Taylor}

\date{}

\usepackage{geometry}
\begin{document}

%

\maketitle

\acknofootnote

\vspace{-6ex}

\renewcommand{\abstractname}{\sffamily Abstract}
\begin{abstract}
Consider the random Cayley graph of a finite group $G$ with respect to $k$ generators chosen uniformly at random, with $1 \ll \log k \ll \log |G|$ (ie $1 \ll k = |G|^{o(1)}$). A conjecture of Aldous and Diaconis~\cite{AD:conjecture} asserts, for $k \gg \log |G|$, that the random walk on this graph exhibits cutoff.

When $\log k \lesssim \log \log |G|$ (ie $k = (\log |G|)^{\mathcal O(1)}$), the only example of a non-Abelian group for which cutoff has been established is the dihedral group. We establish cutoff (as $p \to \infty$) for the group of $d \times d$ unit upper triangular matrices with integer entries modulo $p$ (prime), which we denote $U_{p,d}$, for fixed $d$ or $d$ diverging sufficiently slowly. We allow $1 \ll k \lesssim \log |U_{p,d}|$ as well as $k \gg \log |U_{p,d}|$. The cutoff time is $\max\{\log_k |U_{p,d}|, \: s_0 k\}$, where $s_0$ is the time at which the entropy of the random walk on $\mathbb Z$ reaches $(\log |U_{p,d}^\mathrm{ab}|)/k$, where $U_{p,d}^\mathrm{ab} \cong \mathbb Z_p^{d-1}$ is the Abelianisation of $U_{p,d}$. When $1 \ll k \ll \log |U_{p,d}^\mathrm{ab}|$ and $d \asymp 1$, we find the limit profile.
We also prove highly related results for the $d$-dimensional Heisenberg group over $\mathbb Z_p$.

The Aldous--Diaconis conjecture also asserts, for $k \gg \log |G|$, that the cutoff time should depend only on $k$ and $|G|$. This was verified for all Abelian groups. Our result shows that this is not the case for $U_{p,d}$: the cutoff time depends on $k$, $|U_{p,d}| = p^{d(d-1)/2}$ and $|U_{p,d}^\mathrm{ab}| = p^{d-1}$.

We also show that all but $o(|U_{p,d}|)$ of the elements of $U_{p,d}$ lie at graph distance $M \pm o(M)$ from the identity, where $M$ is the minimal radius of a ball in $\mathbb Z^k$ of cardinality $|U_{p,d}^\mathrm{ab}| = p^{d-1}$. Finally, we show that the diameter is also asymptotically $M$ when $k \gtrsim \log |U_{p,d}^\textrm{ab}|$ and $d \asymp 1$.

The aforementioned results all hold with high probability over the random Cayley graph.
\end{abstract}

\small
\begin{quote}
\begin{description}
	\item [Keywords:]
	cutoff, mixing times, random walk, random Cayley graphs, upper triangular ma\-trices, Heisenberg group, concentration of measure, entropy, typical distance, diameter
	
	\item [MSC 2020 subject classifications:]
	05C12, 05C80, 05C81; 20D15; 60B15, 60C05, 60J27, 60K37
\end{description}
\end{quote}
\normalsize


\vspace{2ex}
\numberingroman

\vfill
\printtoc{\value{tocdepth}}
\vspace*{2ex}

\newpage
\section{Introduction and Statement of Results}
\label{sec-p1:intro}

\subsection{Motivation, Brief Overview of Results and Notation}

\subsubsection{Motivating Conjecture of \citeauthor{AD:conjecture} and Objectives of Paper}

We analyse properties of the random walk (abbreviated \textit{RW}) on a \textit{Cayley graph} of a finite group.
The generators of this graph are chosen independently and uniformly at random.
Precise definitions are given in \S\ref{sec-p2:intro:cayley-def}; for now, let $G$ be a finite group, let $k$ be an integer (allowed to depend on $G$) and denote by $G_k$ the Cayley graph of $G$ with respect to $k$ independently and uniformly random generators.
We consider values of $k$ with $1 \ll \log k \ll \log \abs G$ for which $G_k$ is connected with high probability (abbreviated \textit{whp}), ie with probability tending to 1 as $\abs G$ grows.


\medskip

\textcite{AD:conjecture,AD:shuff-stop} coin the phrase \textit{cutoff phenomenon}:
	this occurs when the total variation distance (TV) between the law of the RW and its invariant distribution drops abruptly from close to $1$ to close to $0$ in a time-interval of smaller order than the mixing time.
The material in this article is motivated by a conjecture of theirs regarding `universality of cutoff' for the RW on the random Cayley graph $G_k$; it is given in \cite[Page~40]{AD:conjecture}, which is an extended version of \cite{AD:shuff-stop}.

\begin{introconj*}[\citeauthor{AD:conjecture}, \citeyear{AD:conjecture}]
	For any group $G$,
	if $k \gg \log \abs G$ and $\log k \ll \log \abs G$, then
	the random walk on $G_k$ exhibits cutoff whp.
	Further, the cutoff time, to leading order, is independent of the algebraic structure of the group: it can be written as a function only of~$k$~and~$\abs G$.
\end{introconj*}

This conjecture spawned a large body of work, including \cite{D:phd,DH:enumeration-rws,H:cutoff-cayley->,H:cutoff-cayley-survey,H:cutoff-cayley-<,R:random-random-walks,W:rws-hypercube}; see \S\ref{sec-p1:intro:previous-work}.
It has been established in the Abelian set-up by \citeauthor{DH:enumeration-rws} \cite{DH:enumeration-rws,H:cutoff-cayley->}; see \S\ref{sec-p1:intro:previous-work:ad-conj} and \cite[\S\ref{sec-p2:conc-rmks:roichman}]{HOt:rcg:abe:cutoff}.
Further, their upper bound is valid for all groups.
When $\log k \gg \log \log \abs G$, this upper bound matches the trivial diameter lower bound of $\log_k \abs G$.
Herein we focus on $1 \ll \log k \lesssim \log \log \abs G$.

\medskip

The purpose of this paper is threefold, all motivated by the Aldous--Diaconis conjecture:
\begin{enumerate}[noitemsep, topsep = \smallskipamount, label = (\roman*), ref = (\roman*)]
	\item \label{itm-p1:intro:purpose:two}
	to study the RW on $G_k$ for two canonical families of nilpotent groups;
	
	\item \label{itm-p1:intro:purpose:gen}
	to develop general techniques which could be applied to general nilpotent groups;
	
	\item \label{itm-p1:intro:purpose:faster}
	to demonstrate that the RW on the random Cayley graph for a non-Abelian group can mix significantly faster than for an Abelian group, even when $k \gg \log \abs G$.
\end{enumerate}
We elaborate briefly on each of these.

\smallskip

\ref{itm-p1:intro:purpose:two}
While much work has been done in the Abelian set-up, little is known when $G$ is non-Abelian.
We prove cutoff for the RW on $G_k$ for two canonical families of nilpotent groups:
	$\UU_{m,d}$, the $d \times d$ unit-upper triangular matrices with entries in $\mbz_m$;
	$\HH_{m,d}$, the $d$-dimensional Heisenberg group over $\mbz_m$.
Further, for the former we study the geometry of $G_k$:
	we establish concentration for the typical distance of a vertex from the origin and determine its asymptotic growth,
	as well as determining the asymptotic growth of the diameter in a certain regime.

\smallskip

\ref{itm-p1:intro:purpose:gen}
We develop a general techniques in our study of cutoff for $\UU_{m,d}$.
The analysis of $\HH_{m,d}$ requires only a small adaptation to that of $\UU_{m,d}$.
But further, the techniques apply even to typical distance and diameter.
We believe our techniques are widely applicable to nilpotent groups.
We briefly elaborate in the concluding remarks.
Further work (in preparation) explores this in detail.

\smallskip
\ref{itm-p1:intro:purpose:faster}
Our results demonstrate that the cutoff time in the non-Abelian set-up can be significantly smaller than in the Abelian and depend on more than just $k$ and $\abs G$.
Precisely,
the mixing time for the RW on $G_k$ is a factor $d$ smaller when $G = \UU_{m,d}$ than when $G$ is an Abelian group of size $\abs{\UU_{m,d}}$ in the regime $\log \abs{\UU_{m,d}} \ll k \le (\log \abs{\UU_{m,d}})^{1 + 2/(d-2)}$, provided $d$ does not diverge too quickly; see \cref{res-p1:intro:cutoff}.
We thus reject the secondary aspect of the Aldous--Diaconis conjecture.

Previous work establishing cutoff in the non-Abelian set-up was limited to $k \gg \log \abs G$ and to `essentially Abelian' groups, in two somewhat complementary senses.
	First,
	\textcite{DH:enumeration-rws} consider certain groups that have an Abelian subgroup of bounded index---in fact, the only concrete example they exhibit is the dihedral group.
	Second,
	cutoff can be derived for general groups with a large Abelian quotient, namely $\log \abs{G / [G,G]} \eqsim \log \abs \ggr$, from the results in \cite{DH:enumeration-rws}.
We establish cutoff in the regime $1 \ll k \lesssim \log \abs G$ for certain groups with a large Abelian quotient in a companion paper \cite{HOt:rcg:abe:cutoff}.
The cutoff time in all the above examples agrees with that of the Abelian case.
Contrastingly, the cutoff time for $\UU_{m,d}$ in the current article has markedly different behaviour.



\subsubsection{Brief Overview of Results}

\renewcommand{\ugr}{\UU}
\renewcommand{\uab}{\UU^\ab}
\renewcommand{\ucom}{\UU_^\com}
\renewcommand{\uk}{\UU_k}

Let $\ugr \cq \UU_{p,d}$ be the group of $d \times d$ unit upper triangular matrices with entries in $\mbz_p$.
The expression $\UU_k$ denotes the random Cayley graph $(\UU_{p,d})_k$, ie $\uk \cq (\UU_{p,d})_k$.
For a group $G$, denote by $\gcom \cq [G, G]$ its \textit{commutator} and by $\gab \cq G/\gcom$ its \textit{Abelianisation}.
We have $\UU_{p,d}^\ab \cong \mbz_p^{d-1}$.

Our focus is on mixing properties of the RW on the random Cayley graph $\UU_k$.
We consider the limit as $n \cq \abs \ugr = p^{d(d-1)/2} \to \infty$ under the assumption that $1 \ll \log k \ll \log \abs \ugr$.
The condition $1 \ll \log k \ll \log \abs G$ is necessary for cutoff on $G_k^\pm$ for all nilpotent $G$; see \cref{rmk-p1:intro:cutoff:nec-con}.

\begin{itemize}[itemsep = 0pt, topsep = \smallskipamount, label = \bcdot]
	\item 
	\textit{Cutoff.}
	For the RW on $\ugr_k$ with $p$ prime and $d \ge 3$ (allowed to diverge suitably),
	we establish cutoff.
	Mixing on $\ugr$ is governed by mixing on the Abelianisation $\uab$ in a quantitative sense.
	Eg,
	for $1 \ll k \ll \log \abs \uab$,
	the mixing times for $\ugr_k$ and $\uab_k \cq (\UU_{p,d}^\ab)_k$ are asymptotically equivalent;
	in this case when $d \asymp 1$ we find the limit profile of the convergence to equilibrium.
	Under suitable conditions, we relax the requirement that $p$ is prime.
	
	We also prove highly related results for the $d$-dimensional Heisenberg group over $\mbz_p$.
	
	\item 
	\textit{Typical Distance and Diameter.}
	Draw $V \sim \Unif(\ugr)$.
	We show that the law of the graph distance $\dist(\id,V)$ concentrates.
	This value is again governed by the Abelianisation: typical distances on $\ugr_k$ and $\uab_k$ are asymptotically equivalent
	for $1 \ll k \ll \log \abs \uab$.
	For $k \gtrsim \log \abs \uab$ and $d \asymp 1$, we show that the diameter and typical distance asymptotically~agree.
\end{itemize}

Introduced by \textcite{AD:conjecture}, there has been a great deal of research into these random Cayley graphs.
Motivation for this model and an overview of historical work is given~in~\S\ref{sec-p1:intro:previous-work}.

\subsubsection{Notation and Terminology}
\label{sec-p1:intro:not-term}

Let $m, d \in \mbn$ with $d \ge 3$.
We denote by $\HH_{m,d}$ the $d$-dimensional Heisenberg group, which can be represented as the set of triples $(x,y,z) \in \mbz_m^{d-2} \times \mbz_m^{d-2} \times \mbz$ with multiplication defined by
\[
	(x, y, z) \circ (x', y', z')
=
	(x + x', y + y', z + z' + x \cdot y').
\]
An element $(x,y,z) \in \HH_{m,d}$ can be represented as a $d \times d$ unit upper triangular matrix:
place
	$x_1, ..., x_{d-2}$ on the top row,
	$y_1, ..., y_{d-2}$ on the right side,
	$z$ on the top-right corner
and
	$0$s elsewhere.
Via this representation, one immediately sees that $\HH_{m,d}$ can be viewed as a subgroup of $\UU_{m,d}$.

This article is devoted primarily to $\UU_{m,d}$, which we term the \textit{upper triangular group}.

\smallskip

Cayley graphs are either directed or undirected; we emphasise this by writing $G_k^+$ and $G_k^-$, respectively.
When we write $G_k$ or $G^\pm_k$, this means ``either $G^-_k$ or $G^+_k$'', corresponding to the undirected, respectively directed, graphs with generators chosen independently and uniformly at~random.

Conditional on being simple, $G^+_k$ is uniformly distributed over the set of all simple degree-$k$ Cayley graphs. Up to a slightly adjusted definition of \textit{simple} for undirected Cayley graphs, our results hold with $G_k$ replaced by a uniformly chosen simple Cayley graph of degree $k$; see \S\ref{sec-p1:intro:rmks:typ-simp}.

Our results are for sequences $(G_N)_\Ninn$ of finite groups with \toinf{\abs{G_N}} \asinf N.
For ease of presentation, we write statements like ``let $G$ be a group'' instead of ``let $(G_N)_\Ninn$ be a sequence of groups''.
Likewise, the quantities $d$, $m$, $p$ and, of course, $k$ appearing in the statements all correspond to sequences, which need not be fixed (or bounded) unless we explicitly say otherwise.
In the same vein, an event holds \textit{with high probability} (abbreviated \textit{\whp}) if its probability tends~to~1.

We use standard asymptotic notation:
	``$\ll$'' or ``$\oh{\cdot}$'' means ``of smaller order'';
	``$\lesssim$'' or $\Oh{\cdot}$'' means ``of order at most'';
	``$\asymp$'' means ``of the same order'';
	``$\eqsim$'' means ``asymptotically equivalent''.

\subsection{Statements of Main Results}
\label{sec-p1:intro:res}

\renewcommand{\ugr}{\UU_{p,d}}
\renewcommand{\uab}{\UU_{p,d}^\ab}
\renewcommand{\ucom}{\UU_{p,d}^\com}
\renewcommand{\uk}{(\UU_{p,d})_k}

We analyse mixing in the \textit{total variation} (abbreviated \textit{TV}) distance.
The uniform distribution on $G$, denoted $\pi_G$, is invariant for the RW.
Let $S = (S(t))_{t\ge0}$ denote the RW on $G_k$; its law is denoted $\pr[G_k]{S(t) \in \cdot}$.
For $t \ge 0$,
denote the TV distance between the law of $S(t)$ and $\pi_G$ by
\[
	d_{G_k}(t)
\cq
	\tvb{ \pr[G_k]{ S(t) \in \cdot } - \pi_G }
=
	\MAX{A \subseteq G}
	\absb{ \pr[G_k]{ S(t) \in A } - \abs A / \abs G }.
\]
Throughout, unless explicitly specified otherwise, we use continuous time: $t \ge 0$ means $t \in [0,\infty)$.

\subsubsection{Cutoff for Random Walk on Random Cayley Graph of the Upper Triangular Group}

\nextresult

We use standard notation and definitions for \textit{mixing} and \textit{cutoff}; see, eg, \cite[\S 4 and \S 18]{LPW:markov-mixing}.

\begin{introdefn*}
	A sequence $(X^N)_\Ninn$ of Markov chains is said to exhibit \textit{cutoff} when, in a short time-interval, known as the \textit{cutoff window}, the TV distance of the distribution of the chain from equilibrium drops from close to $1$ to close to $0$, or more precisely if there exists $(t_N)_\Ninn$ with
	\[
		\LIMINF{\Ninf} d_N\rbb{ t_N (1 - \eps) } = 1
	\Qand
		\LIMSUP{\Ninf} d_N\rbb{ t_N (1 + \eps) } = 0
	\Qforall
		\eps \in (0,1),
	\]
	where $d_N(\cdot)$ is the TV distance of $X^N(\cdot)$ from its equilibrium distribution for each $\Ninn$.
	
	We say that a RW on a sequence of random graphs $(H_N)_\Ninn$ \textit{exhibits cutoff around time $(t_N)_\Ninn$ whp} if,
	for all fixed $\eps$, in the limit $\Ninf$,
	the TV distance at time $(1 + \eps) t_N$ converges in distribution to $0$ and at time $(1 - \eps) t_N$ to $1$,
	where the randomness is over $H_N$.
\end{introdefn*}

The main contribution of this article is to establish cutoff whp for a family of non-Abelian groups in this set-up of random walks on random Cayley graphs.
As a by-product, we reject the secondary aspect of Aldous--Diaconis conjecture.

\smallskip

We use an \emph{entropic method}, which involves defining \emph{entropic times}; see \S\ref{sec-p1:intro:previous-work:generic-ent} for a high-level description of the method and \S\ref{sec-p1:ent:method} for the specific application.
The main idea is to use an auxiliary process $W$ to generate the walk $S$; one then studies the entropy of the process $W$.
Write $Z = [Z_1, ..., Z_k]$ for the (multiset of) generators of the Cayley graph; then $G_k$ corresponds to choosing $Z_1, ..., Z_k \sim^\iid \Unif(G)$.
	Here, $W_i(t)$ is, for each $i$, the number of times generator $Z_i$ has been applied minus the number of times $Z_i^{-1}$ has been applied;
	$W$ is a rate-1 RW on $\mbz^k$.

For undirected graphs, $W$ is the usual simple RW (abbreviated \textit{SRW}):
	a coordinate is selected uniformly at random and incremented/decremented by 1 each with probability $\tfrac12$.
For directed graphs, inverses are never applied, so a step of $W$ is as follows:
	a coordinate is selected uniformly at random and incremented by 1;
	we term this the \textit{directed random walk} (abbreviated \textit{DRW}).

\begin{introdefn}
\label{def-p1:intro:t*}
	Let $t^\pm_0(k, N)$ be the time at which the entropy of rate-1 RW (ie SRW or DRW, as appropriate) on $\mbz^k$ is $\log N$.
	Define
	\(
		t^\pm_*(k,p,d) \cq \max\bra{ t^\pm_0(k, \abs \uab), \: \log_k \abs \ugr }.
	\)
\end{introdefn}

A description of $t^\pm_0(k, \abs \uab)$, up to subleading order terms, can be found in \cref{res-p1:ent:t0}.
Recall that $1 \ll \log k \ll \log \abs G$ is necessary for cutoff for nilpotent $G$, eg $\UU$ or $\HH$;
see \cref{rmk-p1:intro:cutoff:nec-con}.
A more refined statement than the one below,
with quantified conditions on $d$,
is given in \cref{res-p1:cutoff:res}. 

\renewcommand{\ugr}{\UU}
\renewcommand{\uab}{\UU^\ab}
\renewcommand{\uk}{\UU_k}

\begin{introthm}
\label{res-p1:intro:cutoff}
	Let $p$ be prime and $d \ge 3$.
	Abbreviate $\ugr \cq \UU_{p,d}$.
	Suppose that $1 \ll \log k \ll \log \abs \ugr$ and that $d$ is fixed or diverges sufficiently slowly.
	Then the RW on $\uk^\pm$ exhibits cutoff whp at
	\[
		t^\pm_*(k,p,d)
	=
		\max\brb{t^\pm_0(k, \abs \uab), \: \log_k \abs \ugr}
	\eqsim
	\begin{cases}
		t^\pm_0(k, \abs \uab)
			&\text{when}\quad
		k \le \rbr{ \log \abs \uab }^{1 + 2/(d-2)},
	\\
		\log_k \abs \ugr
			&\text{when}\quad
		k \ge \rbr{ \log \abs \uab }^{1 + 2/(d-2)}.
	\end{cases}
	\]
\end{introthm}

\renewcommand{\ugr}{\UU_{p,d}}
\renewcommand{\uab}{\UU_{p,d}^\ab}
\renewcommand{\uk}{(\UU_{p,d})_k}

	The cutoff time for $\ugr$ cannot, for all regimes of $k$, be written as a function only of $k$ and $\abs \ugr = p^{d(d-1)/2}$.
	The only additional information required is the size of the Abelianisation, ie $\abs \uab = p^{d-1}$.
	In Open Question~\ref{oq-p1:cutoff:general-groups} we discuss potential generalisations of this phenomenon.

\bigskip

\nextresult

We adapt the proof of \cref{res-p1:intro:cutoff} to prove some related results.
The proofs are given in \S\ref{sec-p1:ext}.

\begin{description}[itemsep = 0pt, topsep = \smallskipamount, font = {\mdseries}]
	\item [\cref{res-p1:intro:ext:comp}:]
	we remove the condition that $p$ is prime under suitable conditions.
	
	\item [\cref{res-p1:intro:ext:heis}:]
	we replace the upper triangular group $\ugr$ with the Heisenberg group $\hgr$.
	
	\item [\cref{res-p1:intro:ext:profile}:]
	we find the limit profile of the cutoff when $1 \ll k \ll \log \abs \uab$.
\end{description}



Write $\div m \cq \sumt[m]{\ell=1} \one{ \ell \ \text{divides} \ m }$ for the number of divisors of $m \in \mbn$.
It is standard that `typically' $\div m \asymp \log m$; see \cite[\S 18.2]{HW:number-theory}.
Thus
\(
	\log \div m
\le
	\sqrt{\log m}
\)
for $m$ in a density-$1$ subset of~$\mbn$.

We use a subscript $m$, eg in $\UU_{m,d}$, rather than $p$ to emphasise that $m \in \mbn$ need not be prime.

\renewcommand{\ugr}{\UU}
\renewcommand{\uab}{\UU^\ab}
\renewcommand{\uk}{\UU_k}

\begin{introthm}
\label{res-p1:intro:ext:comp}
	Let $\ugr \cq \UU_{m,d}$.
	Suppose that
		$1 \ll \log k \ll \log \abs \ugr$
	and
		\(
			k \gg \max\{ \sqrt{\log m}, \: \log \div m \};
		\)
		in particular, $\log \div m \ll \log m$, so $k \gtrsim \log \abs \uab$ is sufficient for the latter.
	Suppose that $d \ge 3$ is fixed or diverges sufficiently slowly.
	Then the RW on $\uk^\pm$ exhibits cutoff whp~at~$t_*^\pm(k, m, d)$.
\end{introthm}

\renewcommand{\ugr}{\UU_{m,d}}
\renewcommand{\uab}{\UU_{m,d}^\ab}
\renewcommand{\uk}{(\UU_{m,d})_k}



We discuss next the Heisenberg group $\hgr$.
Recall that $\abs{\HH_{m,d}} = m^{2d-3}$ and $\abs{\HH_{m,d}^\ab} = m^{2d-4}$.

\renewcommand{\hgr}{\HH}
\renewcommand{\hab}{\HH^\ab}
\renewcommand{\hk}{\HH_k}

\begin{introthm}
\label{res-p1:intro:ext:heis}
	Let $\hgr \cq \HH_{m,d}$.
	Suppose	that $1 \ll \log k \ll \log \abs \hgr$.
	If $m$ is not prime, then require
	\(
		k \gg \max\bra{ \sqrt{\log m}, \: \log \div m };
	\)
	in particular, $\log \div m \ll \log m$, so $k \gtrsim \log \abs \hab$ suffices.
	Suppose that $d \ge 3$ is fixed or diverges sufficiently slowly.
	Then the RW on $\hk^\pm$ exhibits cutoff whp at
	\[
		\tilde t^\pm_*(k,m,d)
	\cq
		\max\brb{t^\pm_0(k, \abs \hab), \: \log_k \abs \hgr}
	\eqsim
	\begin{cases}
		t^\pm_0(k, \abs \hab)
			&\text{when}\quad
		k \le \rbr{ \log \abs \hab }^{2d-3},
	\\
		\log_k \abs \hgr
			&\text{when}\quad
		k \ge \rbr{ \log \abs \hab }^{2d-3}.
	\end{cases}
	\]
\end{introthm}

\renewcommand{\hgr}{\HH_{m,d}}
\renewcommand{\hab}{\HH_{m,d}^\ab}
\renewcommand{\hk}{(\HH_{m,d})_k}

A highly related theorem holds for diverging $d$.
We establish it as a corollary of a comparison result between mixing times for nilpotent and Abelian groups; see \cite[\cref{res-p2:intro:comp:nil-abe,res-p2:intro:comp:cor:heis}]{HOt:rcg:abe:cutoff}.
We require $k \gg \log m$, but observe that $\log \abs \hab \asymp d \log m$ may be much larger than $\log m$ if $d$ is large or $m$ small; we allow even $m \asymp 1$.
The conditions are weakened in the special case that $m = p$ is prime; see \cite[\cref{res-p2:intro:comp:cor:special}]{HOt:rcg:abe:cutoff}.


\hypersetup{linkcolor = black}
\begin{customrmkt}{\ref{res-p1:intro:ext:comp}/\ref{res-p1:intro:ext:heis}}
\hypersetup{linkcolor = blue}
We strongly believe that one can strengthen \cref{res-p1:intro:ext:comp,res-p1:intro:ext:heis} fairly significantly.
One can define a function $f : \mbn \to [0,\infty)$ which satisfies $f(m) \asymp 1$ for `typical' $m$ and then relax the lower bound $k \gg \max\bra{\sqrt{\log m}, \: \log \div m}$ to $k \gg \log f(m)$.
The term $f(m)$ is a decaying weighted sum of the divisors of $m$.
It is `typically' order $1$: $\ex{ f(U_m) } \asymp 1$ if $U_m \sim \Unif([m])$.

This is work in progress.
It relies on the framework developed here, but requires \emph{significant} additional effort and new ideas.
One needs to refine the analysis of certain gcd terms in \cref{res-p1:ext:comp:k-small:gcd}.
One needs very precise control, namely a local CLT, for a random variable $C_{i,j}$ defined in \eqref{eq-p1:cutoff:3:Cij-def}.
\end{customrmkt}
\hypersetup{linkcolor = blue}

\renewcommand{\ugr}{\UU_{p,d}}
\renewcommand{\uab}{\UU_{p,d}^\ab}
\renewcommand{\uk}{(\UU_{p,d})_k}

Finally, returning to the upper triangular group $\ugr$ with $p$ prime, we can establish the \textit{cutoff profile}
in the regime $1 \ll k \ll \log \abs \uab$.
We use the same techniques in \cite[\S\ref{sec-p2:cutoff1}]{HOt:rcg:abe:cutoff} to find the limit profile when the underlying group $G$ is Abelian.
(There the regime is $1 \ll k \ll \log \abs G = \log \abs \gab$.)

\renewcommand{\ugr}{\UU}
\renewcommand{\uab}{\UU^\ab}
\renewcommand{\uk}{\UU_k}

\begin{introthm}
\label{res-p1:intro:ext:profile}
	Let $p$ be prime and $d \ge 3$.
	Abbreviate $\ugr \cq \UU_{p,d}$.
	Suppose that
		$1 \ll k \ll \log \abs \uab$
	and
		$d \ge 3$ is fixed or diverges sufficiently slowly.
	Then there exist times $(t_\alpha)_{\alpha \in \mbr}$ satisfying
	\[
		t_0 \eqsim k \cdot \tfrac1{2 \pi e} \abs \uab^{2/k},
	\quad
		t_\alpha - t_0 \eqsim \alpha \sqrt 2 t_0 / \sqrt k = \oh{t_0}
	\Qand
		d^\pm_{\uk}(t_\alpha) \eqsim \Psi(\alpha) \ \whp,
	\]
	where $\Psi$ is the standard Gaussian tail, ie $\Psi(\alpha) \cq (2 \pi)^{-1/2} \intt[\infty]{\alpha} e^{-x^2/2} dx$ for $\alpha \in \mbr$.
\end{introthm}

\renewcommand{\ugr}{\UU_{m,d}}
\renewcommand{\uab}{\UU_{m,d}^\ab}
\renewcommand{\uk}{(\UU_{m,d})_k}


\renewcommand{\hgr}{\HH_{p,d}}
\renewcommand{\hab}{\HH_{p,d}^\ab}
\renewcommand{\hk}{(\HH_{p,d})_k}

\renewcommand{\ugr}{\UU_{p,d}}
\renewcommand{\uab}{\UU_{p,d}^\ab}
\renewcommand{\uk}{(\UU_{p,d})_k}

\subsubsection{Remarks on \cref{res-p1:intro:cutoff,res-p1:intro:ext:heis}}


\setcounter{intrormkT}{0}

\begin{subtheorem-num}{intrormkT}
	\label{rmk-p1:intro:cutoff}

\begin{intrormkt}
\label{rmk-p1:intro:cutoff:nec-con}
This article establishes cutoff in a variety of set-ups, but always in the regime $1 \ll \log k \ll \log \abs G$.
This leaves the regimes $k \asymp 1$ and $\log k \asymp \log \abs G$, for which there is no cutoff for any choice of generators:
	when $k \asymp 1$, this holds whenever the group is nilpotent;
	when $\log k \asymp \log \abs G$, this holds for all groups.
The former result is due to \textcite{DSc:growth-rw}; we give a short exposition of this in \cite[\S\ref{sec-p5:const-k}]{HOt:rcg:abe:extra}.
The latter result is proved in \cite[\S\ref{sec-p2:conc-rmks:roichman}]{HOt:rcg:abe:cutoff}; the mixing time~is~order~1.
\textcite[Theorems~3.3.1 and~3.4.7]{D:phd} establishes a more general result for $\log k \asymp \log \abs G$.
\end{intrormkt}
	
\begin{intrormkt}
\label{rmk-p1:intro:cutoff:free-group}
The Abelianisation
$\uab$
is isomorphic to $\mbz_p^{d-1}$; it corresponds to the super-diagonal of the matrices.
Roughly, we split the analysis into ``the mixing of the Abelianisation'' and ``the mixing of the commutator (ie `non-Abelian part')''.
The structure of the proof is the same for all $k$, except in bounding one specific (combinatorial) probability.

We study cutoff when the underlying group $G$ is an arbitrary Abelian group in \cite{HOt:rcg:abe:cutoff}.
The proof goes via lifting the walk $S$ on the Cayley graph to a walk $W$ on the free Abelian group with $k$ generators (namely $\mbz^k$).
The mixing time is then the time at which $W$ has entropy $\log \abs G$.
We perform some analysis on $W$ before projecting back to $S$.

It may seem natural here, then, to lift the random walk to the free nilpotent group of class $d-1$ (ie the nilpotency class of $\ugr$) and to take as the candidate mixing time the time at which this walk has entropy $\log \abs \ugr$.
To the best of our knowledge, the idea of studying RWs on nilpotent groups by lifting the walk to a corresponding free nilpotent group was first used by \textcite{DSc:rw-nilpotent}.
Interestingly, though, instead we still consider a walk on the free \emph{Abelian} group, but now the candidate mixing time is the time at which this walk has entropy $\log \abs \uab$.
At this time, by our results in \cite{HOt:rcg:abe:cutoff}, the walk on the Cayley graph projected to the Abelianisation has~mixed.

Naturally we require the mixing time to be at least $\log_k \abs \ugr$ so that all vertices can be reached with reasonable probability.
We consider the maximum of this entropic time with $\log_k \abs \ugr$.
	%
\end{intrormkt}	

\begin{intrormkt}
\label{rmk-p1:intro:cutoff:heis-to-general}
	Upper triangular matrices and Heisenberg groups are canonical classes of nilpotent groups.
	Our analysis is quite general and extends to other nilpotent groups; see \S\ref{sec-p1:conc-rmks:nil} for a brief overview and \cite[\S\ref{sec-p5:heis-nil}]{HOt:rcg:abe:extra} for more details.
	Hence this article is a first step towards establishing cutoff for other nilpotent groups.
	This is work in progress.
\end{intrormkt}


\end{subtheorem-num}


%

\setcounter{intrormkT}{3}

\begin{intrormkt}
	The fact that the `phase transition' in the mixing time in \cref{res-p1:intro:ext:heis} occurs at exponent $2d-3$ which increases with $d$ gives evidence to the claim that the Heisenberg group $\hgr$ becomes `more Abelian' as $d$ increases.
	In particular, since $\tfrac{2d-3}{2d-4} \to 1$ as $d \to \infty$, if $d$ diverges, then the mixing time for $\hgr$ is asymptotically equivalent to that of an Abelian group of the same size.
	Contrast this with the mixing time for $\ugr$ (from \cref{res-p1:intro:cutoff}) which is of smaller order than this time.
	(Also, $\hgr$ is step-2 nilpotent for all $d \ge 3$ while $\ugr$ is step-$(d-1)$.)
\end{intrormkt}



\subsubsection{Typical Distance and Diameter}

\nextresult

Our next result concerns typical distance in the random Cayley graph.

\begin{introdefn}
	For
		a group $G$,
		$k \in \mbn$
	and
		$\beta \in (0,1)$,
	define the \textit{$\beta$-typical distance} $\mcd_{G_k}(\beta)$ via
	\[
		\mcb^\pm_{G_k}(R) \cq \bra{ x \in G \midb \dist_{G^\pm_k}(\id,x) \le R }
	\Qand
		\mcd^\pm_{G_k}(\beta) \cq \min\brb{ R \ge 0 \midb \abs{ \mcb^\pm_{G_k}(R) } \ge \beta \abs G },
	\]
	with the $\pm$-superscript indicating definitions for both the directed and undirected cases.
\end{introdefn}

Informally,
we show that the mass (in terms of number of vertices) concentrates at a thin `slice', or `shell', consisting of vertices at a distance
$M \pm \oh M$ from the origin, with $M$ explicit.

Investigating this typical distance when $k$ diverges with $\abs G$ was suggested to us by \textcite{B:typdist:private}.
Previous work concentrated on fixed $k$, ie independent of $\abs G$; see \S\ref{sec-p1:intro:previous-work:typdist}.

A more refined statement than the one below is given in \cref{res-p1:dist:res:typdist}.

\renewcommand{\ugr}{\UU}
\renewcommand{\uab}{A}
\renewcommand{\uk}{\UU_k}

\begin{introthm}
\label{res-p1:intro:typdist}
	Let $p$ be prime and $d \ge 3$.
	Let $\UU \cq \UU_{p,d}$ and $A \cq \UU_{p,d}^\ab$.
	Write
	\[
		M^+_k \cq k \abs \uab/e,
	\quad
		M^-_k \cq k \abs \uab^{1/k} / (2e)
	\Qand
		M^*_k \cq \log \abs \uab / \log(k / \log \abs \uab).
	\]
	Assume that $1 \ll \log k \ll \log \abs \ugr$ and that $d$ is either a fixed constant or diverges sufficiently~slowly.
	
	Then,
	for all $\lambda \in (0,\infty)$,
	there exists a constant $\alpha^\pm_\lambda \in (0,\infty)$ so that,
	for all constants $\beta \in (0,1)$,
	the following convergences in probability hold:
	\begin{alignat*}{3}
		\mcd^\pm_{\uk}(\beta) / M^\pm_k &\to^\mbp 1
	&&\Quad{if}&
		1 \ll k &\ll \log \abs \uab;
	\\
		\mcd^\pm_{\uk}(\beta) / \rbb{ \alpha^\pm_\lambda k} &\to^\mbp 1
	&&\Quad{if}&
		k &\eqsim \lambda \log \abs \uab;
	\\
		\mcd^\pm_{\uk}(\beta) / \max\brb{ M^*_k, \: \log_k \abs \ugr } &\to^\mbp 1
	&&\Quad{if}&
		k &\gg \log \abs \uab.
	\end{alignat*}
	Alternatively, the typical distance concentration value can be given by the maximum of $\log_k \abs \ugr$ and the minimal radius of a $k$-dimensional lattice ball of volume at least~$\abs \uab$.
	Note that
	\[
		\max\brb{ M^*_k, \: \log_k \abs \ugr }
	=
		\max\brb{ \tfrac{\rho}{\rho-1}, \: \tfrac12 d } \log \abs \uab.
	\]
\end{introthm}

\renewcommand{\ugr}{\UU_{p,d}}
\renewcommand{\uab}{\UU_{p,d}^\ab}
\renewcommand{\uk}{(\UU_{p,d})_k}

\begin{subtheorem-num}{intrormkT}
	\label{rmk-p1:intro:typdist}

\begin{intrormkt}
\label{rmk-p1:intro:typdist:nil-gen}
	By a classical result, to generate a nilpotent group it is enough that the maps of the generators under $g \mapsto g \gcom : G \mapsto \gab$ generate the Abelianisation $\gab$;
	this follows from the fact that for nilpotent groups $\gcom \le \Phi(G)$, the Frattini subgroup of non-generators of $G$.
	
	We prove a quantitative version of this result, where the typical distance in $G$ is very close to the typical distance in $\gab$ for the Cayley graph with generating multiset $[Z_1\gcom, ... ,Z_k\gcom]$.
	
	See \textcite{EbP:cayley-diam-nil} for a recent different result in the same spirit.
\end{intrormkt}

\begin{intrormkt}
\label{rmk-p1:intro:typdist:interp}
	In \cref{rmk-p1:intro:cutoff:free-group}, we interpret the cutoff time for the RW in the following way:
		if the walk has run for long enough so that the projection to the Abelianisation is mixed
	\emph{and}
		almost every vertex can be reached with reasonable probability,
	then the walk is mixed on the full group.
	\cref{res-p1:intro:typdist} says that the typical distance and mixing time agree when $k \gg \log \abs \uab$; this gives a sense of rigour to the above interpretation.
\end{intrormkt}

\end{subtheorem-num}

\nextresult

In the regime $k \gtrsim \log \abs \uab$, when $d \asymp 1$, we can extend the typical distance argument to determine the \textit{diameter}, ie the maximal distance between pairs of vertices.
In this regime, the two are the same, up to smaller order terms, whp.
For a graph $H$, denote by $\diam H$ its diameter.

\renewcommand{\ugr}{\UU}
\renewcommand{\uab}{A}
\renewcommand{\uk}{\UU_k}

\begin{introthm}
\label{res-p1:intro:diam}
	Let $p$ be prime and $d \ge 3$ fixed.
	Abbreviate $\ugr \cq \UU_{p,d}$.
	Suppose that $1 \ll \log k \ll \log \abs \ugr$.
	For all $\lambda \in (0,\infty)$,
	with $\alpha^\pm_\lambda$ the constant from \cref{res-p1:intro:typdist},
	the following hold:
	\begin{alignat*}{3}
		\rbb{ \diam{} \uk^\pm } / \rbb{ \alpha^\pm_\lambda k} &\to^\mbp 1
	&&\Quad{if}&
		k &\eqsim \lambda \log \abs \uab;
	\\
		\rbb{ \diam{} \uk^\pm } / \max\brb{ M^*_k, \: \log_k \abs \uab } &\to^\mbp 1
	&&\Quad{if}&
		k &\gg \log \abs \uab.
	\end{alignat*}
\end{introthm}

\begin{intrormkt}
	\cref{res-p1:intro:cutoff,res-p1:intro:typdist,res-p1:intro:diam} combined give
	\(
		\tmix(\uk) \asymp \rbr{ \diam{} \uk }^2/k
	\)
	whp.
\end{intrormkt}

\renewcommand{\ugr}{\UU_{p,d}}
\renewcommand{\uab}{\UU_{p,d}^\ab}
\renewcommand{\uk}{(\UU_{p,d})_k}

Interesting is the way we prove \cref{res-p1:intro:typdist}, and by extension \cref{res-p1:intro:diam}.
It is quite common in mixing proofs to use geometric properties of the graph, such as expansion or distance properties.
We, in essence, do the opposite: we adapt the mixing proof to this geometric set-up.
(We give a proof-outline in \S\ref{sec-p1:dist:outline}.)
This is in the same spirit as \cite{LP:ramanujan}; see \S\ref{sec-p1:intro:previous-work:typdist}.

\begin{intrormkt*}
	When $k \gtrsim \log \abs \uab$, for typical distance, and by extension diameter, we can remove the primarily assumption on $p$.
	This involves using ideas in the proof of \cref{res-p1:intro:ext:comp}, ie the extension to non-prime $p$ for cutoff with $k \gtrsim \log \abs \uab$.
	We do not go into detail here.
\end{intrormkt*}

\subsection{Historic Overview}
\label{sec-p1:intro:previous-work}

In this subsection, we give a fairly comprehensive account of previous work on random Cayley graphs, for cutoff and then for typical distance; we compare our results with existing ones.
The occurrence of cutoff in particular has received a great deal of attention over the years.
We also mention, where relevant, other results which we have proved in companion papers;
see also \S\ref{sec-p1:intro:rmks:advert}.

\subsubsection{Universal Cutoff: The Aldous--Diaconis Conjecture}
\label{sec-p1:intro:previous-work:ad-conj}

Since pioneering work of Erd\H{o}s, it has been  understood that the typical behaviour of \emph{random} objects in some class can shed valuable light on the class as a whole.
Thus, when considering some class of combinatorial objects, it is natural to ask questions such as the following.
\begin{itemize}[noitemsep, topsep = \smallskipamount, label = \bcdot] \slshape
	\item What does a typical object in this class `look like'?
	\item If an object is chosen uniformly at random, which properties hold with high probability?
\end{itemize}
\textcite{AD:conjecture} applied this philosophy to the study of random walks on groups.

\smallskip

\textcite[Page~40]{AD:conjecture} stated their conjecture for $k \gg \log \abs G$.
A more refined version is given by \textcite[Conjectures~3.1.2 and~3.4.5]{D:phd}; see also \cite{H:cutoff-cayley->,R:random-random-walks}.
An informal, more general, variant was reiterated by \citeauthor{D:group-rep} in \cite[Chapter~4G, Question~8]{D:group-rep};
he gave some related open questions recently in \cite[\S 5]{D:some-things-learned}.
Towards the conjecture, an upper bound, valid for arbitrary groups, was established by \textcite[Theorem~1]{DH:enumeration-rws} and later \textcite[Theorems~1 and 2]{R:random-random-walks}, who simplified their argument.
A matching lower bound, valid only for Abelian groups, was given by \textcite[Theorem~3]{H:cutoff-cayley->}; see also \textcite[Theorem~5]{H:cutoff-cayley-survey}.
\textcite[Theorem~4]{DH:enumeration-rws} modify the proof of \cite[Theorem~3]{H:cutoff-cayley->} to extend from Abelian groups to some families of groups with irreducible representations of bounded degree.
Combined, this established the Aldous--Diaconis conjecture for Abelian groups and such groups with low degree irreducible representations.
Moreover, the cutoff time was determined explicitly:
it is at
\[
	T(k, \abs G)
\cq
	\log \abs G / \log(k / \log \abs G)
=
	\tfrac{\rho}{\rho-1} \log_k \abs G
\Qwhere
	\text{$\rho$ is defined by $k = \rbr{ \log \abs G }^\rho$}.
\]
(To have $k \gg \log \abs G$, one needs $\rho - 1 \gg 1/\log \log \abs G$.)
See also \textcite{D:phd,H:cutoff-cayley-survey}.

To extend the conjecture to $1 \ll k \lesssim \log \abs G$, one naturally needs some conditions. For example, if $k < d$ and $G \cq \mbz_2^d$, then the group is not generated, and so no mixing can occur.
Subject to suitable conditions, which are nearly optimal in a quantitative sense, in \cite[\cref{res-p2:intro:tv}]{HOt:rcg:abe:cutoff} we establish cutoff for all Abelian groups when $1 \ll k \lesssim \log \abs G$.
From these conditions, one sees that the cutoff time depends only on $k$ and $\abs G$ when $d(G) \ll \log \abs G$ and $k - d(G) \asymp k \gg 1$.~%
However, one can also exhibit a class of groups with $d(G) \asymp \log \abs G$ which exhibit cutoff, but, when $k \asymp \log \abs G$, the time depends on the algebraic structure of $G$ not just on $k$ and $\abs G$.
See \cite[\cref{rmk-p2:intro:cutoff:ad}]{HOt:rcg:abe:cutoff}.

There is a trivial diameter-based lower bound of $\log_k \abs G$.
If $\rho \gg 1$, ie $k$ is super-polylogarithmic in $\abs G$,
then
\(
	T(k, \abs G) \eqsim \log_k \abs G.
\)
Thus cutoff is established for all groups~for~such~$k$.

Even with all this progress for Abelian groups, no-one had established cutoff for any non-Abelian group, unless $k$ grew super-polylogarithmically in $\abs G$; in fact, no lower bound, other than the diameter-based lower bound of $\log_k \abs G$, was known.

\subsubsection{Random Walks on Upper Triangular Matrices}
\label{sec-p1:intro:previous-work:rws-heisenberg}

Random walks on upper triangular matrices have been the focus of a great deal of attention; focus has primarily been on $3 \times 3$ matrices.
(See in particular \cite[\S 2.1]{BDHMW:heisenberg} and \cite[\S 1.1]{PS:heisenberg}, upon which we have based the description below.)
The probabilistic study of random walks on $\ugr$ was initiated by \textcite{Z:heisenberg}; she interpreted the walk in terms of random number generation.
Using a specific generating set of size 4, the correct order of mixing was established by \textcite{DSc:growth-rw,DSc:rw-nilpotent,DSc:nash} using geometric theoretical tools.
Further proofs were given by 	\textcite{D:threads-group,S:rw-evals-gen,S:rw-evals-burnside,S:rw-matrix,BDHMW:heisenberg}.

Moving even further from the realm of Abelian groups, consider $p$ fixed and general $d \ge 3$.
One can consider a simple walk on $\ugr$: a row is chosen uniformly and added to or subtracted from the row above.
\textcite{E:upper-triangular} studied the diameter of the associated Cayley graph, with $d$ growing, subsequently improving this in \textcite{ET:unipotent}.
\textcite{S:rw-matrix} gave mixing bounds via analysis of eigenvalues.
\citeauthor{CP:heisenberg}~\cite{CP:heisenberg,P:heisenberg} look directly at mixing.
This has then been further studied, improved upon and generalised by \textcite{PS:heisenberg,N:rw-upper-tri,NS:heisenberg}; \textcite{NS:heisenberg} are the first to optimise bounds for $p$ and $d$ simultaneously.

In a recent impressive work, \textcite{DH:rw-heisenberg} introduced a new method for proving a central limit theorem for random walks on nilpotent groups.
They illustrate the method on $\ugr$, obtaining some extremely precise results on the rate at which individual coordinates mix as a function of their distance from the diagonal.
They show that the greater the distance from the diagonal is, the faster the mixing time of the coordinate is.
In the same spirit, we show that, in many cases, the bottleneck for mixing is the super-diagonal coordinates, while in the rest of the cases, the cutoff time is given by the diameter-based lower bound $\log_k \abs \ugr$.

\smallskip

Related work includes analysis of the spectrum of a random walk on the upper triangular matrices by \textcite{BVZ:heisenberg-spectrum}.
There has also been work on a random walk on a $3 \times 3$ Heisenberg group with entries in $\mbr$. See, for example, \textcite{B:lie,B:heisenberg}.

\subsubsection{Cutoff for `Generic' Markov Chains and the Entropic Method}
\label{sec-p1:intro:previous-work:generic-ent}

A common theme in the study of mixing times is that `generic' instances often exhibit the cutoff phenomenon.
Moreover, this can often be handled via the \emph{entropic method}, outlined briefly above.
Further details can be found in our companion article, \cite[\S\ref{sec-p2:intro:previous-work:generic-ent}]{HOt:rcg:abe:cutoff}; there the auxiliary process is the same as here.
A detailed description of our application is given in \S\ref{sec-p1:ent:method}.

\subsubsection{Typical Distance and Diameter}
\label{sec-p1:intro:previous-work:typdist}

As well as determining cutoff for these random Cayley graph, we study a geometric property of a diameter flavour; recall the concept of \emph{typical distance} from \cref{res-p1:intro:typdist}.
Previous work on distance metrics (detailed below) had concentrated on the case where the number of generators $k$ is a \emph{fixed} number.
The results establish (non-degenerate)
limiting laws.
This restricts the (sequences of) groups which can be studied; eg, in order for it to be even possible to generate the group---never mind having independent, uniform generators do so whp---one needs $d(G) \le k \asymp 1$.
We discuss generation of groups further in \cite[\S\ref{sec-p5:gen}]{HOt:rcg:abe:extra}; see in particular \cite[\S\ref{sec-p5:gen:k=1}]{HOt:rcg:abe:extra} where we describe adaptations made in order to obtain connected graphs in the references given below.

Our results are in a different direction:
	for us, \toinf k \asinf{\abs G}, and we establish concentration of the observables.
This allows us to consider a much wider range of groups, in particular with $d(G)$ diverging with $\abs G$.
This line of enquiry was suggested to us by \textcite{B:typdist:private}.

\smallskip

\textcite{AGg:diam-cayley-Zq} studied the diameter of the random Cayley graph of cyclic groups of prime order.
They prove that the diameter is order $\abs G^{1/k}$; see \cite[Theorems~1 and 2]{AGg:diam-cayley-Zq}.
They conjecture that the diameter divided by $\abs G^{1/k}$ converges in distribution to some non-degenerate distribution \asinf{\abs G}; see \cite[Conjecture~3]{AGg:diam-cayley-Zq}.

\textcite{MS:diameter-cayley-cycle} consider, as a consequence of a quite general framework, the diameter of the random Cayley graph of $\mbz_n$ with respect to a fixed number $k$ of random generators, for a random $n$, without any primality assumption.
They derive distributional limits for the diameter, the average distance (defined with respect to various $L_p$ metrics) and the girth.
They determine limit distributions for each of these, and in some cases derive explicit formulas.

\textcite{SZ:diam-cayley} build on the framework of \textcite{MS:diameter-cayley-cycle}, again for fixed $k$; they are able to consider non-random $\abs G$, as well as Abelian groups of arbitrary (fixed) rank, instead of only cyclic groups.
In particular, they verify the conjecture of \textcite[Conjecture~3]{AGg:diam-cayley-Zq};
they additionally work with average distance and girth.

\textcite{EbP:cayley-diam-nil} determine, for fixed $k$, the diameter of various finite nilpotent groups of fixed step, including the upper triangular group.
They derive a general inequality, showing that the diameter of a nilpotent group is governed by the diameter of its Abelianisation.
This proves a conjecture of ours from an earlier version of this paper.

\smallskip

\textcite{LP:ramanujan} derive an analogous typical distance result for $n$-vertex, $d$-regular Ramanujan graphs:
	whp all but $\oh n$ of the vertices lie at a distance $\log_{d-1} n \pm \Oh{\log \log n}$;
	they establish this by proving cutoff for the non-backtracking random walk at time $\log_{d-1} n$.

Related work on the diameter of random Cayley graphs, including concentration of certain measures, can be found in \cite{LPSW:cayley-digraphs,S:diam-cayley}.

\smallskip

In a companion article, \cite[\cref{res-p3:intro:typdist}]{HOt:rcg:abe:geom}, we show concentration of typical distance for general Abelian groups.
The Aldous--Diaconis conjecture for mixing can be transferred naturally to typical distance:
	the mass should concentrate at a distance $M$, which can be written as a function only of $k$ and $\abs G$;
	ie there is concentration of mass at a distance independent of the algebraic structure of the group.
For Abelian groups, we find conditions for this to hold, in the regime $1 \ll k \lesssim \log \abs G$, which are almost exactly the same as those for the original (mixing) conjecture (see \cite[\cref{res-p2:intro:tv}]{HOt:rcg:abe:cutoff}).
\cref{res-p1:intro:typdist} shows that the typical distance conjecture does not hold for the upper triangular group.

\subsection{Additional Remarks}
\label{sec-p1:intro:rmks}

\subsubsection{Precise Definition of Cayley Graphs}
\label{sec-p1:intro:cayley-def}

Consider a finite group $G$.
Let $Z$ be a multisubset of $G$.
We consider geometric properties,
namely through distance metrics and the spectral gap,
of the \textit{Cayley graph} of $(G,Z)$;
we call $Z$ the \textit{generators}.
The
	\textit{undirected}, respectively \textit{directed},
\textit{Cayley graph of $G$ generated by $Z$}, denoted
	$G^-(Z)$, respectively $G^+(Z)$,
is the multigraph whose vertex set is $G$ and whose edge multiset is
\[
	\sbb{ \bra{ g,g \cdot z } \mid g \in G, \, z \in Z },
\Quad{respectively}
	\sbb{ \rbr{ g,g \cdot z } \mid g \in G, \, z \in Z }.
\]
If the walk is at $g \in G$, then a step in $G^+(Z)$, respectively $G^-(Z)$, involves choosing a generator $z \in Z$ uniformly at random and moving to $g z$, respectively one of $g z$ or $g z^{-1}$ each with probability~$\tfrac12$.

We focus attention on the \emph{random} Cayley graph defined by choosing $Z_1, ..., Z_k \sim^\iid \Unif(G)$; when this is the case, denote $G^+_k \cq G^+(Z)$ and $G^-_k \cq G^-(Z)$.
While we do not assume that the Cayley graph is connected (ie, $Z$ may not generate $G$), when $G = \UU_{m,d}$ the random Cayley graph $G_k$ is connected whp whenever $k - d \gg 1$; see \cref{rmk-p1:intro:typdist:nil-gen} and \cite[\cref{res-p5:gen:dichotomy}]{HOt:rcg:abe:extra}.


The graph depends on the choice of $Z$.
Sometimes it is convenient to emphasise this;
we use a subscript, writing $\pr[G(z)]{\cdot}$ if the graph is generated by the group $G$ and multiset~$z$.
Analogously, $\pr[G_k]{\cdot}$ stands for the \emph{random} law $\pr[G(Z)]{\cdot}$ where $Z = [Z_1, ..., Z_k]$ with $Z_1, ..., Z_k \sim^\iid \Unif(G)$.

\subsubsection{Typical and Simple Cayley Graphs}
\label{sec-p1:intro:rmks:typ-simp}

The directed Cayley graph $G^+(z)$ is simple if and only if no generator is picked twice, ie $z_i \ne z_j$ for all $i \ne j$.
The undirected Cayley graph $G^-(z)$ is simple if in addition no generator is the inverse of any other, ie $z_i \ne z_j^{-1}$ for all $i,j \in [k]$.
In particular, this means that no generator is of order 2, as any $s \in G$ of order 2 satisfies $s = s^{-1}$---this gives a multiedge between $g$ and $g s$ for each $g \in G$.

The RW on $G^-(z)$ is equivalent to an adjusted RW on $G^+(z)$ where,
when a generator $s \in z$ is chosen,
instead of applying a generator $s$,
either $s$ or $s^{-1}$ is applied, each with probability $\tfrac12$.
Abusing terminology, we relax the definition of simple Cayley graphs to allow order 2 generators, ie remove the condition $z_i \ne z_i^{-1}$ for all $i$.

Given a group $G$ and an integer $k$,
we are drawing the generators $Z_1, ..., Z_k$ independently and uniformly at random.
It is not difficult to see that the probability of drawing a given multiset depends only on the number of repetitions in that multiset.
Thus, conditional on being simple, $G_k$ is uniformly distributed on all simple degree-$k$ Cayley graphs.
Since $k \ll \sqrt{\abs G}$, the probability of simplicity tends to 1 \asinf{\abs G}.
So when we say that our results hold ``\whp (over $Z$)'', we could equivalently say that the result holds ``for almost all degree-$k$ simple Cayley graphs of $G$''.

Our asymptotic evaluation does not depend on the particular choice of $Z$, so the statistics in question depend very weakly on the particular choice of generators for almost all choices.
In many cases, the statistics depend only on $G$ via $\abs G$ and $d(G)$.
This is a strong sense of `universality'.

\subsubsection{Overview of Random Cayley Graphs Project}
\label{sec-p1:intro:rmks:advert}

This paper is one part of an extensive project on random Cayley graphs.
There are
	three main articles \cite{HOt:rcg:abe:cutoff,HOt:rcg:matrix,HOt:rcg:abe:geom}
	(including the current one \cite{HOt:rcg:matrix}),
	a technical report \cite{HOt:rcg:abe:extra}
and
	a supplementary document \cite{HOt:rcg:supp}.
\textit{Each main article is readable independently.}

The main objective of the project is to establish cutoff for the random walk and determining whether this can be written in a way that, up to subleading order terms, depends only on $k$ and $\abs G$; we also study universal mixing bounds, valid for all, or large classes of, groups.
Separately, we study the distance of a uniformly chosen element from the identity, ie typical distance, and the diameter; the main objective is to show that these distances concentrate and to determine whether the value at which these distances concentrate depends only on $k$ and $\abs G$.

\begin{itemize}[noitemsep, topsep = \smallskipamount, label = \bcdot]
	\item [\cite{HOt:rcg:abe:cutoff}]
	Cutoff phenomenon (and Aldous--Diaconis conjecture) for general Abelian groups; also, for nilpotent groups, expander graphs and comparison of mixing times with Abelian groups.
	
	\item [\cite{HOt:rcg:abe:geom}]
	Typical distance, diameter and spectral gap for general Abelian groups.
	
	\item [\cite{HOt:rcg:matrix}]
	Cutoff phenomenon and typical distance for upper triangular matrix groups.
	
	\item [\cite{HOt:rcg:abe:extra}]
	Additional results on cutoff and typical distance for general Abelian groups.
	
	\item [\cite{HOt:rcg:supp}]
	Deferred technical results mainly regarding random walk on $\mbz$ and the volume of lattice balls.
\end{itemize}

\subsubsection{Acknowledgements}
\label{sec-p1:intro:ackno}

This whole random random Cayley graphs project has benefited greatly from advice, discussions and suggestions of many of our peers and colleagues.
We thank a few of them specifically here.

\begin{itemize}[itemsep = 0pt, topsep = \smallskipamount, label = $\bcdot$]
	\item 
	Allan Sly for suggesting the entropy-based approach taken in this paper.
	
	\item 
	Justin Salez for reading this paper in detail and giving many helpful and insightful comments as well as stimulating discussions ranging across the entire random Cayley graphs project.
	
	\item 
	P\'eter Varj\'u for multiple insightful discussions on mixing for the upper triangular group and for more general nilpotent groups, particularly other step-2 nilpotent groups.
	
	\item 
	Itai Benjamini for discussions on typical distance.
	
	\item 
	Evita Nestoridi and Persi Diaconis for general discussions, consultation and advice.
\end{itemize}

\section{Entropic Method and Times}
\label{sec-p1:entropic}

\subsection{Description of Entropic Methodology}
\label{sec-p1:ent:method}

We use an `entropic method', as mentioned in \S\ref{sec-p1:intro:previous-work}; cf \cite{BLPS:giant-mixing,BCS:cutoff-entropic,BL:cutoff-entropic-covered,Ck:cutoff-lifts}.
The method is fairly general; we now explain the specific application in a little more depth.

\medskip

To study the mixing time, we define an auxiliary random process $(W(t))_{t\ge0}$, recording how many times each generator has been used.
More precisely, for $t \ge 0$, for each $i = 1,...,k$, write $W_i(t)$ for the number of times that generator $Z_i$ has been picked by time $t$ for the DRW, and this minus the number of times $Z_i^{-1}$ has been picked for the SRW.
By independence of the coordinates, $W(\cdot)$ forms a rate-1 SRW/DRW on $\mbz^k$ in the un/directed case.

\smallskip

When the group is Abelian, given the choice $Z$ of generators, $S(t)$ is a function of $W(t)$ alone.
When the group is not Abelian, this is not the case.
A key ingredient is to analyse the \textit{commutator} $\ucom$ and \textit{Abelianisation} $\uab$ separately; see \S\ref{sec-p1:cutoff:outline}.
We use the entropic method to analyse the Abelianisation (an Abelian group); we explain the high-level picture now.

Recall that the invariant distribution is uniform, regardless of the group.
For an Abelian group $G$, we propose as the mixing time the time at which the auxiliary process $W$ obtains entropy $\log \abs G$.
The reason for this is the following:
	using the equivalence $- \log \mu \ge \log \abs G$ if and only if $\mu \le 1/\abs G$,
	`typically' $W(t)$ takes values to which it assigns probability smaller than $1/\abs G$;
	informally, this means that $W(t)$ is `well spread out'.
If we could immediately deduce that $S(t)$ typically takes values to which it assigns probability approximately $1/\abs G$, we would be basically done.
However, one could have two independent copies $S$ and $S'$ (using the same generators $Z$) with $S(t) = S'(t)$ but $W(t) \ne W'(t)$; the uniformity of the generators will show that, on average, this is unlikely.
We thus deduce that $S(t)$ is well spread out, ie well mixed.
In contrast, if the entropy is much smaller than $\log \abs G$, then $W(t)$ is not well spread out: it is highly likely to live on a set of size $\oh{1/\abs G}$. The same must then be true for $S(t)$; hence $S(t)$ is not mixed.

\smallskip

For a non-Abelian group, as noted above, just looking at $W(t)$ loses information.
We decompose $\ugr$ into its Abelianisation $\uab$ and commutator $\ucom$.
The above heuristics for Abelian groups suggest that the walk projected to the Abelianisation $\uab$ should be mixed at the time at which $W$ has entropy $\log \abs \uab$.
We then need to check that the walk on the commutator $\ucom$ is mixed at this time.
The Abelianisation $\uab$ corresponds to the super-diagonal.
\textcite{DH:rw-heisenberg} showed that coordinates mix faster the farther they are from the diagonal.
It is thus natural then to expect the commutator to mix faster than the Abelianisation, at least for $d$ not too large.
We need to make sure that all elements of the group can be reached with reasonable probability, and so need to run for at least $\log_k \abs \ugr$.
This suggests $\max\bra{ t_0(k, \abs \uab), \: \log_k \abs \ugr }$ as the mixing~time.

\medskip

To study typical distance, we define a related auxiliary variable, $A$, corresponding to the number of times each generator is used: $A$ is uniformly distributed on a $k$-dimensional lattice ball of a certain radius.
We apply the chosen generators in a uniformly random order.
We do not apply an entropic method here, per se, but the underlying principles of the proof are extremely similar.

\subsection{Definition of Entropic Times and Concentration}
\label{sec-p1:ent:definitions}

In this section, we define the notion of \textit{entropic times}.
For $t \ge 0$, write $\mu_t$, respectively $\nu_t$, for the law of $W(t)$, respectively $W_1(t)$;
so $\mu_t = \nu_t^{\otimes k}$.
Also, for each $i = 1,...,k$, define
\[
	Q_i(t) \cq - \log \nu_{t/k}\rbb{W_i(t)}
\Quad{and set}
	Q(t) \cq - \log \mu_t\rbb{W(t)} = \sumt[k]{1} Q_i(t).
\]
So $\ex{Q(t)}$ and $\ex{Q_1(t)}$ are the entropies of $W(t)$ and $W_1(t)$, respectively. Observe that $t \mapsto \ex{Q(t)} : [0,\infty) \to [0,\infty)$ is a smooth, increasing bijection.

\begin{defn}
\label{def-p1:ent:t0}
	For $k,N \in \mbn$,
	define the \textit{entropic time} $t_0(k,N)$ so that $\ex{Q_1(t_0(k,N))} = \log N/k$.
	
	We apply this with $N \cq p^{d-1} = n^{2/d}$; abbreviate $t_0 \cq t_0(k,p^{d-1}) = t_0(k,n^{2/d})$.
\end{defn}

Direct calculation, with the SRW and Poisson laws gives the following relations.
We sketch the argument in \S\ref{sec-p1:ent:find}; full, rigorous details are given in \cite[\cref{res-p0:se:t0a} and \S\ref{sec-p0:se:entropic-time-calc}]{HOt:rcg:supp}.
Recall that the $+$-superscript corresponds to the DRW and the $-$-superscript to the SRW.

\begin{prop}
\label{res-p1:ent:t0}
Assume that $1 \ll \log k \ll \log N$.
Write $\kappa \cq k/\log N$.
For all $\lambda \in (0,\infty)$,
the following relations hold, for some continuous, decreasing bijection $f^\pm : (0,\infty) \to (0,\infty)$:
\begin{subequations}
	\label{eq-p1:ent:t0}
\begin{alignat*}{2}
	t^\pm_0(k,N) &\eqsim k \cdot N^{2/k}/(2 \pi e)&
\Qwhere
	k &\ll \log N;
\label{eq-p1:ent:t0:<}
\nt
\\
	t^\pm_0(k,N) &\eqsim k \cdot f^\pm(\lambda)&
\Qwhere
	k &\eqsim \lambda \log N;
\label{eq-p1:ent:t0:=}
\nt
\\
	t^\pm_0(k,N) &\eqsim k \cdot 1/(\kappa \log \kappa)&
\Qwhere
	k &\gg \log N.
\label{eq-p1:ent:t0:>}
\nt
\end{alignat*}
\end{subequations}
	%
\end{prop}

By a standard argument considering appropriate subsequences, to cover the general case $k \asymp \log N$, it suffices to assume that $k/\log N$ actually converges, say to $\lambda \in (0,\infty)$.

\smallskip

Since $Q = \sumt[k]{1} Q_i$ is a sum of $k$ iid random variables, $Q(t_0)$ concentrates around $\log N$.
One can show that if the time is multiplied by a factor $1 + \xi$ for any constant $\xi > 0$ then the entropy increases by a significant amount; similarly, if $\xi < 0$ then the entropy decreases by a significant amount.
Further, the change is by an additive term of larger order than the standard deviation $\sqrt{\Var{Q(t_0)}}$.
Thus $Q\rbr{ (1 + \xi) t_0 }$ concentrates around this new value.
In particular, the following~hold:
\[
	\mu_{(1 + \xi) t_0}\rbr{ W\rbr{ (1 + \xi) t_0 } } &= \exp{ - Q\rbr{ (1 + \xi) t_0 } }
\ll
	1/n
\quad
	\whp;
\\
	\mu_{(1 - \xi) t_0}\rbr{ W\rbr{ (1 - \xi) t_0 } } &= \exp{ - Q\rbr{ (1 - \xi) t_0 } }
\gg
	1/n
\quad
	\whp.
\]

The following proposition quantifies this change in entropy and this concentration.
For rigorous details, see \cite[\cref{def-p0:se:t0a,res-p0:se:t0a,res-p0:se:CLT}]{HOt:rcg:supp} in the supplementary material.



\begin{prop}
\label{res-p1:ent:concentration}
	Assume that $k$ satisfies $1 \ll \log k \ll \log N$.
	Then
		$\Var{ Q(t_0) } \gg 1$
	and further, for $\xi \in \mbr \setminus \bra{0}$, writing $v \cq \Var{Q_1(t_0)}$ and $\omega \cq \Var{Q(t_0)}^{1/4} = (vk)^{1/4}$, we have
	\[
		\pr{ Q( (1 + \xi) t_0) \ge \log N \pm \omega }
	\to
		\one{ \xi > 0 }.
	\label{eq-p1:ent:concentration}
	\nt
	\]
	(There is no specific reason for choosing this $\omega$; we just need some $\omega$ with $1 \ll \omega \ll (vk)^{1/2}$.)
\end{prop}

\subsection{Finding the Entropic Times: Sketch Evaluation}
\label{sec-p1:ent:find}

In this subsection, we sketch details towards a proof of \cref{res-p1:ent:t0}.
The full, rigorous details can be found in \cite[\cref{res-p0:se:t0a} and \S\ref{sec-p0:se:entropic-time-calc}]{HOt:rcg:supp},
where all the approximations below are justified.

\medskip

\cref{eq-p1:ent:t0:<}.
When $k \ll \log N$, the target entropy for the rate-$1/k$ RW $W_1$ on $\mbz$ is $\log N / k \gg 1$.
Hence $t_0(k,N) / k \gg 1$.
When a rate-1 RW on $\mbz$ is run for time $s \gg 1$, it approximates a normal distribution with variance $s$; the SRW has mean 0 and the DRW mean $s$.
Direct calculation shows that the entropy of such a normal distribution is precisely
\(
	\tfrac12 \log(2 \pi e s).
\)
Assuming that we can approximate the entropy of the RW by that of the normal distribution sufficiently well (which is precisely what we show in \cite[\cref{res-p0:se:ent:s>1}]{HOt:rcg:supp}),
the claim now follows.

\smallskip

\cref{eq-p1:ent:t0:=}.
When $k \asymp \log N$, the target entropy is $\log N / k \asymp 1$.
So we consider a rate-1 RW on $\mbz$ run for time order 1.
The claim now follows
with $f(\lambda) \cq H^{-1}(1/\lambda)$, where $H(s)$ is the entropy of the rate-1 RW on $\mbz$ run for time $s$.
See \cite[\cref{res-p0:se:ent:s=1|k=logn}]{HOt:rcg:supp} for a more formal treatment.

\smallskip

\cref{eq-p1:ent:t0:>}.
When $k \gg \log N$, the target entropy is $\log N / k \ll 1$.
The rate-1 RW on $\mbz$ run for time $s \ll 1$ is approximated by a Bernoulli distribution with success probability $1 - e^{-s} \approx s$ (along with a uniformly chosen sign for the SRW).
Such a (possibly signed) Bernoulli distribution has entropy $s \log(1/s) + \Oh{s}$.
Again assuming that this approximation can be suitably justified (which is precisely what we show in \cite[\cref{res-p0:se:ent:s<1}]{HOt:rcg:supp}),
the claim now follows.

\section{Cutoff for Random Walks on Upper Triangular Groups}
\label{sec-p1:cutoff}

In this section, we consider mixing for the random walk on the random directed Cayley graph of the upper triangular group $\ugr$.
We take $p$ prime and ``$\equiv$'' means ``equivalent modulo $p$''.

Recall that we denote by
$\ucom = [\ugr, \ugr]$
the \textit{commutator} and
$\uab = \ugr / [\ugr, \ugr]$
the \textit{Abelianisation} (noting that the latter is an Abelian group).
It is straightforward to see that
$\uab \cong \mbz_p^{d-1}$
(corresponding to the super-diagonal terms).
Set $n \cq \abs \ugr = p^{d(d-1)/2}$.

Informally, we show that there is competition between mixing of the Abelianisation and of the commutator.
Which part governs the mixing depends on the regime of $k$:
	for $k \ll \logn[1+2/(d-2)]$, it is the Abelianisation, meaning that the overall mixing time is the same as that for $\mbz_p^{d-1}$;
	for $k \gg \logn[1+2/(d-2)]$, it is non-Abelian part, and the overall mixing time is given by the standard diameter-based lower bound of $\log_k \abs \ugr$;
	see \cref{def-p1:cutoff:t*,res-p1:cutoff:res}.

Throughout this section, we use the following notation:\par
{\centering \smallskip
	\(
		k = \rbb{ \log \abs \uab }{}^\rho,
	\quad
		\tfrac12 d = \rbb{ \log \abs \ugr }{}^\nu
	\Quad{and}
		n = \abs \ugr = p^{d(d-1)/2};
	\) \par\smallskip}\noindent
the choice of $\nu$ is so that $\log \abs \ugr = (\log \abs \uab)^{1+\nu}$, so also $k = \rbr{ \log \abs \ugr }^{\rho/(1+\nu)}$.

\subsection{Precise Statement and Remarks}

In this section, we state the more refined version of \cref{res-p1:intro:cutoff}.

\begin{defn}
\label{def-p1:cutoff:t*}
Define
\(
	\tdiam(k,n) \cq \log_k n.
\)
Define
\[
	t^\pm_*(k,p,d)
\cq
	\max\brb{ t^\pm_0(k, \abs \uab), \: \tdiam(k, \abs \ugr) }.
\]
Abbreviate $\tdiam \cq \tdiam(k,p^{d(d-1)/2})$ and $t^\pm_0 \cq t_0(k,p,d)$.
\end{defn}

The following proposition determines $t_*$ up to a $1 \pm \oh1$ factor; it follows easily from \cref{res-p1:ent:t0,def-p1:cutoff:t*}, using $N \cq \abs \uab = \abs \ugr^{2/d}$.

\begin{prop}
\label{res-p1:cutoff:t*}
We have the following approximation to $t_*$:
\begin{subequations}
	\label{eq-p1:cutoff:t*}
\begin{empheq}
[ left = { t_* \eqsim \empheqlbrace} ]%
{alignat = 3}
	&k \cdot \tfrac1{2 \pi e} \abs \uab^{2/k}&
&\Qwhen&
	1 &\eqmathsbox[\mathrel]{A}{\ll} k \eqmathsbox[\mathrel]{B}{\ll} \log \abs \uab;
\label{eq-p1:cutoff:t*:<}
\\
	&k \cdot f(\lambda)&
&\Qwhen&
	&\mathrel{\phantom{\eqmathsbox[\mathrel]{A}{\ll}}} k \eqmathsbox[\mathrel]{B}{\eqsim} \lambda \log \abs \uab;
\label{eq-p1:cutoff:t*:=}
\\
	&\tfrac{\rho}{\rho - 1} \tfrac2d \log_k \abs \ugr&
&\Qwhen&
	\log \abs \uab &\eqmathsbox[\mathrel]{A}{\ll} k \eqmathsbox[\mathrel]{B}{\le} (\log \abs \uab)^{1+2/(d-2)};
\label{eq-p1:cutoff:t*:>}
\\
	&\log_k \abs \ugr&
&\Qwhen&
	(\log \abs \uab)^{1+2/(d-2)} &\eqmathsbox[\mathrel]{A}{\le} k, \; \log k \ll \log \abs \ugr;
\label{eq-p1:cutoff:t*:>>}
\end{empheq}
\end{subequations}
here $f$ is the function from \cref{res-p1:ent:t0}.
(The third regime is empty if $(\log \abs \uab)^{1/d} \asymp 1$, ie $d \gtrsim \log \log p$; in this case, the lower bound in the fourth regime becomes $k \gg \log \abs \uab$.)
\end{prop}

There are some simple conditions that the parameters must satisfy for our proof to be valid.
The conditions will be assumed throughout the remainder of the section, often not explicitly stated.

\begin{hyp}
\label{hyp-p1:cutoff}
	The triple $(k,p,d)$ satisfies \textit{\cref{hyp-p1:cutoff}} if the following conditions hold:
	\begin{itemize}[noitemsep, topsep = \smallskipamount, label = \bcdot]
		\item 
		$1 \ll \log k \ll \log \abs \ugr$;
		
		\item 
		if $k \ll \log \abs \uab$, then $d^3 \ll k$ and $k \ll \log \abs \uab / \log d$ (equivalently $k \ll d \log p / \log d$);
		
		\item 
		if $k \gtrsim \log \abs \uab$, then $\log d \ll \log \log p$ (equivalently $\log d \ll \log \log \abs \ugr$).
	\end{itemize}
	(Recall that implicitly we consider sequences $(k_N, p_N, d_N)_\Ninn$.)
\end{hyp}

\begin{rmkt*}
	\cref{hyp-p1:cutoff} holds when $d \asymp 1$ and $1 \ll \log k \ll \log \abs \uab$.
	As noted in the introduction, there is no cutoff for $k$ outside this regime.
	Thus our conditions are optimal when $d \asymp 1$.
\end{rmkt*}

%
%
%
%

We now state the main result of this section; it is in essence a restatement of \cref{res-p1:intro:cutoff}.

\renewcommand{\ugr}{\UU_{p,d}}
\renewcommand{\uab}{\UU_{p,d}^\ab}
\renewcommand{\uk}{(\UU_{p,d})_k}

\begin{thm}[Cutoff]
\label{res-p1:cutoff:res}
	Let $k,p,d \in \mbn$ with $p$ prime and $d \ge 3$, jointly satisfying \cref{hyp-p1:cutoff}.
	Then the RW on $\uk^\pm$ exhibits cutoff at time $t^\pm_*(k,p,d)$, given in \cref{def-p1:cutoff:t*}, \whp.
	Moreover, the implicit lower bound on mixing holds \forallZ and requires only $1 \ll \log k \ll \log \abs \ugr$.
\end{thm}

\begin{rmkt*}
For ease of presentation, consider for the moment $d$ independent of $n$.
Define
\[
	T(\rho,N)
\cq
	\tfrac{\rho}{\rho - 1} \log_k N
=
	\tfrac{\rho}{\rho - 1} \tdiam(k,N).
\]
Simple algebraic manipulations give
\(
	T(\rho,N) \eqsim t^\pm_0\rbr{ \rbr{ \log N }^\rho, N }
\)
when $\rho > 1$ is bounded away~from 1.
This is the universal mixing time upper bound for a group of size $N$ from \S\ref{sec-p1:intro:previous-work:ad-conj} when $k = \rbr{ \log N }^\rho$ with $\rho > 1$.
(It is tight for Abelian groups.)
Recall that the Abelianisation
has size $\abs \uab = p^{d-1}$.

Consider $k = \rbr{ \log \abs \uab }^\rho$ with $\rho > 1$.
Hence the walk projected to the Abelianisation has cutoff at $t_0(k, \abs \uab)$.
Our proof shows that the random walk on the whole group exhibits cutoff, with time given by the maximum of this and the diameter lower bound, $\tdiam(k, \abs \ugr)$.

This heuristic is only valid when $\rho > 1$.
From \cite[\cref{res-p2:cutoff1:res}]{HOt:rcg:abe:cutoff}, one sees that the walk projected to the Abelianisation has cutoff at $t_0(k, \abs \uab)$ for $\rho \le 1$ too; in this regime, $t_0(k, \abs \uab) \gg \tdiam(k, \abs \ugr)$.
We show that the mixing time is given by $t_0(k, \abs \uab)$ when $\rho \le 1$.
\end{rmkt*}

The fact that the mixing time is a maximum of two quantities suggests some sort of `competition' between the Abelianisation and the rest of the group; this leads to a `phase transition' in the mixing time, which has an interesting consequence for the Aldous--Diaconis conjecture.

\begin{rmkt*}
Consider for the moment $\rho \cq 1 + \tfrac1d$.
If $d \to \infty$ sufficiently slowly, then
\[
	\log \abs \uab
\ll
	(\log \abs \uab)^{1+1/d}
\ll
	(\log \abs \uab)^{1+2/(d-2)}.
\]
As stated previously, $T(\rho, \abs G)$ is a universal mixing bound (valid for all groups) which is tight for Abelian groups.
However, \cref{res-p1:cutoff:t*,res-p1:cutoff:res} show that the mixing time $t_*$ for the upper triangular group satisfies
\(
	t_* \eqsim T(\rho, \abs \uab) = \tfrac2d T(\rho, \abs \ugr)
\)
in this regime.
Hence the time $T(\rho, \abs \ugr)$ is off by a factor of $\tfrac2d$;
it thus does not even capture the correct order of the mixing (since we allow $d \to \infty$).
The upper triangular groups thus provide a counter-example to the claim that $T(\rho, \abs G)$ is the mixing time for \emph{all} groups in the regime $k = (\log \abs G)^\rho \gg \log \abs G$.
\end{rmkt*}


Recall that cutoff is already established (for all groups) when $k$ grows super-polylogarithmically in $n$, ie $\log k \gg \log \log \abs \ugr$.
Below assume that $\log k \lesssim \log \log \abs \ugr$, ie $k = \rbr{ \log \abs \ugr }^{\Oh1}$.

\subsection{Outline of Proof}
\label{sec-p1:cutoff:outline}

We now give a high-level description of our approach, introducing notations and concepts along the way.
No results or calculations from this section will be used in the remainder of the document.
Recall the definitions from the previous sections.

\renewcommand{\ugr}{\UU}
\renewcommand{\uab}{A}
\renewcommand{\ucom}{\UU^\com}

For ease of notation, we suppress the $p$ and $d$ dependence from $\UU_{p,d}$, writing just $\ugr$.
Similarly we write $\uab \cq \UU_{p,d}^\ab$ for the Abelianisation and $\ucom \cq \UU_{p,d}^\com$ for the commutator.

\medskip

We start by discussing the lower bound.
In \cite[\S\ref{sec-p2:cutoff1:lower}]{HOt:rcg:abe:cutoff} we consider an analogous entropic lower bound but where the underlying group is Abelian.
To apply this, we simply project the walk from $\ugr$ to $\uab$, which is an Abelian group.
Projection cannot increase the TV distance.

If $Q$ is sufficiently small, then $W$, and hence also $S$, is restricted to a small set.
Indeed,
\(
	Q \le \log \abs \uab - \omega
\Quad{if and only if}
	\mu\rbr{W} \ge \abs \uab^{-1} e^\omega,
\)
and thus if this is the case then
\(
	W \in \bra{ w \mid \mu(w) \ge \abs \uab^{-1} e^\omega }.
\)
Write $S^\ab$ for $S$ projected to the Abelianisation $\uab$.
Since $\uab$ is an Abelian group, $S^\ab(t)$ depends only on $W(t)$ (not additionally any $W(t')$ for $t' < t$).
It is thus also the case that
\[
	S^\ab(t)
\in
	E
\cq
	\bra{ a \in \uab \mid \pr{ S^\ab(t) = a} \ge \abs \uab^{-1} e^\omega }.
\]
But clearly $\abs{E} \le e^{-\omega} \abs \uab$.
Choosing the time $t$ slightly smaller than the entropic time $t_0$ and $\omega \gg 1$ suitably, the event
\(
	\bra{ Q(t) \le \log \abs \uab - \omega }
\)
will hold whp.
Thus, whp, $S^\ab(t)$ is restricted to a set of size $\oh{\abs \uab}$.
It hence cannot be mixed.
This heuristic applies for any choice of generators.

Precisely,
we show
	for any $\omega$ with $1 \ll \omega \ll \log \abs \uab$,
	all $t$
and
	all $Z = [Z_1, ..., Z_k]$,
that
\[
	d_{G_k}(t)
\ge
	\pr{ Q(t) \le \log \abs \uab - \omega } - e^{-\omega}.
\]
Observe that the probability on the right-hand side is independent of $Z$. Thus we are naturally interested in the fluctuations of $Q(t)$ for $t$ close to $t_0$.
Using the concentration of $Q$, ie \cref{res-p1:ent:concentration} with $\xi < 0$ and $\omega \cq \Var{Q(t_0)}^{1/4}$, we deduce the lower bound in \cref{res-p1:cutoff:res}.

\medskip

We now turn to discussing the upper bound.
We use a \emph{modified $L_2$ calculation}; see \cref{res-p1:cutoff:prelim:D-equiv}.
If $S$ and $S'$ are independent copies, with auxiliary processes $W$ and $W'$, respectively, and $\mcw \subseteq \mbz^k$ is some set, then applying first the triangle inequality then a standard $L_2$ calculation gives
\[
	\ex{ \tvb{ \pr{ S(t) \in \cdot } - \pi_G } }
\le
	\tfrac12 \sqrt{ \abs \ugr \, \pr{ S(t) = S'(t) \mid W(t), W'(t) \in \mcw } - 1 } + \pr{ W(t) \notin \mcw }.
\]
We choose $\mcw$ so that $\pr{ W(t) \notin \mcw } = \oh1$ and think of $\mcw$ as the set of `typical $W(t)$'.
By imposing some mild typicality conditions, we show that $S'(t) S(t)^{-1}$ is uniformly distributed on $\ugr$ when $W(t) \ne W'(t)$ with both $W(t)$ and $W'(t)$ typical.
It thus remains to show that
\[
	\abs \ugr \, \pr{ S(t) = S'(t), \: W(t) = W'(t) \midb W(t), W'(t) \in \mcw } - 1
=
	\oh{1/\abs G}.
\]
This is where we analyse the Abelianisation and commutator separately.
Drop the $t$ from the notation.
Also write
\(
	\widebar \mbp\rbr{ \cdot }
\cq
	\pr{ \, \cdot \mid W(t), W'(t) \in \mcw }.
\)

For the regime in which the mixing time is the entropic time $t_0(k, \abs \uab)$,
for the Abelianisation, we add an entropic condition to typicality:
	if $w \in \mcw$, then $\mu(w) \le \abs \uab^{-1} e^{-\omega}$.
By definition of the entropic time, like in the lower bound, we have
\[
	\widebar \mbp\rbb{ W = W' }
\le
	e^{-\omega} \abs \uab^{-1}
\ll
	\abs \uab^{-1};
\]
see \cref{res-p1:cutoff:prelim:W=W'}.
If $W = W'$, then $S' S^{-1} \in \ucom$.
Given $W = W'$, we desire $S' S^{-1}$ to be approximately uniformly distributed on $\ucom$, say with modal probability order $\abs \ucom^{-1}$,
as then
\[
	\widebar \mbp\rbb{ S = S', \: W = W' }
=
	\widebar \mbp\rbb{ S = S' \midb W = W' }
\cdot
	\widebar \mbp\rbb{ W = W' }
\ll
	\abs \uab^{-1} \cdot \abs \ucom^{-1}
=
	\abs \ugr^{-1},
\]
as desired.
When the mixing time is $\log_k \abs \ugr$, then we perform an analogous analysis, but this time we calculate the entropy shortly after $\log_k \abs \ugr$, which is larger than $t_0(k, \abs \uab)$; see \cref{res-p1:cutoff:prelim:h0,res-p1:cutoff:prelim:W=W'}.
We can then relax the ``approximate uniformity'' appropriately.

\smallskip

We now briefly explain how to establish this ``approximate uniformity'' of $S' S^{-1}$ on $\ucom$ given $W = W'$.
We explain the method for $d = 3$; general $d$ imposes some additional technical hurdles.

We define a combinatorial event $\mce$ in terms of $W$ and $W'$,
which depends on the order in which the generators are chosen not just on the final counts $W(t)$ and $W'(t)$.
We can describe this event in terms of a random walk on a free nilpotent group.
Let $N_k$ be a free nilpotent group of step 2 with $k$ generators.
Let $\widetilde S$ be the RW on $N_k$; assign to it auxiliary process $W$.
Let $(\widetilde S', W')$ be an independent copy of $(\widetilde S, W)$.
Let $w \in \mbz^k$.
Given $W(t) = w = W'(t)$, the combinatorial event $\mce$ is exactly the event
\(
	\bra{ \widetilde S'(t) \widetilde S(t)^{-1} \in [N_k, N_k] }
\setminus
	\bra{ \widetilde S'(t) \widetilde S(t)^{-1} \in [N_k, [N_k, N_k]] \setminus \bra{\id} }.
\)
It is interesting that this condition turns out to be the relevant condition for arguing that $S'(t) S(t)^{-1}$ is roughly a uniformly distributed commutator $\ucom = [\ugr, \ugr]$.
We plan to investigate this further in the context of general step 2 nilpotent groups in future work.


It remains to control the probability of $\mce$; we again use typicality conditions for this.
This is the only place in which the method differs according to the regime of $k$; see \cref{res-p1:cutoff:3:prod-bound}.

When the mixing time is the entropic time, this `error probability' will be smaller than $1/\abs \ucom$, meaning that
\(
	\widebar \mbp\rbr{ S = S', \: W = W' } = \oh{1/\abs \ugr},
\)
as described above.
When the mixing time is $\log_k \abs \ugr$, the error is larger, but combined with $\pr{W = W'}$ gives $\oh{1/\abs \ugr}$; see \cref{res-p1:cutoff:d:prod-bound}.

\subsection{Lower Bound}
\label{sec-p1:cutoff:lower}

\renewcommand{\ugr}{\UU}
\renewcommand{\uab}{A}
\renewcommand{\ucom}{\UU^\com}

The lower bound is relatively straightforward to prove: we project onto the Abelianisation, then use the lower bound for Abelian groups from our companion article \cite{HOt:rcg:abe:cutoff}, specifically \cite[\S\ref{sec-p2:cutoff1:lower}]{HOt:rcg:abe:cutoff}.

For ease of notation, we suppress the $p$ and $d$ dependence from $\UU_{p,d}$, writing just $\ugr$.
Similarly we write $\uab \cq \UU_{p,d}^\ab$ for the Abelianisation.
Also write $n \cq \abs \ugr = p^{d(d-1)/2}$.

\begin{Proof}[Proof of Lower Bound in \cref{res-p1:cutoff:res}]
We assume that $Z$ is given, and suppress it.

For any $\eps > 0$, a lower bound is given by $(1 - \eps) \log_k n$:
	in $m$ steps the support of the random walk is (at most) $k^m$ and hence the walk cannot be mixed in this many steps; cf \cite[Fact~2.1]{LP:ramanujan}.

Write $\Pi : \ugr \to \uab$ for the canonical projection.
Write $N \cq \abs \uab = p^{d-1}$.
Let $\eps > 0$ and let $t \cq (1 - \eps) t_0(k,N)$.
Write
\[
	\mce
\cq
	\brb{ \mu\rbb{W(t)} \ge N^{-1} e^\omega }
=
	\brb{ Q(t) \le \log N - \omega },
\]
with $\mu$, $Q$ and $\omega \gg 1$ from \S\ref{sec-p1:ent:definitions}.
By \cref{res-p1:ent:concentration}, we have $\pr{\mce} = 1 - \oh1$.

For $w \in \mbz_+^k$ and $z_1, ..., z_k \in \ugr$,
write
\(
	z^w
\cq
	z_1^{w_1} \cdots z_k^{w_k}.
\)
Recall that reordering the terms corresponds to multiplication by a particular element of the commutator.
Consider the set
\[
	E
\cq
	\brb{ x \in \uab \midb \exists \, w \in \mbz_+^k \st \mu_t(w) \ge N^{-1} e^\omega \text{ and } x = \Pi(Z^w) }
\subseteq
	\uab.
\]
Since we use $W$ to generate $S$, we have $\pr{ \Pi\rbr {S(t) } \in E \mid \mce } = 1$.
Every element $x \in E$ satisfies $x = \Pi(w_x \bcdot Z)$ for some $w_x \in \mbz_+^k$ with $\mu_t(w_x) \ge N^{-1} e^\omega$.
Hence, for all $x \in E$, we have
\[
	\pr{ \Pi(S(t)) = x }
\ge
	\pr{ W(t) = w_x }
=
	\mu_t(w_x)
\ge
	N^{-1} e^\omega.
\]
Taking the sum over all $x \in E \subseteq \uab$,
we deduce that
\[
	1
\ge
	\sumt{x \in E}
	\pr{ \Pi(S(t)) = x }
\ge
	\abs E \cdot N^{-1} e^\omega,
\Quad{and hence}
	\abs E/N \le e^{-\omega} = \oh1.
\]
Finally we deduce the lower bound from the definition of TV distance:
\[
	\tvb{ \pr{ S(t) \in \cdot \mid Z } - \pi_G }
\ge
	\pr{ S(t) \in \Pi^{-1}(E) } - \pi_G\rbb{\Pi^{-1}(E)}
\ge
	\prt{ \mce } - \tfrac1N \abs{E}
\ge
	1 - \oh1.
\qedhere
\]
\end{Proof}

\begin{rmkt}
\label{rmk-p1:cutoff:general-lower}
For the entropic lower bound, all that we used was the size of the Abelianisation.
The same argument shows,
for
	all finite groups $G$,
	all $k \gg 1$
and
	all multisubsets $Z$ of $G$ of size $k$,
that
\(
	\max\bra{ t_0(k, \abs \gab), \log_k \abs G }
\)
is a lower bound on the mixing time.

The lower bound
here
can be used to determine the profile of the convergence to equilibrium; this is done in \S\ref{sec-p1:ext:window}, and in \cite[\S\ref{sec-p2:cutoff1}]{HOt:rcg:abe:cutoff}.
Another lower bound is proved in \cite[\S\ref{sec-p2:cutoff2:lower}]{HOt:rcg:abe:cutoff}; this cannot be used to determine the profile.
For many groups these lower bounds will be equivalent, but for some the latter captures the correct mixing time while the former is a constant factor too~small.
\end{rmkt}


\subsection{Upper Bound Preliminaries}
\label{sec-p1:cutoff:prelim}

\renewcommand{\ugr}{\UU}
\renewcommand{\uab}{A}
\renewcommand{\ucom}{\UU^\com}

For ease of notation, we suppress the $p$ and $d$ dependence from $\UU_{p,d}$, writing just $\ugr$.
Similarly we write $\uab \cq \UU_{p,d}^\ab$ for the Abelianisation.
Also write $n \cq \abs \ugr = p^{d(d-1)/2}$.

We first prove the upper bound for $d = 3$.
	The majority of the ideas are exposed in this case, while the technical details involved in the general $d$ case somewhat obscure the ideas.
Note that the conditions on $d$ from \cref{hyp-p1:cutoff} is always satisfied when $d = 3$ (or, in fact, any fixed $d$).
Similarly, we first analyse the DRW (ie directed graphs); in \S\ref{sec-p1:cutoff:undir} we then describe the (simple) adaptations to the proof required for the SRW.

\medskip

Before doing so, we need some preliminary results (for both $d = 3$ and general $d$).
First, we need a concept of `typicality' for the auxiliary random variable $W(t)$:
later in the proof, we define a set $\mcw \subseteq \mbz_+^k$ (dropping the $t$-dependence from the notation) with the property
\[
	\mcw
\subseteq
	\brb{ w \in \mbz_+^k \midb \mu_t(w) \le e^{-h}, \: \maxt{i} w_i < \tfrac12 p },
\label{eq-p1:cutoff:prelim:typ:subset}
\nt
\]
where $h$ is roughly the entropy of $W(t)$; see \cref{def-p1:cutoff:prelim:h0} for the precise definition of $h$, and also \cref{res-p1:cutoff:prelim:h0} for the relation to the entropy.
This set will satisfy $\pr{W \in \mcw} = 1 - \oh1$, hence the name `typical'; see \cref{res-p1:cutoff:prelim:typ}.

It is often easier to work with $L_2$, rather than TV, distance, since it has a nice explicit representation; on the other hand, with TV one can condition on `typical' events.
We combine the two with a `modified $L_2$ calculation':
	let $S$ and $S'$ be independent copies (given $Z$), and let $W$ and $W'$ be their associated auxiliary random variables; write $\typ \cq \bra{ W,W' \in \mcw }$.

\begin{lem}
\label{res-p1:cutoff:prelim:D-equiv}
	Assume that
	\(
		\pr{ W((1+\eps)t_*) \notin \mcw } = \oh1
	\)
	for all constants $\eps > 0$.
	Then the upper bound in \cref{res-p1:cutoff:res} is established by showing, for all constants $\eps > 0$, that
	\[
		D(t) \cq \abs \ugr \, \pr{ S = S' \mid \typ } - 1
	\Quad{satisfies}
		D\rbb{ (1+\eps) t_* } = \oh1.
	\]
\end{lem}

\begin{Proof}
Using the triangle inequality and then Cauchy-Schwarz inequality, we obtain the following:
\begin{gather*}
	\tvb{ \pr{ S \in \cdot \mid Z } - \pi_G }
\le
	\tvb{ \pr{ S \in \cdot \mid Z, \: W \in \mcw } - \pi_G }
+	\pr{ W \notin \mcw };
\\
	\ex{ 2 \, \tvb{ \pr{ S \in \cdot \mid Z, \: W \in \mcw } - \pi_G } }^2
\le
	\abs \ugr \, \pr{ S = S' \mid \typ } - 1
=
	D(t).
\end{gather*}
Combining these with the assumption $\pr{ W \notin \mcw } = \oh1$, and Markov's inequality, gives
\[
	\tvb{ \pr{ S \in \cdot \mid Z } - \pi_G }
=
	\oh1
\quad
	\whp.
\qedhere
\]
	%
\end{Proof}

To upper bound $D \cq D(t)$, we separate into cases according to whether or not $W = W'$:
\[
	\pr{ S = S' \mid \typ }
&
=
	\pr{ S = S' \mid W = W', \: \typ } \pr{ W = W' \mid \typ }
\\&
+
	\pr{ S = S' \mid W \ne W', \: \typ } \pr{ W \ne W' \mid \typ }.
\]
Were the underlying group Abelian, $W = W'$ would imply $S = S'$.
This is not the case for non-Abelian groups; in fact estimating $\pr{ S = S' \mid W = W', \: \typ }$ is the main part of the proof.

\smallskip

First we control $\pr{W = W' \mid \typ}$.
To do this, we must estimate the entropy shortly after the proposed mixing time.
Recall
	that $W(\cdot)$ is a RW on $\mbz_+^k$,
	that $\mu_t$ is the law of $W(t)$
and
	that $Q(t) = - \log \mu_t(W(t))$.
Denote by
\[
	h(t) = \ex{Q(t)},
\quad
	\text{the entropy of $W(t)$}.
\]
Recall that
	$\tdiam = \log_k \abs \ugr$,
	$\tfrac12 d = (\log \abs \uab)^\nu$,
	$k = \rbr{ \log \abs \uab }^\rho = \rbr{ \log \abs \ugr }^{\rho/(1+\nu)}$
and
	$\abs \ugr = \abs \uab^{d/2}$.

\begin{lem}
\label{res-p1:cutoff:prelim:h0}
Let $\xi > 0$.
Then, for any $\omega \ll \min\bra{k, \: \log \abs \uab}$, the following lower bounds hold.
\begin{itemize}[noitemsep, topsep = \smallskipamount, label = $\bcdot$]
	\item 
	For $t \ge (1 + \xi) t_0(k, \abs \uab)$,
	we have
	\(
		h(t) = \ex{Q(t)} \ge \log \abs \uab + 2 \omega.
	\)
	
	\item 
	For $t \ge (1 + \xi) \tdiam$,
	if $\rho \ge 1 + 2/(d-2)$,
	then
	\(
		h(t) = \ex{Q(t)} \ge (1 - \tfrac1\rho) \log \abs \ugr + 2 \omega.
	\)
\end{itemize}
\end{lem}

To prove this lemma, we use the following result, which will be used independently later.

\begin{lem}
\label{res-p1:cutoff:prelim:omega-choice}
	Let $N \gg 1$ be any diverging integer.
	Let $t_0$ and $t_{2\omega}$ be the entropic times for entropy $\log N$ and $\log N + 2 \omega$, respectively.
	Then we have $t_{2\omega} \eqsim t_0$ if $\omega \ll \min\bra{k, \log N}$.
\end{lem}

We defer the proof of \cref{res-p1:cutoff:prelim:omega-choice} to \cite[\cref{res-p0:se:entropic-time-calc:+omega}]{HOt:rcg:supp}.
We now prove \cref{res-p1:cutoff:prelim:h0}.

\begin{Proof}[Proof of \cref{res-p1:cutoff:prelim:h0}]
Consider first time $t_0(k, \abs \uab)$.
We have $h(t) \ge \log \abs \uab$ by definition of the entropic time.
The $+ 2 \omega$ additive term then follows immediately from \cref{res-p1:cutoff:prelim:omega-choice}.

Consider now the time $\tdiam$.
Recall from \S\ref{sec-p1:ent:find}
that the entropy of the rate-1 RW on $\mbz$ at time $s \ll 1$ satisfies
\(
	H(s) \eqsim s \log(1/s).
\)
Take $s \cq \tdiam/k$.
Direct calculation gives $s \ll 1$.
%
Thus
\[
	h(\tdiam) / \log n
&
\eqsim
	\tdiam \log(k/\tdiam) / \log n
=
	\tfrac1{\log k} \log\rbb{ \tfrac{k \log k}{\log n} }
\\&
=
	\tfrac1{\rho \log \log \abs \uab } \rbb{ (\rho - 1 - \nu) \log \log \abs \uab + \log \log( \rbr{ \log \abs \uab }^\rho ) }
\ge
	1 - \tfrac{1+\nu}\rho.
\]
For $\xi \in (0,1)$ fixed,
\(
	h\rbr{ (1 + \xi) \tdiam }
\eqsim
	(1 + \xi) h\rbr{ \tdiam }.
\)
The conditions on $d$ gives
\(
	1 - 1/\rho \ge 2/d \gg \nu.
\)
Hence the claim is true for all $\xi > 0$ (fixed) provided $k \ll \tfrac2d \log n$, ie $k \ll \log \abs \uab$.
\end{Proof}

Motivated by this lemma, recalling that $t_* = \max\bra{t_0, \tdiam}$, we make the following definition.

\begin{defn}
\label{def-p1:cutoff:prelim:h0}
Define $h_0$ as follows:
\begin{empheq}
[ left = {h_0 \cq \empheqlbrace} ]
{alignat = 3}
	&\log \abs \uab
		&\Qwhen&
	k \le (\log \abs \uab)^{1 + 2/(d-2)};
\nonumber
\\
	&(1 - \tfrac1\rho) \log \abs \ugr
		&\Qwhen&
	k \ge (\log \abs \uab)^{1 + 2/(d-2)}.
\nonumber
\end{empheq}
Fix some $\omega$ such that $1 \ll \omega \ll \min\bra{k, \: \log \abs \uab}$, and set $h \cq h_0 + \omega$.
\end{defn}

Not only does the entropy satisfy this lower bound, but the $Q$ random variable, which is defined, in \S\ref{sec-p1:ent:definitions}, so that $\ex{Q(t)} = h(t)$, concentrates, giving the following result.


\begin{lem}
\label{res-p1:cutoff:prelim:typ}
	Assume that $\omega \ll \min\bra{k, \: \log \abs \uab}$.
	Let $\eps > 0$ and $t \ge (1 + 3\eps) \max\bra{t_0, \tdiam}$.
	Then
	\[
		\pr{ Q(t) \ge h } = \pr{ \mu_t\rbb{W(t)} \le e^{-h} } = 1 - \oh1.
	\]
\end{lem}

\begin{Proof}
	Rearrange the inequality
	\(
		\mu \le e^{-h}
	\)
	into
	\(
		Q \cq - \log \mu \ge h,
	\)
	use \cref{res-p1:cutoff:prelim:h0},
	the definition of $h$ and $h_0$ from \cref{def-p1:cutoff:prelim:h0}
and
	apply the concentration result \cref{res-p1:ent:concentration}.
\end{Proof}

We now control $\pr{ W = W' \mid \typ }$.
This is where the typicality condition in \cref{eq-p1:cutoff:prelim:typ:subset} comes in.

\begin{lem}
\label{res-p1:cutoff:prelim:W=W'}
	Recall $h$ as defined in \cref{def-p1:cutoff:prelim:h0}.
	We have
	\[
		\pr{ W = W' \mid \typ }
	\le
		e^{-h} / \prt{\typ}.
	\]
\end{lem}

\begin{Proof}
By direct calculation,
since $W$ and $W'$ are independent copies,
we have
\[
	\pr{ W = W', \: \typ }
=
	\pr{ W = W', \: W \in \mcw }
=
	\sumt{w \in \mcw}
	\pr{ W = w } \pr{ W' = w }
\le
	e^{-h},
\]
using the fact that $\sumt{w \in \mcw} \pr{W = w} \le 1$ and $\pr{W' = w} \le e^{-h}$ for all $w \in \mcw$.
\end{Proof}

Consideration of $\pr{ S = S' \mid W \ne W', \: \typ }$ is the topic of the next two subsections (\S\ref{sec-p1:cutoff:3} for $d = 3$ and \S\ref{sec-p1:cutoff:d} for general $d$).
The main ingredient is the following lemma, which follows immediately from the fact that $p$ is prime. (Recall that $[m] = \bra{1, ..., m}$ for $m \in \mbn$.)

\begin{lem}
\label{res-p1:a.X-unif}
	Let $\ell \in \mbn$, $X_1, ..., X_\ell \sim^\iid \Unif(\mbz_p)$ and $a_1, ..., a_\ell \in [p-1]$.
	Then
	\(
		\sumt[\ell]{1} a_i X_i \sim \Unif(\mbz_p).
	\)
\end{lem}

We can extend this to general $p \in \mbn$: we have $\sumt[\ell]{1} a_i X_i \sim \Unif(\mfgcd \mbz_p)$ where $\mfgcd \cq \gcd(a_1, ..., a_\ell, p)$ and $\mfgcd \mbz_p = \bra{g, 2g, ..., p}$; see \cite[\cref{res-p2:cutoff1:unif-gcd}]{HOt:rcg:abe:cutoff}.
(These statements are in $\mbz_p$, ie modulo $p$.)

\subsection{Upper Bound for $3 \times 3$ Upper Triangular Matrices}
\label{sec-p1:cutoff:3}

\renewcommand{\ugr}{\UU}
\renewcommand{\uab}{A}
\renewcommand{\ucom}{\UU^\com}

As in \S\ref{sec-p1:cutoff:prelim}, we use the abbreviations
\(
	\ugr \cq \UU_{p,d},
\)
\(
	\uab \cq \UU_{p,d}^\ab
\)
and
\(
	n \cq \abs{\UU_{p,d}} = p^{d(d-1)/2}.
\)
Here (\S\ref{sec-p1:cutoff:3}), we study only $d = 3$.
In the $3 \times 3$ case, we only have three terms to deal with.
We abbreviate
\[
	\text{a matrix $M \in \UU_{p,3}$}
\Quad{by}
	\rbr{ M_{1,2}, M_{2,3}, M_{3,3} }.
\]
For matrices $M_1, M_2, ... \in \UU_{p,3}$, writing $M_j \cq (a_j, b_j, c_j)$ for each $j$, we have
\[
\begin{gathered}
	\prodt[t]{1} M_s
=
	\rbB{ \sumt[t]{1} a_s, \: \sumt[t]{1} b_s, \: \sumt[t]{1} c_s + f\rbb{ (a_s)_1^t, \: (b_s)_1^t } }
\\
\text{where}\quad
	f\rbb{ (a_j)_1^t, \: (b_j)_1^t }
\cq
	\sumt[t]{s=1} b_s \sumt[s-1]{r=1} a_r,
\end{gathered}
\label{eq-p1:cutoff:3:product}
\nt
\]
Note that the first two terms are `Abelian' (and correspond to the Abelianisation): we can reorder the product $M_1 \cdots M_t$ in any way we desire, and the first two terms are unchanged; also, so is the first part of the third term, but the polynomial $f$ is not.

We have $k$ generators $Z = [Z_1, ..., Z_k]$; write $Z_i \cq (A_i, B_i, C_i)$ for each $i$.
Recall that $W$ is a DRW on $\mbz_+^k$.
Suppose that $N \cq N(t)$ steps are taken.
Write $(\alpha_1, \beta_1, \gamma_1)$, ..., $(\alpha_N, \beta_N, \gamma_N)$ for the steps taken by $S$.
Write $G_m$ for the generator index chosen at step $m \in [N]$, ie $G_m = i$ if $(\alpha_m, \beta_m, \gamma_m) = (A_i, B_i, C_i)$.
Write $\bm\alpha \cq (\alpha_m)_1^N$ and $\bm\beta \cq (\beta_m)_1^N$.
Let $S'$ be an independent copy of $S$, and make similar definitions.
From \cref{eq-p1:cutoff:3:product},~we~have
\[
	S(t)
=
	\rbb{ \sumt[k]1 A_i W_i(t), \: \sumt[k]1 B_i W_i(t), \: \sumt[k]1 C_i W_i(t) + f\rbr{ \bm\alpha, \bm\beta } }.
\label{eq-p1:cutoff:3:S}
\nt
\]

Recall that we write ``$\equiv$'' to mean ``equivalent modulo $p$''.

\begin{Proof}[Proof of \cref{res-p1:cutoff:res} \mdseries\unboldmath (when $d = 3$)]
First, we claim that
\[
	\pr{ S = S' \mid W \ne W', \: \typ }
=
	1/n
=
	1/p^3.
\label{eq-p1:cutoff:3:W=/W'}
\nt
\]
Indeed, for any $v \in \mbz_p^k \setminus \bra{0}$, by \cref{res-p1:a.X-unif}, each of $\sum_1^k A_i v_i$, $\sum_1^k B_i v_i$ and $\sum_1^k C_i v_i$ is an independent $\Unif(\mbz_p)$;
also, $f(\bm \alpha, \bm \beta)$ is independent of $\sum_1^k C_i W_i(t)$.
Note also that, by typicality, $\abs{W_i - W'_i} < p$ for all $i$, and so $W_i \equiv W'_i$ mod $p$ if and only if $W_i = W'_i$.
Hence conditioning on $W$ and $W'$ and then using \cref{eq-p1:cutoff:3:S} establishes the claimed \cref{eq-p1:cutoff:3:W=/W'}.
Next, recall from \cref{res-p1:cutoff:prelim:W=W'}
that
\[
	\pr{ W = W' \mid \typ }
\le
	e^{-h} / \pr{\typ}
=
	e^{-\omega} e^{-h_0} / \pr{\typ}.
\label{eq-p1:cutoff:3:W=W'}
\nt
\]

It remains to consider the case that $W(t) = W'(t) = w$, for some $w \in \mcw$.
In particular, $S$ and $S'$ take the same number of steps: $N = N'$.
Note, by \cref{eq-p1:cutoff:3:S}, that $S-S' = (0, 0, f(\bm\alpha,\bm\beta)-f(\bm\alpha',\bm\beta'))$.

Expanding the definition of $f$ in \cref{eq-p1:cutoff:3:product}, we may write
\[
	f\rbr{ \bm\alpha, \bm\beta }
=
	\sumt[k]{i,j=1}
	C_{i,j} A_i B_j,
\label{eq-p1:cutoff:3:f-Cij}
\nt
\]
for appropriate $\bra{C_{i,j}}_{i,j=1}^k$;
specifically, for $i,j \in [k]$, we have
\[
	C_{i,j}
\cq
	\sumt[N]{\ell=1}
	\oneb{ G_\ell = j }
	\sumt[\ell-1]{m=1}
	\oneb{ G_m = i };
\Quad{write}
	\bm C \cq \rbr{ C_{i,j} \mid i,j \in [k] }.
\label{eq-p1:cutoff:3:Cij-def}
\nt
\]
Define $C'_{i,j}$ and $\bm C'$ analogously with respect to $W'$.
The body of the proof will be controlling the probability that $\bm C \equiv \bm C'$ conditional on $W(t) = W'(t) = w$, for some typical $w \in \mcw$.
Write
\[
	\mce
\cq
	\brb{ \bm C \equiv \bm C' }
=
	\brb{ C_{i,j} \equiv C'_{i,j} \: \forall \, i,j \in [k] }.
\label{eq-p1:cutoff:3:E-def}
\nt
\]

Trivially,
\(
	\pr{ S = S' \mid W = W' = w, \: \mce }
\le
	1.
\)
We now argue that
\[
	\pr{ S = S' \mid W = W' = w, \: \mce^c }
\le
	2/p.
\label{eq-p1:cutoff:3:Ec}
\nt
\]
Write $D_{i,j} \cq C_{i,j} - C'_{i,j}$.
On the event $\mce^c$, there exist $i',j' \in [k]$ with $D_{i',j'} \not\equiv 0$.
Then
\[
	f(\bm\alpha,\bm\beta) - f(\bm\alpha',\bm\beta')
=
	A_{i'} \rbb{ D_{i',j'} B_{j'} + \sumt{j \ne j'} D_{i,j} B_j }
+	\sumt{i \ne i'} A_i \sumt{j} D_{i,j} B_j.
\label{eq-p1:cutoff:3:f-f'}
\nt
\]
We can write this final expression (with the natural association) as
\[
	U \rbr{ V + X } + Y.
\label{eq-p1:cutoff:3:U(V+X)+Y:def}
\nt
\]
Since $D_{i',j'} \not\equiv 0$ (by choice of $i'$ and $j'$) and $p$ is prime, $U,V \sim^\iid \Unif(\mbz_p)$.
Moreover, $U$ is jointly independent of $X$ and $Y$ and $V$ is independent of $X$ (but not of $Y$); hence $V+X \sim \Unif(\mbz_p)$, independent of $U$, and so $U(V+X) \sim \Unif(\mbz_p)$ and is independent of $Y$ on the event $\bra{V+X \not\equiv 0}$.
(These independence statements are all conditional on $W = W' = w$.)
Thus
\[
	\pr{ U(V+X) + Y \equiv 0 }
\le
	\maxt{u} \pr{ U \equiv u }
+	\pr{ V+X \equiv 0 }
=
	2/p.
\label{eq-p1:cutoff:3:U(V+X)+Y:pr}
\nt
\]
This establishes \cref{eq-p1:cutoff:3:Ec}.

Combining these results,
recalling that $\mce = \bra{ \bm C \equiv \bm C' }$,
writing
\[
	q(t)
\cq
	\MAX{w \in \mcw} \,
	\pr{ \mce \mid W = W' = w },
\]
recalling that $w$ is an arbitrary (fixed) element of $\mcw$,
we find that
\[
	\pr{ S = S' \mid W = W' = w, \: \typ }
\le
	2/p + q(t).
\label{eq-p1:cutoff:3:decomp}
\nt
\]
Once we average over $w \in \mcw$, recalling \cref{eq-p1:cutoff:3:W=W'}, we obtain
\[
	\pr{ S = S', \: W = W' \mid \typ }
\le
	2 e^{-h} \rbb{ 1/p + q(t) } / \pr{\typ}.
\label{eq-p1:cutoff:3:final}
\nt
\]
It remains to make an appropriate definition of typicality, ie of $\mcw$: we require
	that it satisfies \cref{eq-p1:cutoff:prelim:typ:subset},
	that $\pr{W \in \mcw} = 1 - \oh1$, and hence $\pr{\typ} = 1 - \oh1$,
and
	that $e^{-h} \rbr{2/p + q(t)} = \oh{1/n}$.
This is done in \cref{res-p1:cutoff:3:prod-bound} below, recalling that $h = h_0 + \omega$ and $\omega \gg 1$; it is the main technical part of the proof.
Once this is done, combining \cref{eq-p1:cutoff:3:W=/W',eq-p1:cutoff:3:final} gives
\[
	n \, \pr{ S = S' \mid \typ } - 1 = \oh1.
\]
The upper bound in \cref{res-p1:cutoff:res} then follows from \cref{res-p1:cutoff:prelim:D-equiv}, modulo \cref{res-p1:cutoff:3:prod-bound}.
\end{Proof}

It remains to appropriately upper bound $q(t)$ so that the right-hand side of \cref{eq-p1:cutoff:3:decomp} is $\oh{e^h/n}$.

\begin{lem}
\label{res-p1:cutoff:3:prod-bound}
	Suppose that $1 \ll \log k \ll \log n$.
	There exists a $\mcw \subseteq \mbz_+^k$, satisfying \cref{eq-p1:cutoff:prelim:typ:subset}, so that
	\[
		\pr{ W \in \mcw } = 1 - \oh1
	\Qand
		n e^{-h_0} \rbb{ 1/p + q(t) } \lesssim 1.
	\]
\end{lem}

For this proof, let $\eps > 0$, and assume that it is as small as required (but independent of $n$).
Recall that here $d = 3$, so $\log \abs \uab \asymp \log n$ and $1 + \tfrac2{d-2} = 3$.
Hence there are three main regimes:\par\smallskip
{\centering
	$k \ll \log \abs \uab$,\quad%
	$\log \abs \uab \lesssim k \le (\log \abs \uab)^3$%
\QUAD{and}%
	$k \ge (\log \abs \uab)^3$.\par\smallskip}

\begin{Proof}[Proof of \cref{res-p1:cutoff:3:prod-bound} when $k \ge (\log \abs \uab)^3$]
We have $t \ge (1 + 3\eps) \tdiam = (1 + 3\eps) \log_k \abs G$.
Recall from \cref{def-p1:cutoff:prelim:h0} that, in this regime, we take
\(
	h_0 \cq \rbr{1 - \tfrac1\rho} \log n.
\)
Hence $e^{-h_0} = n^{-1+1/\rho}$.

Since $t \ll k$, almost all the generators are picked at most once whp.
For $w \in \mbz_+^k$, define
\[
	\mcj(w)
\cq
	\brb{ i \in [k] \mid w_i = 1 }
\Qand
	J(w) \cq \abs{\mcj(w)}.
\label{eq-p1:cutoff:3:>:pick-once}
\nt
\]
Using this, we make precise our definition of typicality:
\[
	\mcw
\cq
	\brb{ w \in \mbz_+^k \midb \mu_t(w) \le e^{-h}, \: \abs{ J(w) - t e^{-t/k}} \le \tfrac12 \eps t e^{-t/k}, \: \maxt{i} w_i < \tfrac12 p },
\label{eq-p1:cutoff:3:>:typ}
\nt
\]
satisfying \cref{eq-p1:cutoff:prelim:typ:subset}.
Using the conditions of \cref{hyp-p1:cutoff}, we have
\(
	\log_k n \ll \sqrt p
\)
and
\(
	1 \ll t_0 \lesssim k \ll \sqrt p.
\)
Thus the condition $\bra{\maxt{i} w_i < \tfrac12 p}$ holds with probability $1 - \oh1$.
By Binomial concentration and \cref{res-p1:cutoff:prelim:typ}, we then have
\(
	\pr{W \in \mcw} = 1 - \oh1.
\)

Let $w \in \mcw$.
We now argue that
\[
	\pr{\mce \midb W = W' = w, \: \abs{\mcj(w)} = J} \le 1/J!.
\label{eq-p1:cutoff:3:1/J!}
\nt
\]
This holds since, conditional on $W = W' = w$, \emph{different} (relative) orderings, between $S$ and $S'$, of the coordinates chosen once, ie in $\mcj(w)$, must result in some pair $(i,j)$ such that $C_{i,j} = 1$ and $C'_{i,j} = 0$.
There are $J!$ different orderings.

Applying \cref{eq-p1:cutoff:3:1/J!},
using the condition $\abs{J(w) - t e^{-t/k}} \le \tfrac12 \eps t e^{-t/k}$ for $w \in \mcw$,
gives
\[
	q(t)
\le
	1/\rbb{(1-\eps)t}!
\label{eq-p1:cutoff:3:1/t!}
\nt
\]
Note that
	$k = \rbr{ \log \abs \ugr }^{\rho/(1+\nu)}$
and so
	$\tdiam = \log_k n = \tfrac1\rho (1 + \nu) \log n / \log \log n$.
Since $t \ge (1+3\eps) \tdiam$,
direct calculation using \cref{eq-p1:cutoff:3:1/t!} and Stirling's approximation gives
\[
	q\rbr{t}
\le
	\rbb{(1+\eps)\tdiam/e}^{-(1+\eps)\tdiam}
\le
	n^{-1/\rho}.
\label{eq-p1:cutoff:3:q:>}
\nt
\]
Recalling that $e^{-h_0} = n^{-1+1/\rho}$, the proof is completed in the regime $k \ge (\log \abs \uab)^3$:
\[
	n e^{-h_0} \rbb{1/p + q(t)}
\le
	2 n \cdot  n^{-1+(1+\nu)/\rho} \cdot n^{-(1+\nu)/\rho}
=
	2.
\qedhere
\]
\end{Proof}

\begin{Proof}[Proof of \cref{res-p1:cutoff:3:prod-bound} when $\log \abs \uab \lesssim k \le (\log \abs \uab)^3$]
We have $t \ge (1 + 3\eps) t_0$.
Recall from \cref{def-p1:cutoff:prelim:h0} that, in this regime, we take
\(
	h_0 \cq \log \abs \uab.
\)
Hence $e^{-h_0} = \abs \uab^{-1}$.

Since $d \asymp 1$, we have $1 \ll t \lesssim k$.
We use the same definition of typicality here as for $k \ge (\log \abs \uab)^3$.
Since $1 \ll t \lesssim k$, we have $\pr{W \in \mcw} = 1 - \oh1$.

\begin{subequations}
	\label{eq-p1:cutoff:3:q:=:ab}
Since $t \ge (1+3\eps) t_0$, direct calculation using \cref{eq-p1:ent:t0:>,eq-p1:cutoff:3:1/t!} and Stirling's approximation gives
\[
	q\rbb{(1+3\eps)t_0}
\le
	\rbb{(1+\eps)t_0/e}^{-(1+\eps)t_0}
\le
	\abs \uab^{1/(\rho-1)}
\Qwhere
	k \gg \log \abs \uab.
\label{eq-p1:cutoff:3:q:=:>>}
\nt
\]
For $k \eqsim \lambda \log \abs \uab$, with $\lambda \in (0,\infty)$,
we have
\(
	t_0 \eqsim f(\lambda) k \eqsim \lambda f(\lambda) \log \abs \uab
\)
by \cref{eq-p1:cutoff:t*:=}, and thus
\[
	\ex{\abs\mcj} \eqsim \lambda f(\lambda) e^{-f(\lambda)} \log \abs \uab.
\]
Applying \cref{eq-p1:cutoff:3:1/J!},
using the condition $\abs{J(w) - t e^{-t/k}} \le \tfrac12 \eps t e^{-t/k}$ for $w \in \mcw$,
gives
\[
	q(t)
\le
	1/\rbb{ (1 - \eps) \lambda f(\lambda) e^{-f(\lambda)} \log \abs \uab }!;
\]
applying Stirling's approximation, it is easy to see that
this decays super-polynomially, ie
\[
	\log\rbb{ 1 / q\rbb{(1+3\eps)t_0} } \gg \log \abs \uab
\Qwhere
	k \eqsim \lambda \log \abs \uab,
\label{eq-p1:cutoff:3:q:=:=}
\nt
\]
provided $\eps$ is sufficiently small.
\end{subequations}
Hence
\[
	q(t) \le \abs \uab^{-1/(\rho-1)}
\Qwhere
	\log n \lesssim k \le (\log \abs \uab)^3.
\label{eq-p1:cutoff:3:q:=}
\nt
\]
Recall that we want to compare $q(t)$ with $1/p = \abs \uab / \abs \ugr$.
Recall that $\abs \ugr = \abs \uab^{2/d}$.
Some simple algebra then shows that
$q(t) \le 1/p$ when $\rho \le 1 + \tfrac2{d-2} = 3$.
Recalling that $e^{-h_0} =\abs \uab^{-1}$, the proof is completed in the regime $\log \abs \uab \lesssim k \le (\log \abs \uab)^3$:
\[
	n e^{-h_0} \rbb{ 2/p + q(t) }
\le
	\abs \ugr \cdot \abs \uab^{-1} \cdot 3\abs \uab / \abs \ugr
=
	3.
\qedhere
\]
\end{Proof}

\begin{Proof}[Proof of \cref{res-p1:cutoff:3:prod-bound} when $1 \ll k \ll \log \abs \uab$]
We have $t \ge (1 + 3\eps) t_0$.
Recall from \cref{def-p1:cutoff:prelim:h0} that, in this regime, we take
\(
	h_0 \cq \log \abs \uab.
\)
Hence $e^{-h_0} = \abs \uab^{-1}$.

Since $d \asymp 1$, we have
\(
	t_0 \asymp k \abs \uab^{2/k} = k p^{4/k} \gg k.
\)
Hence the same generator is picked lots of times, and so we need a new approach for calculating $q(t)$.
The expected number of times a generator is picked is $s \cq t/k \gg 1$.
As part of our typicality requirements, we ask that `most' pairs $(2i,2i-1)$, with $i \in \bra{1, ..., \floor{k/2}}$, are picked between $\eta s$ and $\eta^{-1} s$ times, for a small positive constant $\eta$, to be chosen later; for the moment, let $\eta \in (0,1)$.
\begin{subequations}
	\label{eq-p1:cutoff:3:<:pick-gen}
For $w \in \mbz_+^k$, write
\[
	\mcc(w)
\cq
	\brb{ i \in \bra{1, ..., \floor{k/2}} \midb \eta s \le \min\bra{w_{2i},w_{2i-1}} \le \max\bra{w_{2i},w_{2i-1}} \le \eta^{-1} s }.
\label{eq-p1:cutoff:3:<:pick-gen:def}
\nt
\]
Then, for $\eta$ sufficiently small (but still a constant), we have
\[
	\pr{ \abs{\mcc(W)} \ge \tfrac25 k } = 1 - \oh1.
\label{eq-p1:cutoff:3:<:pick-gen:prob}
\nt
\]
\end{subequations}
(We could replace $\tfrac25$ by any constant less than $\tfrac12$, at the cost only of making $\eta$ a smaller constant.)
We use this to make precise our definition of typicality for this regime:
\[
	\mcw
\cq
	\brb{ w \in \mbz_+^k \midb \mu_t(w) \le e^{-h}, \: \abs{\mcc(w)} \ge \tfrac25 k, \: \maxt{i} w_i < \tfrac12 p }.
\label{eq-p1:cutoff:3:<:typ}
\nt
\]
Then, like before and additionally using \cref{eq-p1:cutoff:3:<:pick-gen:prob}, we have
\(
	\pr{W \in \mcw} = 1 - \oh1.
\)
If $i \in \mcc(w)$, then
\(
	\max\bra{C_{2i,2i-1},C'_{2i,2i-1}} \le w_i^2 \lesssim s^2 \ll p
\)
as $d \asymp 1$,
so
\(
	\bra{C_{2i,2i-1} \equiv C'_{2i,2i-1}} = \bra{C_{2i,2i-1} = C'_{2i,2i-1}}.
\)

We claim that it is sufficient to fix an arbitrary $w \in \mcw$ and prove the bound
\[
	\maxt{i} \rho_i \le p^{-3/k}
\Qwhere
	\rho_i
\cq
	\maxt{x}
	\pr{ C_{2i,2i-1} = x \midb W = w } \oneb{ i \in \mcc(w) }.
\label{eq-p1:cutoff:3:<:mode-i}
\nt
\]
To see this, first make the simple observation that, for any $\mci \subseteq \bra{1,...,\floor{k/2}}$, we have
\[
	\brb{ \rbr{C_{i,j}}_{i,j\in[k]} = \rbr{C'_{i,j}}_{i,j\in[k]} }
\subseteq
	\brb{ \rbr{C_{2i,2i-1}}_{i\in\mci} = \rbr{C'_{2i,2i-1}}_{i\in\mci} }.
\label{eq-p1:cutoff:3:<:restrict-i}
\nt
\]
Given $W = W' = w$, the event $\bra{C_{i,j} = C'_{i,j}}$ is determined by the relative order in which the generators $i$ and $j$ are chosen.
Hence, since the pairs $(2i,2i-1)$ are disjoint, the events $\bra{ C_{2i,2i-1} = C'_{2i,2i-1} }$ are independent for different $i$, conditional on $W = W' = w$.
Take $\mci \cq \mcc(w)$, which has size at least $\tfrac25 k$.
By the aforementioned independence, given \cref{eq-p1:cutoff:3:<:mode-i}, we have
\[
	\pr{ \bm C = \bm C' \midb W = W' = w }
\le
	\rbr{ \maxt{i} \rho_i }^{2k/5}
\le
	p^{(-3/k)(2k/5)}
=
	p^{-6/5}
\ll
	1/p,
\label{eq-p1:cutoff:3:<:q<<1/p}
\nt
\]
and hence $q(t) \ll 1/p$.
The proof is then completed as in the regime $\log \abs \uab \lesssim k \le (\log \abs \uab)^3$.

It remains to prove \cref{eq-p1:cutoff:3:<:mode-i}.
For simplicity of notation, we assume that $1 \in \mcc(w)$ and set $i \cq 1$.
	Let $r \cq w_1 + w_2$,
	let $N \cq \sumt[k]{1} w_i$ be the number of steps taken
and
	write $G_\ell \in [k]$ for the generator index chosen in the $\ell$-th step.
We can then write the random word as $S = Z_{G_1} \cdots Z_{G_N}$.
Let $J_1 < \cdots < J_r$ be the (random) indices with $G_{J_\ell} \in \bra{1,2}$.
Now define the vector $I \in \bra{1,2}^r$ by $I_\ell \cq G_{J_\ell}$.
Thus $I$ encodes the relative order between the different occurrences (with multiplicities) of the generators labelled by $\bra{1,2}$ in the word $S$.
By typicality, $2 \eta s \le r \le 2 \eta^{-1} s$.

Let $I'$ be the random vector obtained from $I$ by picking a 2 uniformly at random and omitting it from $I$.
(Eg, if $I = (2,2,1,1,2,1)$ and we pick the last 2, then $I' = (2,2,1,1,1)$.)
Importantly, we are omitting elements of the \emph{relative} order of appearances of $Z_1$ and $Z_2$, not the \emph{absolute} locations of the corresponding generators.

By the definition of $C_{1,2}$, given in \cref{eq-p1:cutoff:3:Cij-def}, given $W = w$, the value of $C_{1,2}$ is a function only of the relative locations $I$.
Hence, given $I'$ also, it is a function only of the location of the omitted 2.
It is constant on the set of locations which give rise to the same $I$:
	two different placements of the omitted 2 give rise to the same $I$ if and only if they both lie in the same (possibly empty) interval of consecutive 2s.
(Eg, if $I' = (2,2,1)$, then there are three locations in which we can insert a 2 to get $I = (2,2,2,1)$, namely the first, second and third positions, and only one to get $I = (2,2,1,2)$, namely the fourth position; the first three give rise to $C_{1,2} = 0$ and the fourth to $C_{1,2} = 1$.)

Hence, writing $L(I')$ for the longest interval of 2s in $I'$, we have
\[
	\maxt{x}
	\pr{ C_{1,2} = x \midb W = w, \: I' }
\le
	\rbb{ L(I') + 1 } / r.
\]
By \cref{res-p1:cutoff:3:na} below, with $m = 2$, we find that $L(I')/(C \log r) \preceq \Geom(\tfrac12)$ given $W = w$ for a sufficiently large constant $C$, and so $\ex{L(I') \mid W = w} \lesssim \log r$.
Hence
\[
	\maxt{x}
	\pr{ C_{1,2} = x \midb W = w }
\le
	\ex{ L(I') + 1 \midb W = w } / r
\lesssim
	(\log r) / r.
\]
Since $2 \eta s \le r \le 2 \eta^{-1} s$, as $w \in \mcw$, and $\eta$ is a (small) constant, this last expression is $\oh{1/s^{3/4}}$. (In fact, it is $\Oh{\log s}/s$.)
Recalling that $s \asymp p^{4/k}$ establishes \cref{eq-p1:cutoff:3:<:mode-i}.
This completes the proof.
\end{Proof}

It remains to state and prove the claim regarding $\ex{L(I')}$.
We actually state and prove a slightly more general claim, that we are then able to use in the analysis of the $d \times d$ matrices.

\begin{clm}
\label{res-p1:cutoff:3:na}
	Let $m \in \mbn$ and $\eta \in (0,1)$.
	Let $\bra{w_1, ..., w_m}$ be arbitrary positive integers satisfying $w_i/w_j \in [\eta^2, \eta^{-2}]$ for all $i, j \in [m]$.
	For each $k \in \bra{1, ..., m}$, let there be $w_k$ balls of colour $k$; write $r \cq \sum_{k=1}^m w_k$ for the total number of balls.
	Choose a uniform permutation of the balls on positions $\bra{1, ..., r}$.
	For each $k \in \bra{1, ..., m}$, let $L_k$ be the longest interval \emph{without} any balls of colour $k$.
	Then, for each $k$, we have the stochastic domination
	\[
		L_k / (\eta^{-2} m \log r) \preceq \Geom(\tfrac12).
	\]
\end{clm}

\begin{Proof}
Without loss of generality, take $k \cq 1$ and write $L \cq L_1$.
By assumption, $w_i/w_j \in [\eta^2, \eta^{-2}]$ for all $i$ and $j$, and $\eta \in (0,1)$ is a constant.
Hence $r \le m \eta^{-2} w_1$, and so $w_1 \ge \eta^2 r/m$.
Let $\ell \in \mbn$ to be chosen shortly; write
\(
	[1,r] \subseteq [1,\ell] \cup [2,\ell+1] \cups [r,r+\ell-1].
\)
By direct calculation,
\[
	\pr{ L > \ell }
&
\le
	r \, \pr{ \text{no 1 in the interval $[1,\ell]$} }
\\&
=
	r \cdot \rbb{1 - \tfrac{w_1-1}{r-1}} \rbb{1 - \tfrac{w_1-1}{r-2}} \cdots \rbb{1 - \tfrac{w_1-1}{r-\ell}}
\\&
\le
	r \rbb{1 - \tfrac{w_1-1}{r-1}}^\ell
\le
	r \expb{ - \ell w_1/r }
\le
	r \expb{ - \ell \eta^2 / m },
\]
where for the penultimate inequality we used the fact that $\tfrac{w_1-1}{r-1} \le \tfrac{w_1}{r}$, which holds since $w_1 \le r$.
Choosing $\ell \cq (k + 1) \eta^{-2} m \log r$ gives
\[
	\pr{ L > (k + 1) \eta^{-2} m \log r }
&
\le
	r \expb{ - (k+1) \log r }
=
	r^{-k}.
\]
Thus we may stochastically dominate
\[
	L / (\eta^{-2} m \log r)
\preceq
	\Geom(1 - 1/r)
\preceq
	\Geom(\tfrac12).
\qedhere
\]
	%
\end{Proof}

%

\subsection{Upper Bound for $d \times d$ Upper Triangular Matrices}
\label{sec-p1:cutoff:d}

\renewcommand{\ugr}{\UU_{p,d}}
\renewcommand{\uab}{\UU_{p,d}^\ab}

The high-level ideas of the proof will be similar to the $d = 3$ case, but there are a number of subtleties which need to be navigated.
Analogously to \cref{res-p1:cutoff:3:prod-bound}, there will be a certain probability that requires bounding, and the argument for bounding this will differ depending on $k$; the specific reference will be Lemmas \ref{res-p1:cutoff:d:prod-bound} and \ref{res-p1:cutoff:prelim:omega-choice}, and comes at the end of the section.

We also use the same preliminaries (see \S\ref{sec-p1:cutoff:prelim}), and in particular consider
\[
	D(t) = \abs \ugr \, \pr{ S = S' \mid \typ } - 1.
\]

The analogues of \cref{eq-p1:cutoff:3:product,eq-p1:cutoff:3:S} are different for general $d$ than for $d = 3$: they have the same basic structure, but with the addition of `higher order' terms (given by $g_{a,b}$ in the lemma below).
The following lemma is for both the DRW and the SRW; take $\sigma_\ell \cq 1$ for all $\ell$ to reduce to the~DRW.

\begin{lem}
\label{res-p1:cutoff:d:matrix-product}
Let $Z_1,...,Z_k \in \ugr$.
Let
	$\gamma \in [k]^L$ and $\sigma \in \bra{\pm1}^L$.
For $i,j \in [k]$,
set
\[
	C_{i,j}(\gamma, \sigma)
\cq
	\sumt[L]{\ell=0}
	\sumt[\ell-1]{m=0}
	\sigma_m \sigma_\ell \one{\gamma_m = i, \: \gamma_\ell = j }
+	\one{i = j}
	\sumt[L]{\ell=0}
	\one{\gamma_\ell = i, \: \sigma_\ell = -1}.
\]
Set
\(
	M \cq Z_{\gamma_1}^{\sigma_1} \cdots Z_{\gamma_L}^{\sigma_L}.
\)
Then,
for all $a \in [d]$,
we have
\[
	M(a,a) = 1
\Qand
	M(a,a+1) = \sumt[L]{\ell=1} \sigma_{\gamma_\ell} Z_{\gamma_\ell}(a,a+1),
\]
and,
for all $a,b \in [d]$ with $b \ge a+2$,
we have
\[
	M(a,b)
=
	\sumt{\ell \in [L]}
	Z_{\gamma_\ell}(a,b)
+	\sumt{i,j \in [k]}
	C_{i,j}(\gamma, \sigma) Z_i(a,a+1) Z_j(a+1,b)
+	g_{a,b}(\gamma, \sigma; Z_1,...,Z_k),
\label{eq-p1:cutoff:d:matrix-product}
\nt
\]
where $g_{a,b}(\gamma, \sigma; Z_1, ..., Z_k)$ is a polynomial in
\(
	\rbr{ Z_i(x,y) : i \in [k], \: x \in [d-1], \: y > x }.
\)
Further, in this polynomial, each monomial contains the term $Z_i(a,a+1)$ either 0 times or exactly once and no monomial contains a term of the form $Z_i(a,a+1)Z_j(a+1,b)$ for $i,j \in [k]$.
\end{lem}

We give a sketch of the argument here; the rigorous details are deferred to \cite[\cref{res-p0:heisenberg:matrix-product}]{HOt:rcg:supp}.

\begin{Proof}[Proof Sketch of \cref{res-p1:cutoff:d:matrix-product}]
The fundamental idea is to write a matrix $M_\ell \in \ugr$ as $I + N_\ell$ where $N_\ell$ is \emph{strictly} upper triangular.
For $M_1, ..., M_L \in \ugr$
written like this,
one then has
\[
	M_1 \cdots M_L
=
	\prodt[L]{\ell=1}
	\rbr{ 1 + N_\ell }
=
	I
+	\sumt[L]{\ell=1}
	N_\ell
+	\sumt{m_1 < m_2}
	N_{m_1} N_{m_2}
+	\sumt[L]{\ell=3}
	\sumt{m_1 < \cdots < m_\ell}
	\prodt[\ell]{r=1}
	N_{m_r},
\]
where the indices $m_r$ run over $[L]$.
Further,
for $(a,b)$ with $b - a \ge 2$,
one can write
\[
	\rbb{ N_{m_1} N_{m_2} }(a,b)
&
=
	\sumt{c \in [1,d]}
	N_{m_1}(a,c) N_{m_2}(c,b)
=
	\sumt{c \in [a+1,b-1]}
	M_{m_1}(a,c) M_{m_2}(c,b)
\\
&=
	M_{m_1}(a,a+1) M_2(a+1,b)
+	\sumt[b-1]{c=a+2}
	M_{m_1}(a,c) M_2(c,b).
\]
We consider this latter sum along with all products of degree at least 3 as `higher order' terms.
Writing $\sumt{m_1 < m_2}$ as $\sumt[L]{m_2=1} \sumt[m_2-1]{m_1=1}$, the formula for $C_{i,j}$ follows for the DRW (ie $\sigma_\ell \cq 1$ for all~$\ell$).

The SRW analysis is similar. Since $N_\ell$ is strictly upper triangular, $N_\ell^d = 0$.
Thus
\[
	M_\ell^{-1}
=
	(1 + N_\ell)^{-1}
=
	I
-	N_\ell
+	N_\ell^2
+	\sumt[d]{t=3}
	(-1)^t N_\ell^t.
\]
Separating out `higher order' terms similarly, we deduce the formula for the SRW.
\end{Proof}

\renewcommand{\ugr}{\UU}
\renewcommand{\uab}{A}

As previously, we use the abbreviations
\(
	\ugr \cq \UU_{p,d},
\)
\(
	\uab \cq \UU_{p,d}^\ab
\)
and
\(
	n \cq \abs{\UU_{p,d}} = p^{d(d-1)/2}.
\)

\begin{Proof}[Proof of \cref{res-p1:cutoff:res} \mdseries \unboldmath (general $d$)]
When $W(t) \ne W'(t)$, the same argument as for $d = 3$, using \cref{res-p1:a.X-unif}, applies, replacing \cref{eq-p1:cutoff:3:S} by \cref{eq-p1:cutoff:d:matrix-product}:
\[
	\pr{ S = S' \mid W \ne W', \: \typ }
=
	1/n
=
	1/p^{d(d-1)/2}
=
	p^{-(d-1)(d-2)/2} \cdot p^{-(d-1)}.
\label{eq-p1:cutoff:d:W=/W'}
\nt
\]
This is the analogue of \cref{eq-p1:cutoff:3:W=/W'}.
Next, recall from \cref{res-p1:cutoff:prelim:W=W'}
that
\[
	\pr{ W = W' \mid \typ }
\le
	e^{-h} / \pr{\typ}
=
	e^{-\omega} e^{-h_0} / \pr{\typ}.
\label{eq-p1:cutoff:d:W=W'}
\nt
\]

\smallskip

Now suppose that $W(t) = W'(t) = w$, where $w$ is some fixed element of $\mcw$ (yet to be defined fully). Then the `Abelian' parts of $S$ and $S'$, corresponding to the first term in the right-hand side of \cref{eq-p1:cutoff:d:matrix-product}, cancel (as was the case when $d = 3$).
Write $\bm C \cq (C_{i,j})$ and $\bm C' \cq (C'_{i,j})$ for the $C(\gamma)$ in \cref{res-p1:cutoff:d:matrix-product} generated by $S$ and $S'$, respectively.
Write $\mce \cq \bra{ \bm C = \bm C' }$.
On $\mce$, the middle terms of \cref{eq-p1:cutoff:d:matrix-product} cancel, leaving only the higher-order terms; upper bound $\pr{ S = S' \mid W = W', \: \mce } \le 1$.

Now suppose that $\mce$ does not hold; choose, and fix, $(i',j')$ so that $C_{i',j'} \not\equiv C'_{i',j'}$.
By the condition \cref{eq-p1:cutoff:prelim:typ:subset} which $\mcw$ must satisfy and the definition of $C_{i,j}$, this implies that $C_{i',j'} \neq C'_{i',j'}$.
Analogously to \cref{eq-p1:cutoff:3:f-f',eq-p1:cutoff:3:U(V+X)+Y:def}, where $d$ was equal to 3, letting
\[
	U_{a,b} \cq Z_{i'}(a,a+1)
\Qand
	V_{a,b} \cq \rbr{ C_{i',j'} - C'_{i',j'} } Z_{j'}(a+1,b),
\label{eq-p1:cutoff:d:UV:def}
\nt
\]
we can, for some random variables $X_{a,b}$ and $Y_{a,b}$, write
\[
	\sumt[k]{i,j}
	\rbr{ C_{i,j} - C'_{i,j} } Z_i(a,a+1) Z_j(a+1,b)
\Quad{naturally as}
	U_{a,b} (V_{a,b} + X_{a,b}) + Y_{a,b}.
\]
For the moment, fix $(a,b)$.
Analogously to the $d = 3$ case, ie \cref{eq-p1:cutoff:3:Ec,eq-p1:cutoff:3:f-f',eq-p1:cutoff:3:U(V+X)+Y:def,eq-p1:cutoff:3:U(V+X)+Y:pr},
the following hold:
	$U_{a,b},V_{a,b} \sim \Unif(\mbz_p)$;
	$U_{a,b}$ is independent of $(V_{a,b}, X_{a,b}, Y_{a,b})$;
	$V_{a,b}$ is independent of $X_{a,b}$ (but not of $Y_{a,b}$).
Thus $U_{a,b} (V_{a,b} + X_{a,b}) \sim \Unif(\mbz_p)$ is independent of $Y$ on the event $V_{a,b} + X_{a,b} \not\equiv 0$.
Hence
\[
	\maxt{r} \pr{ U_{a,b}(V_{a,b}+X_{a,b}) + Y_{a,b} \equiv r }
\le
	\maxt{u} \pr{ U_{a,b} \equiv u } + \pr{ V_{a,b} + X_{a,b} \equiv 0 }
\le
	2/p;
\label{eq-p1:cutoff:d:U(V+X)+Y-pr}
\nt
\]
Now compare $S_{a,b}$ and $S'_{a,b}$.
	Since $W = W'$, the `Abelian' part cancels;
	we are left with the $U_{a,b}(V_{a,b}+X_{a,b}) + Y_{a,b}$ part and the higher-order terms, given by the $g_{a,b}$ polynomials in \cref{eq-p1:cutoff:d:matrix-product}.
	These two parts are independent, by the conditions of \cref{res-p1:cutoff:d:matrix-product}.
Hence \cref{eq-p1:cutoff:d:U(V+X)+Y-pr} implies that
\[
	\prb{ S_{a,b} = S'_{a,b} \midb W = W' = w, \: \mce^c }
\le
	2/p.
\label{eq-p1:cutoff:d:Sab=S'ab}
\nt
\]

Now, the random variables $\bra{X_{a,b}, Y_{a,b}}_{a,b}$ are not independent.
Also, $U_{a,b} = Z_{i'}(a,a+1)$ does not depend on $b$, and so $\bra{ U_{a,b}, V_{a,b} \mid b \ge a+2 }_{a,b}$ are not independent either.
However, if we fix $b$ then $\bra{ U_{a,b}, V_{a,b} \mid b \ge a+2 }_{a}$ is a collection of independent variables.
We exploit this.

Partition the $[k]$ generators into $d-2$ sets $(P_3, ..., P_d)$.
For each (fixed) $b \in \bra{3,...,d}$, we use generators only from $P_b$; this will give independence when we consider all $b$.
(Note that for $b \in \bra{1,2}$ there are no terms above the super-diagonal.)
Then for the $(a,b)$-th coordinate we try to get $C_{i',j'} \not\equiv C'_{i',j'}$ for some $(i',j')$ with $i', j' \in P_b$.
Now, for each $b$,
using this pair $(i',j')$ in the definition \cref{eq-p1:cutoff:d:UV:def} of $U_{a,b}$ and $V_{a,b}$,
the random variables $\bra{U_{a,b}, V_{a,b} \mid b \ge a+2}_{a}$ are independent, since they depend on a disjoint set of generators.

For each $b \in \bra{3,...,d}$, write
\[
	\bm C_b \cq \rbr{C_{i,j}}_{i,j \in P_b},
\quad
	\bm C'_b \cq \rbr{C'_{i,j}}_{i,j \in P_b}
\Qand
	\mce_b
\cq
	\bra{ \bm C_b = \bm C'_b }.
\label{eq-p1:cutoff:d:CbEb}
\nt
\]
We wish to get an analogue of \cref{eq-p1:cutoff:3:decomp}, for general $d$.
For $w \in \mcw$ and $b \in \bra{3, ..., d}$,
write
\[
	\prt[w]{\cdot} \cq \pr{ \, \cdot \mid W = W' = w}
\Qand
	S_{:,b} \cq \rbr{ S_{a,b} \mid a = 1,...,b-2 },
\]
ie the $b$-th column above the super-diagonal.
Herein, in $\sum_3^d$ and $\prod_3^d$, the implicit index is $b$.
Then
\[
	\pr[w]{ S = S' }
=
	\prodt[d]{3}
	\pr[w]{ S_{:,b} = S'_{:,b} \midb S_{:,b'} = S'_{:,b'} \: \forall \, b' = 3,...,b-1 }
\]
Using \cref{eq-p1:cutoff:d:Sab=S'ab}, and noting that $S_{:,b}$ has $b-2$ entries, we obtain
\[
	\pr[w]{ S_{:,b} = S'_{:,b} \midb S_{:,b'} = S'_{:,b'} \: \forall \, b' = 3,...,b-1 }
\le
	(2/p)^{b-2} + \pr[w]{ \mce_b };
\]
this uses the aforementioned independence between columns, guaranteed by the partitioning of the generators.
Combining these two equations, we obtain
\[
	\pr[w]{ S = S' }
\le
	2^{d^2/2} \prodt[d]{3}
	\rbb{ 1/p^{b-2} + q_b(t) }
\Qwhere
	q_b(t)
\cq
	\MAX{w \in \mcw} \, \prodt[d]{3} \prt[w]{\mce_b}.
\label{eq-p1:cutoff:d:decomp}
\nt
\]

It remains to make an appropriate definition of typicality, ie of $\mcw$, and choose the partition $(P_3, ..., P_d)$ appropriately.
For reasons explained later, we end up choosing $P_b$ so that $R_b \cq \abs{P_b}/k = (b-2)/\binom{d-1}2$, omitting floor/ceiling signs.
(Note that $\sumt[d]3 R_b = 1$, as required.)
We justify the omission of floor/ceiling signs by the fact that $\abs{P_b} \asymp (b-2) kd^{-2} \gg 1$ (as $d^2 \ll k$).

This is all done in \cref{res-p1:cutoff:d:prod-bound}, which gives the following bound:
\[
	n \, \pr[w]{ S = S' }
\equiv
	n \, \pr{ S = S' \mid W = W' = w }
\le
	e^{h_0} 2^{d^2}.
\label{eq-p1:cutoff:d:S=S'|W=W'}
\nt
\]
Combined with \cref{eq-p1:cutoff:d:W=/W',eq-p1:cutoff:d:W=W'} this implies that
\[
	n \, \pr{ S = S' \mid \typ } - 1
\le
	2 \cdot e^{-\omega} 2^{d^2},
\label{eq-p1:cutoff:d:L2:final}
\nt
\]
where we shall choose $\typ$ so that $\pr{\typ} = 1 - \oh1$.
If we can show that we can choose $\omega \gg d^2$, then the upper bound in \cref{res-p1:cutoff:res} then follows from \cref{res-p1:cutoff:prelim:D-equiv}, modulo \cref{res-p1:cutoff:d:prod-bound}.

It remains to prove that we can choose $\omega \gg d^2$.
\cref{res-p1:cutoff:prelim:omega-choice} says that we can choose any $\omega \ll \min\bra{k, \: \log \abs \uab}$.
\cref{hyp-p1:cutoff} implies that
\(
	d^2 \ll \min\bra{k, \: \log \abs \uab}
\)
is satisfied, as required:
\[
	d^3 \ll k \text{ when } k \ll \log \abs \uab
\Quad{and}
	d \ll \log \log p \text{ when } k \gtrsim \log \abs \uab.
\qedhere
\]
\end{Proof}

It remains to appropriately bound $q_b(t)$, defined in \cref{eq-p1:cutoff:d:decomp}.

\begin{lem}
\label{res-p1:cutoff:d:prod-bound}
	Let $\eps > 0$ and set $t \cq (1 + 3\eps) t_*$.
	Assume the conditions of \cref{hyp-p1:cutoff}.
	Then there exists a $\mcw \subseteq \mbz_+^k$, satisfying \cref{eq-p1:cutoff:prelim:typ:subset}, so that
	\[
		\pr{W \in \mcw} = 1 - \oh1
	\Quad{and}
		\abs \ugr e^{-h_0} \prodt[d]{3} \rbb{ 1/p^{b-2} + q_b(t) }
	\lesssim
		2^{d^2/2}.
	\]
\end{lem}

Recall the condition on $\mcw$ given by \cref{eq-p1:cutoff:prelim:typ:subset}.
Since $t \ge (1 + 3\eps)t_0$, by \cref{res-p1:cutoff:prelim:typ}, this condition is satisfied with probability $1 - \oh1$.
Hence we need only check that any additional constraints are also satisfied with probability $1 - \oh1$.
Recall that we use the notation
\[
	k = \rbb{ \log \abs \uab }{}^\rho
\Quad{and}
	\tfrac12 d = \rbb{ \log \abs \ugr }{}^\nu,
\Quad{so}
	k = \rbb{ \log \abs \ugr }{}^{\rho/(1+\nu)},
\Quad{and}
	n = \abs \ugr = p^{d(d-1)/2}.
\]

%

\begin{Proof}[Proof of \cref{res-p1:cutoff:d:prod-bound} for $k \gtrsim \log \abs \uab$]
Let $\eps > 0$ and set $t \cq (1 + 3\eps) t_*$;
write $s \cq t/k$.

\medskip

\emph{Typicality.}
As when $d = 3$, when $k \gg \log \abs \uab$ almost all the generators are picked at most once; when $k \asymp \log \abs \uab$, a constant proportion are.
As part of our typicality requirement ($\typ$), we ask that at least $(1-\eps) t e^{-t/k}$ generators are picked exactly once---ie at least $(1-\eps)$ times the expected number.
Given this, we can then choose our partition so that, for each $b \in \bra{3,...,d}$, writing $R_b \cq \abs{P_b}/k$, at least $(1-\eps) t e^{-t/k} R_b$ generators from $P_b$ are picked exactly once.

We can hence use the same definition of typicality, for $k \gtrsim \log \abs \uab$, as when $d = 3$:
\[
	\mcw
\cq
	\brb{
		w \in \mbz_+
	\midb
		\mu_t(w) \le e^{-h}, \:
		\abs{J(w) - s e^{-s} k} \le \tfrac12 \eps s e^{-s} k, \:
		\maxt{i} w_i < \tfrac12 p
	},
\label{eq-p1:cutoff:d:typ:>=}
\nt
\]
satisfying \cref{eq-p1:cutoff:prelim:typ:subset},
recalling that $J(w) = \sumt[k]{1} \one{w_i = 1}$.
As previously, we have
\(
	\pr{W \in \mcw} = 1 - \oh1.
\)

Analogously to \cref{eq-p1:cutoff:3:1/t!}, when $k \gg \log \abs \uab$, we have
\[
	q_b(t)
\le
	1/\rbb{(1 - 2\eps)t R_b}!,
\label{eq-p1:cutoff:d:1/t!}
\nt
\]
where we have absorbed the $e^{-t/k} = 1 - \oh1$ term into the $(1-2\eps)$; we consider $k \asymp \log \abs \uab$ later.

	%

\medskip

\emph{Regime $k \ge (\log \abs \uab)^{1 + 2/(d-2)}$.}
We have
\(
	t \ge (1 + 3 \eps) \tdiam = (1 + 3 \eps) \log_k \abs \uab.
\)
Direct calculation, analogous to \cref{eq-p1:cutoff:3:q:>}, using \cref{eq-p1:cutoff:t*:>>} and \cref{eq-p1:cutoff:d:1/t!} and Stirling's approximation gives
\[
	q_b(t)
\le
	1/\rbb{ (1 - 2 \eps) \cdot (1 + 3 \eps) \tdiam \cdot R_b }!
\le
	n^{ - R_b 1/\rho }.
\]

In \cref{eq-p1:cutoff:3:q:>}, we upper bounded $q(t) \le n^{-1/\rho}$, and this term was dominant in the sum $1/p + n^{-1/\rho}$.
Here, we compare $q_b(t) \le n^{-R_b/\rho}$ and $1/p^{b-2}$.
It is thus natural to choose $R_b \propto b-2$, ie $R_b \cq (b-2) / \binom{d-1}2$, for $b \in \bra{3,...,d}$.
Observe that
\[
	1/p^{b-2} \le n^{-R_b/\rho}
\Quad{if and only if}
	\rho (b-2) / \binomt d2 \ge R_b = (b-2)/\binomt{d-1}2;
\]
hence we need
\(
	\rho \ge \binom d2 / \binom{d-1}2 = 1 + \tfrac2{d-2},
\)
which is precisely the regime which we are considering.

Combining the upper bounds just developed, we deduce that
\[
	\prodt[d]{3} \rbb{ 1/p^{b-2} + q_b(t) }
\le
	2^d \prodt[d]{3} q_b(t)
\le
	2^d n^{-1/\rho},
\]
since $\sumt[d]{3} R_b = 1$.
Recalling that $h_0 = (1 - \tfrac1\rho) \log n$ in this regime, we deduce the desired bound.

\medskip

\emph{Regime $\log \abs \uab \lesssim k \ll (\log \abs \uab)^{1 + 2/(d-2)}$.}
We have $t \ge (1 + 3\eps) t_0$.
Recall from \cref{def-p1:cutoff:prelim:h0} that, in this regime, we take $h_0 \cq \log \abs \uab$.
Hence $e^{-h_0} = \abs \uab^{-1}$.
We subdivide the regime.

\smallskip

\textit{Consider first $k \gg \log \abs \uab$.}
Direct calculation, analogous to \cref{eq-p1:cutoff:3:q:=:>>}, using \cref{eq-p1:cutoff:t*:>} and \cref{eq-p1:cutoff:d:1/t!} and Stirling's approximation, gives
\begin{subequations}
	\label{eq-p1:cutoff:d:q:>=}
\[
	q_b(t)
\le
	1/\rbb{ (1-2\eps) \cdot (1 + 3\eps) t_0 \cdot R_b }!
\le
	\expb{ - \tfrac2d \tfrac1{\rho-1} R_b \log n };
\label{eq-p1:cutoff:d:q:>=:>}
\nt
\]
again, this crucially uses the fact that $d = \logn[\oh1]$ and $\eps > 0$ is a constant.

In \cref{eq-p1:cutoff:3:q:=:>>}, we upper bounded $q(t) \le \abs \uab^{-1/(\rho-1)}$, and this term was subdominant in the sum $1/p + \abs \uab^{1/(\rho-1)}$.
Here, we compare $q_b(t) \le \abs \uab^{R_b/(\rho-1)}$ with $1/p^{b-2}$.
Again, it is thus natural to choose $R_b \propto b-2$, ie $R_b \cq (b-2) / \binom{d-1}2$, for $b \in \bra{3,...,d}$.
Observe that
\[
	1/p^{b-2} \ge \expb{ - \tfrac2d \tfrac1{\rho-1} R_b \log n }
\Quad{if and only if}
	(\rho - 1) (b-2) / (d-1) \le R_b = (b-2) / \binomt{d-1}2;
\]
hence we need
\(
	\rho \le 1 + \tfrac2{d-2},
\)
which is precisely the regime that we are considering.

Combining the upper bounds just developed, we deduce that
\[
	\prodt[d]{3} \rbb{ 1/p^{b-2} + q_b(t) }
\le
	2^d \prodt[d]{3} 1/p^{b-2}
\le
	2^d p^{-\binom{d-1}2}
=
	2^d \abs \uab / \abs \ugr.
\]
Recalling that $h_0 = \log \abs \uab$ in this regime, we deduce the desired bound.

\smallskip

\textit{Consider now $k \asymp \log \abs \uab$.}
Suppose that $k \eqsim \lambda \log \abs \uab$ with $\lambda \in (0,\infty)$.
Direct calculation, analogous to \cref{eq-p1:cutoff:3:q:=:=}, using \cref{eq-p1:cutoff:t*:=} and \cref{eq-p1:cutoff:d:1/t!} and Stirling's approximation gives
\[
	\log\rbb{ 1/q_b(t) }
\asymp
	d^{-1} \log \log \abs \uab \cdot (b-2) \log p
\gg
	(b-2) \log p,
\label{eq-p1:cutoff:d:q:>=:=}
\nt
\]
\end{subequations}
using the conditions on $d$.
The proof is then completed in exactly the same way as above.
\end{Proof}

\begin{Proof}[Proof of \cref{res-p1:cutoff:d:prod-bound} for $k \ll \log \abs \uab$]
Set $t \cq (1 + 3\eps) t_* \ge (1 + 3\eps) t_0$.
Recall from \cref{def-p1:cutoff:prelim:h0} that, in this regime, we take
\(
	h_0 \cq \log \abs \uab.
\)
Hence $e^{-h_0} = \abs \uab^{-1}$.
Then $s \cq t/k \asymp \abs \uab^{2/k} \gg 1$, by \cref{res-p1:ent:t0} and the assumption $k \ll \log \abs \uab$.

As noted in the $d = 3$ case, neither the actual value of $t$ nor the fact that $W$ and $W'$ are independent DRWs is of much consequence.
Even the particular form of $s$ is not important: it can be changed, subject to changing the conditions on $d$ appropriately.

In the case $d = 3$, we looked at (adjacent) pairs of indices $(2i,2i-1)$.
For general $d$, this is insufficient; instead, we look at $m$-tuples, where $m$ is a (growing) function of $d$.

In this regime, the same generator is picked lots of times, with expectation $s = t/k \gg 1$.
For the moment, let $\eta \in (0,1)$.
For $w \in \mbz_+^k$, write
\begin{subequations}
	\label{eq-p1:cutoff:d:pick-gen}
\[
	\mcc(w)
\cq
	\brb{ i \in [k] \midb \eta s \le w_i \le \eta^{-1} s }.
\label{eq-p1:cutoff:d:pick-gen:def}
\nt
\]
Then, for $\eta$ sufficiently small (but still a constant), we have
\[
	\pr{ \abs{\mcc(W)}/k \ge \tfrac45 } = 1 - \oh1.
\label{eq-p1:cutoff:d:pick-gen:prob}
\nt
\]
\end{subequations}
(We could replace $\tfrac45$ by any constant less than 1.)
This will form part of our typicality requirements:
\[
	\mcw
\cq
	\brb{ w \in \mbz_+^k \midb \mu_t(w) \le e^{-h}, \: \abs{\mcc(w)} \ge \tfrac45 k, \: \maxt{i} w_i < \tfrac12 p }.
\label{eq-p1:cutoff:d:typ:<}
\nt
\]
Then, like before and additionally using \cref{eq-p1:cutoff:3:<:pick-gen:prob}, we have
\(
	\pr{W \in \mcw} = 1 - \oh1.
\)
For $i,j \in \mcc(w)$,
we have
\(
	\bra{C_{i,j} \equiv C'_{i,j}} = \bra{C_{i,j} = C'_{i,j}},
\)
since
\(
	\max\bra{C_{i,j},C'_{i,j}} \le w_i w_j \lesssim s^2 \asymp p^{2(d-1)/k} \ll p.
\)

Now recall the partition $(P_3, ..., P_d)$ of $k$, and the definition $R_b = \abs{P_k}/k = (b-2)/\binom{d-1}2$.
Let $m$ be an integer (allowed to depend on other parameters) with $m \ll \mint{b} \abs{P_b} = k/\binomt{d-1}2 \asymp k/d^2$.
By exchangeability of the generators, for each $b \in \bra{3, ..., d}$ assume that the first $\tfrac45 \abs{P_b}$ entries $i$ of $P_b$ satisfy $\eta s \le w_i \le \eta^{-1} s$.
Our aim is to show that the mode of the vector $\bm C_m \cq \rbr{ C_{i,j} }_{i,j \in [m]}$, conditional on $W = w$, which we denote $\rho_m$, is bounded above by $s^{-c m^2}$, for some absolute constant $c > 0$.
We prove this in \cref{res-p1:cutoff:d:prod-bound:<:clm} below; for now, assume that claim.
Recall that $s \asymp p^{2(d-1)/k}$.

Recall that we showed in the $d = 3$ case that the mode for the single pair $(1,2)$ was upper bounded by $s^{-3/4}$; thus
\(
	- \log \rho_2
\gtrsim
	\log s.
\)
We are now analysing $(i,j)$ for all $i,j \in [m]$ simultaneously, which is $m^2$ pairs.
If the $C_{i,j}$ were independent over (not necessarily disjoint) pairs $(i,j)$, then the claim would follow immediately from the bound on the mode of $C_{1,2}$.
Unfortunately this is not the case. In some imprecise sense, we establish `approximate independence':
\(
	- \log \rho_m
\gtrsim
	m^2 \log s.
\)

Partition $\bra{1, ..., \tfrac45 \abs{P_b}}$ into
\(
	N
\cq
	\floor{ \tfrac45 \abs{P_b}/m }
\ge
	\frac34 \abs{P_b}/m
\)
sequential intervals of length $m$, say $I_{1,b}, ..., I_{N,b}$.
This allows us to decompose
\[
	\mce_b = \brb{ C_{i,j} = C'_{i,j} \: \forall \, i,j \in P_b }
\subseteq
	\cap_{\ell = 1}^N \brb{ C_{i,j} = C'_{i,j} \: \forall \, i,j \in I_{\ell,b} }.
\]
Moreover, the events in the intersection are independent.
We upper bound each using the mode:
\[
	q_b(t)
\le
	\rho_m^N
\le
	s^{- c m^2 N}
\le
	p^{- (d-1)k^{-1} \cdot c m^2 \cdot (3/4) \abs{P_b} / m }
\le
	p^{- d R_b \cdot \frac12 c m },
\label{eq-p1:cutoff:d:pb:penultimate}
\nt
\]
as $s \asymp p^{2(d-1)/k}$, so in particular $s \ge p^{(d-1)/k}$---recall that we had said that the exact value of $s$ would be unimportant.
We have $d R_b = (b-2) d/\binom{d-1}2 \ge 2(b-2)/d$, so this becomes
\[
	q_b(t)
\le
	p^{-c (b-2) m / d }.
\label{eq-p1:cutoff:d:qb:final}
\nt
\]
Set $m \cq \ceil{d/c}$. Then $c m / d \ge 1$.
Hence $q_b(t) \le 1/p^{b-2}$.

We still need $m \ll k/\binom{d-1}2 \asymp k/d^2$; since $m \asymp d$, this is equivalent to requiring $d^3 \ll k$.
Finally, to apply \cref{res-p1:cutoff:d:prod-bound:<:clm}, we need
\(
	\log m \le \tfrac43 \log \abs \uab / k.
\)
But $m \asymp d$ and $k \ll \log \abs \uab$, so this is implied by the condition $\log d \ll \tfrac32 \log \abs \uab / k$ from \cref{hyp-p1:cutoff}.

This establishes the desired bound, as it did for $\log \abs \uab \lesssim k \le (\log \abs \uab)^{1 + 2/(d-2)}$.
\end{Proof}

\begin{clm}
\label{res-p1:cutoff:d:prod-bound:<:clm}
	In the notation and under the assumptions of the above proof,
	there exists an absolute positive constant $c$ so that,
	assuming that $\log m \le \tfrac43 \log \abs \uab / k$ (so that $m \ll s$),
	we have
	\[
		\rho_m \le s^{-cm^2}.
	\]
\end{clm}

\begin{Proof}
First, note that
\(
	\log m \le \tfrac43 \log \abs \uab / k
\)
implies that
\(
	m \le \abs \uab^{(4/3)/k} \ll \abs \uab^{2/k} \asymp s.
\)

For this proof, we use the following notation:
	for $w \in \mcw$, write $\pr[w]{\cdot} \cq \pr{\, \cdot \mid W = W' = w}$ and $\ex[w]{\cdot}$ similarly;
	often we consider events that depend on $W$ but not on $W'$, in which case we ignore the conditioning on $W'$ (noting that $W$ and $W'$ are independent).
	Recall that we write $G_\ell \in [k]$ for the \emph{index} of the generators chosen in the $\ell$-th step.

We consider pairs $(i,j) \in [m]^2$ such that $i \ne j$, even when we do not say so explicitly.
Write $[m]^2_* \cq \bra{(i,j) \in [m]^2 \mid i \ne j}$.
Let $x_{i,j} \in \mbn_0$ for all $(i,j) \in [m]_*^2$; set $\bm x_m \cq (x_{i,j})_{(i,j) \in [m]^2_*}$.
Take an arbitrary ordering of all $K \cq m(m-1)$ pairs $(i,j) \in [m]^2_*$.
Write the $\kappa$-th term in this ordering as $\ell_\kappa$.
We are interested in controlling $\pr{ \bm C_m = \bm x_m \mid W = w }$; cf \cref{eq-p1:cutoff:3:<:mode-i}.
Then $\rho_m$ is the mode of this distribution, with the maximum taken over choices of $\bm x_m$.
We do this by sequentially estimating the conditional probabilities that $C_{i,j} = x_{i,j}$.
For each $\kappa \in [K]$, let $\chi_\kappa \cq \one{C_{\ell_\kappa} = x_{\ell_\kappa}}$, recalling that $(\ell_1, ..., \ell_K)$ is the chosen ordering of all pairs $(i,j) \in [m]^2_*$.
Then we want to bound
\(
	\ex[w]{ \chi_1 \cdots \chi_K }.
\)

To do this, we use the following general bound, which is an immediate consequence of the tower property for conditional expectation:
	for random variables $V_1, ..., V_K$, 
	we have
	\[
		\ex{ f_1(V_1, ..., V_K) \cdots f_K(V_1, ..., V_K) }
	&
	\le
		\MAX{v_1, ..., v_K}
		\ex{ f_1(V_1, v_2, ..., v_K) \cdots f_K(v_1, ..., v_{K-1}, V_K) }
	\\&\hspace{-3em}
	=
		\MAX{v_1, ..., v_K}
		\exB{ \prodt[K]{\kappa=1}
			\ex{ f_\kappa(v_1, ..., v_{\kappa-1}, V_\kappa, v_{\kappa+1}, ..., v_K \midb V_1, ..., V_{\kappa-1} }
		}.
	\]

The following argument is analogous to that given in the $d = 3$ case; see the end of the proof of \cref{res-p1:cutoff:3:prod-bound}.
Let $r \cq \sumt[m]{1} w_i$, and write our random word as $S = Z_{G_1} \cdots Z_{G_N}$; here $N = \sumt[k]1 w_i$ is the number of steps taken and $G_\ell$ is the generator index chosen in the $\ell$-th step.
Let $J_1 < \cdots < J_r$ be the (random) indices with $G_{J_\ell} \in [m]$.
Now define the vector $I \in \bra{1,...,m}^r$ by $I_\ell \cq G_{J_\ell}$.
Thus $I$ encodes the relative order between the different occurrences (with multiplicities) of the generators labelled by $[m]$ in the word $S$.
By typicality, $m \eta s \le r \le m \eta^{-1} s$.

Sequentially, and without replacements, for each pair $(i,j)$, or index $\kappa$, choose a uniformly random element of $\bra{\ell \mid G_\ell = j}$; call this $U_{i,j}$, or $U_\kappa$.
(This can be done since $\abs{\bra{\ell \mid G_\ell = j}} \gtrsim s \gg m$ by assumption.)
For each pair $(i,j)$, let $V_{i,j}$ be the location in $I$ of the random element $U_{i,j}$.
Now define the vector $I'$ from $I$ by omitting the $K = m(m-1)$ locations $\bra{ V_{i,j} }_{(i,j) \in [m]^2_*}$; so $I' \in \bra{1,...,m}^{r-K}$.
Importantly, we are omitting elements of the \emph{relative} order, not the \emph{absolute}~order.

By definition of $C_{i,j}$, given in \cref{eq-p1:cutoff:3:Cij-def}, given $W = w$, the value of $C_{i,j}$ is a function only of the relative locations $I$.
Hence, given $I'$ also, it is a function only of the location of the omitted $j$.
It is constant on the set of locations which give rise to the same relative order between choices of $i$ and $j$:
	two different placements of the omitted $j$ give rise to the same relative order between choices of $i$ and $j$ if and only if they both lie in the same (possibly empty) interval in which there are no $i$-s;
	this is analogous to the case $d = 3$.
(Recall that for the pair $(i,j)$ we omitted the location of a $j$.)

Recall that the random variables $\bra{ V_\kappa }$ are drawn uniformly at random without replacement from $[r]$; hence the distribution of $V_\kappa$ given $V_1, ..., V_{\kappa-1}$ is uniform on $[r] \setminus \bra{V_1, ..., V_{\kappa-1}}$.
Hence, writing $L_i \cq L_i(I')$ for the longest interval in $I'$ without an $i$ in it, we obtain
\[
	\ex[w]{ \chi_\kappa(v_1, ..., v_{\kappa-1}, V_\kappa, v_{\kappa+1}, ..., v_K) \midb V_1, ..., V_{\kappa-1}, \: I' }
\le
	\rbb{ L_\kappa + m(m-1) + 1 } / (r - \kappa+1).
\]

Hence, applying the above formula and noting that $m(m-1)+1 \le m^2$, we obtain

\[
	\ex[w]{ \chi_1 \cdots \chi_K \midb I' }
\le
	\frac
		{ \prodt[m]{1} \rbr{ L_i + m^2 }^{m-1} }
		{ r (r-1) \cdots (r - m(m-1) + 1) }.
\]
We now average over $I'$.
To bound this expectation, we use the generalisation of H\"older's inequality to the product of $m$ variables:
	for non-negative random variables $X_1, ..., X_m$,
	we have
	\[
		\ex{ X_1 \cdots X_m }
	\le
		\rbb{ \ex{X_1^m} \cdots \ex{X_m^m} }^{1/m}.
	\]
In our application of this, we take $X_i = \rbr{L_i + m^2}^{m-1}$.
This gives
\[
	\ex[w]{ \prodt[m]{1} \rbr{ L_i + m^2 }^{m-1} }
\le
	\maxt{i} \ex[w]{ \rbr{ L_i + m^2 }^{m(m-1)} }.
\]
Since $r \ge m \eta s$ and $s \gg m$, we have $m^2 \le \tfrac12 r$, and so the denominator is at least $2^{-m^2} r^{-m(m-1)}$; recall also that $r \ge m \eta s$.
In summary, we have proved that
\[
	\ex[w]{ \chi_1 \cdots \chi_K }
\le
	2^m
	(2 \eta^{-1} m^{-1})^{m(m-1)}
	s^{-m(m-1)}
	\maxt{i} \ex[w]{ \rbr{ L_i + m^2 }^{m^2} }.
\label{eq-p1:cutoff:d:chi-prod}
\nt
\]
Recall also that $s \asymp p^{2(d-1)/k}$.
It remains to bound this latter expectation.

To do bound this expectation, recall \cref{res-p1:cutoff:3:na}, which, for given $w$, states that
\[
	L_i/(\eta^{-2} m \log r) \preceq \Geom(\tfrac12).
\]
Using the the inequality $(a + b)^\ell \le 2^{\ell-1}(a^\ell + b^\ell)$, valid for $a,b \ge 0$ and $\ell \in \mbn$, we obtain
\[
	\rbb{ L_i + m^2 }^{m^2}
\le
	\rbb{ 2 \eta^{-2} m \log r }^{m^2} \rbb{ L_i/(\eta^{-2} m \log r }^{m^2}
+	\rbr{ 2 m^2 }^{m^2}.
\]
If $X \sim \Geom(\tfrac12)$, then one can show, for $\ell \ge 3$, that
\(
	\ex{ X^\ell } \le \ell^\ell.
\)
(This follows by comparison with the exponential-$(\log2)$ distribution.)
We apply this with $X \cq L_i/(\eta^{-2} m \log r)$ and $\ell \cq m^2$:
\[
	\ex[w]{ \rbr{ L_i + m^2 }^{m^2} }
\le
	\rbb{ 2 \eta^{-2} m \log r }^{m^2} \cdot (m^2)^{m^2}
+	(4 m^2)^{m^2}
\le
	2 \rbb{ 2 \eta^{-2} m^3 \log r }^{m^2}.
\]

Plugging this back into \cref{eq-p1:cutoff:d:chi-prod}, we obtain
\[
	- \log \ex[w]{ \chi_1 \cdots \chi_K }
\asymp
	m^2 \log s + m^2 \log m + m^2 \log r.
\]
Recalling that $r \asymp m s$ and
\(
	\log s
\asymp
	\log \abs \uab / k,
\)
we obtain
\[
	- \log \ex[w]{ \chi_1 \cdots \chi_K }
\asymp
	m^2 \rbb{ \log s + \log m }
\asymp
	m^2 \rbb{ \log \abs \uab / k + \log m}
\asymp
	m^2 \log \abs \uab / k,
\label{eq-p1:cutoff:d:chi-log}
\nt
\]
using the condition $\log m \le \tfrac43 \log \abs \uab / k$.
Recall that we desire
\[
	\ex[w]{ \chi_1 \cdots \chi_K }
\le
	s^{-c m^2},
\]
for some constant $c > 0$.
This holds since we established in \cref{eq-p1:cutoff:d:chi-log} that
\[
	- \log \ex{ \chi_1 \cdots \chi_K \midb W = w }
\asymp
	m^2 \log \abs \uab / k
\asymp
	m^2 \log s.
\qedhere
\]
\end{Proof}

\subsection{Undirected Cayley Graphs}
\label{sec-p1:cutoff:undir}

Thus far, we have always been assuming that the Cayley graph is directed.
Here we describe the required adaptations to consider \emph{undirected}, rather than \emph{directed}, graphs.
We still use an auxiliary process $W$ to generate the walk; it is now a SRW on $\mbz^k$, rather than a DRW on $\mbz_+^k$.

Recall that in the directed case the mixing time was the maximum of the time $t_0^+$ at which the entropy of $W$, ie a DRW on $\mbz_+^k$, reaches $\log \abs \uab$ and the diameter-based bound of $\log_k \abs \ugr$; see \cref{def-p1:intro:t*}.
The undirected case is completely analogous:
	the mixing time is the maximum of the time $t_0^-$ at which the entropy of $W$, now a SRW on $\mbz^k$, reaches $\log \abs \uab$ and the diameter-based lower bound of $\log_k \abs \ugr$.
(The directed graphs are $k$-regular, while the undirected graphs are $2k$-regular; but $\log(2k) \eqsim \log k$, so $\log_{2k} \abs \ugr \eqsim \log_k \abs \ugr$, when $k \gg 1$.)
We still use exactly the same outline, namely we use a modified $L_2$ calculation and `typicality'.

\medskip

The argument for the lower bound on the mixing time is completely unchanged, other than replacing the DRW with a SRW and the entropic time $t_0^+$ with the new entropic time $t_0^-$. See \S\ref{sec-p1:cutoff:lower}.

\renewcommand{\ugr}{\UU}
\renewcommand{\uab}{A}

We now focus on the upper bound.
Abbreviate $\ugr \cq \UU_{p,d}$ and $\uab \cq \UU_{p,d}^\ab$.
We first describe more general adaptations to the method.
The main difference is how we handle the entropy of $W$, which is now a SRW, and how this in turn affects the definition of typicality.
We then describe the adaptations to the regime $k \gtrsim \log \abs \uab$, which is the simpler case.
Finally we analyse $k \ll \log \abs \uab$.

Throughout the upper bound analysis, we assume that $\eps > 0$ is given and that
\[
	t = (1 + \eps) t_*
\Qwhere
	t_* = \max\brb{ t^+_0, \: \log \abs \ugr }.
\]

\subsubsection{General Adaptations}

\begin{Proof}[Expression for $S(t)$]
\qedtriangle
The primary difference comes from our expression for $S(t)$ given the generators chosen:
	for the DRW there were no inverses.
The expression for $S(t)$ is given by \cref{res-p1:cutoff:d:matrix-product}, which already allows for the application of inverses.
The only minor differences between the forms of $S(t)$ in the directed and undirected cases are the form of `remainder polynomial' $g$ and the expression for $C_{i,i}$. However, neither of these forms were needed for our proof.
The SRW generalises the DRW in some sense: if no inverses are applied, then the formulas for the SRW are exactly those for the DRW.
The entropic results of \S\ref{sec-p1:entropic} are analogous for both the DRW and SRW.

Let $W'$ be an independent copy of $W$, which is a SRW.
Previously $W$ was a DRW.
Nevertheless, $V \cq W - W'$ is still a rate-$2$ SRW on $\mbz^k$.
To evaluate $C_{i,j}$ in the undirected case, it is not enough to simply know the relative order of the choices of indices $(i,j)$: one needs to know which choice had which sign, ie know the relative order of the appearances of $(Z_i, Z_i^{-1}, Z_j, Z_j^{-1})$ in $S$.
(This corresponds to the `signs' $\sigma_\ell \in \bra{\pm1}$ in \cref{res-p1:cutoff:d:matrix-product}.)
We thus define
	$W_{+,i}(t)$ to be the number of times that $Z_i$ is applied
and
	$W_{-,i}(t)$ to be the number of times that $Z_i^{-1}$ is applied.
Clearly, $W_i(t) = W_{+,i}(t) - W_{-,i}(t)$.
Write $W_+(t) \cq (W_{+,i}(t))_{i \in [k]}$ and define $W_-(t)$, $W'_+(t)$ and~$W'_-(t)$~similarly.

We tend to suppress $t$ from the notation, eg writing $W_{-,i}$ for $W_{-,i}(t)$.
\end{Proof}

\begin{Proof}[Typicality Adaptations]
\qedtriangle
When considering $\pr{S = S' \mid W = W', \: \typ }$ previously,
we fixed some typical $w \in \mcw$ and conditioned on $W = w = W'$.
We now condition on the signed versions:
\[
	\rbr{ W_+, W_- } = (w_+, w_-)
\Quad{and}
	\rbr{ W'_+, W'_- } = (w'_+, w'_-)
\Quad{with}
	w_+ - w_- = w'_+ - w'_-.
\]
We need these to be `typical' in the appropriate sense.
We first define the entropic part:
\[
	\mcw
\cq
	\brb{ w \in \mbz^k \midb \mu_t(w) \le e^{-h}, \: \maxt{i} \abs{w_i} < \tfrac12 p };
\]
cf \cref{eq-p1:cutoff:prelim:typ:subset}.
Now $\mu_t$ is the law of SRW on $\mbz^k$ run for time $t = (1 + \eps) t_*$.
We recall that $h = h_0 + \omega$ where $h_0$ is (an approximation to) the entropy of $W(t_*)$ and $1 \ll \omega \ll \min\bra{k, \log \abs \uab}$;
see \cref{res-p1:cutoff:prelim:W=W',res-p1:cutoff:prelim:typ,def-p1:cutoff:prelim:h0}.
We require $w \cq w_+ - w_-$ and $w' \cq w'_+ - w'_-$ to be in $\mcw$.
This is the fundamental part of typicality.
It allowed us to control the term $\pr{ W = W' \mid \typ }$ in \cref{res-p1:cutoff:prelim:W=W'}:
\[
	\pr{ W = W' \mid \typ }
\le
	e^{-h} / \pr{\typ}.
\]

There are additional technical requirements, which differ between regimes.
For these we define a new set $\mcw'$ and require $(w_+, w_-), (w'_+, w'_-) \in \mcw'$---this is in contrast with $\mcw$ for which we only require $w,w' \in \mcw$.
\textit{Typicality} is
\(
	\typ
\cq
	\bra{W,W' \in \mcw} \cap \bra{(W_+, W_-), (W'_+, W'_-) \in \mcw'}.
\)
\end{Proof}

The analysis of $\pr{ S = S' \mid W \ne W', \: \typ }$ is exactly the same as in the directed case.
Where differences arise is in the analysis of $\pr{ S = S' \mid W = W', \: \typ }$.
This is largely due to the fact that knowing $W_i(t) = W_{+,i}(t) - W_{-,i}(t)$ does not determine $(W_{+,i}(t), W_{-,i}(t))$ for the SRW.

This completes the general adaptations.
We now move onto the specifics for the two regimes.

\subsubsection{Adaptations for $k \gtrsim \log \abs \uab$}

The adaptations in this regime are relatively simple.
In the directed case we already restricted to generators chosen once in $W$.
In the undirected case such a generator has $W_{+,i} + W_{-,i} = 1$.
The proof is thus almost analogous to the directed case.
We now proceed formally.

\medskip

We restrict to generators chosen at most once in $W$.
For $(w_+, w_-) \in \mbz_+^k \times \mbz_+^k$, write
\[
	\mcj(w_+, w_-)
\cq
	\brb{ i \in [k] \midb w_{+,i} + w_{-,i} = 1 }
\Qand
	J(w_+, w_-)
\cq
	\abs{ \mcj(w_+, w_-) }.
\]
Thus if $(W_+, W_-) = (w_+, w_-)$ then $\mcj(w_+, w_-)$ is the set of indices chosen precisely once.
Set
\[
	\mcw'
\cq
	\brb{
		(w_+, w_-) \in \mbz_+^k \times \mbz_+^k
	\midb
		\abs{ J(w_+, w_-) - s e^{-s} k } \le \tfrac12 \eps s e^{-s} k,
	\:	\maxt{i} \bra{w_{+,i}, w_{-,i}} < \tfrac12 p
	},
\]
where $s = t/k$ as always;
cf \cref{eq-p1:cutoff:d:typ:>=}.
This restriction means that
	if $\abs{W_i} = 1$ then $W_{+,i} = 1$ and $W_{-,i} = 0$ or vice versa
and
	if $\abs{W_i} = 0$ then $W_{+,i} = 0$ and $W_{-,i} = 0$.

When we proved \cref{res-p1:cutoff:d:prod-bound} for $k \gtrsim \log \abs \uab$ in the directed case, we also restricted to generators chosen at most once.
This is most clearly seen in the analysis for $d = 3$ (\S\ref{sec-p1:cutoff:3}) but is the case for general $d$ too.
Consider first $d = 3$.
There is at most one relative order of the appearances of
\(
	\rbr{ Z_i^{W_i} \mid i \in \mcj(W_+, W_-) }
\)
in $S$ such that
\(
	\rbr{ C_{i,j} \mid i,j \in \mcj(W_+, W_-), \: i \ne j }
\)
takes a particular value.
Exactly the same arguments as in the directed case then establish analogues of (\ref{eq-p1:cutoff:3:1/J!},~\ref{eq-p1:cutoff:3:q:=:ab}).

We partitioned the generators into blocks $P_3, ..., P_d$ to analyse general $d$.
We then upper bound
\[
	\brb{ C_{i,j} = C'_{i,j} \: \forall i,j \in P_b }
\subseteq
	\brb{ C_{i,j} = C'_{i,j} \: \forall i,j \in P_b \cap \mcj(W_+, W_-), \, i \ne j },
\]
for $b \in \bra{3, ..., d}$.
We choose the partition so that $P_b \cap \mcj(W_+, W_-)$ is typical for each $b$.
(We do exactly this in \S\ref{sec-p1:cutoff:d}.)
Again, there is at most one relative order of the appearances of
\(
	\rbr{ Z_i^{W_i} \mid w \in P_b \cap \mcj(W_+, W_-) }
\)
in $S$ such that
\(
	\rbr{ C_{i,j} \mid i,j \in P_b \cap \mcj(W_+, W_-), \: i \ne j }
\)
takes a particular value.
Exactly the same arguments as in the directed case then establish analogues of
(\ref{eq-p1:cutoff:d:1/t!},~\ref{eq-p1:cutoff:d:q:>=}).

The proof is concluded from these results analogously to the directed case.

\subsubsection{Adaptations for $k \ll \log \abs \uab$}

Now that we are in the regime $k \ll \log \abs \uab$, we can no longer assume that many generators are chosen at most once. Significantly more care is~needed.
The place where the proof needs to be adapted in this regime is at
	\cref{res-p1:cutoff:d:prod-bound:<:clm}.
This is more involved.
We needed to control a certain mode; see \cref{eq-p1:cutoff:d:pb:penultimate}, \cref{res-p1:cutoff:d:prod-bound:<:clm} and the surrounding material.
We give precise details below.

We first outline the approach, describing what we prove.
We then give a brief sketch of the ideas behind the proof.
The implicit details behind this sketch may not be completely clear initially.
Rather the sketch is intended as a `warm up' for the rigorous details, which come immediately after.


\begin{Proof}[Outline of Approach]
\qedtriangle
We suppress conditioning on `typicality' in this outline.

Initially we study the mode of the distribution of $C_{1,2}$, where the randomness is over the relative order in which the generators $(Z_1, Z_1^{-1}, Z_2, Z_2^{-1})$ are applied in $S$.
Denote this mode $\rho_2$.
Then
\[
	\pr{ C_{1,2} = C'_{1,2} \midb W = W' }
\le
	\rho_2.
\]
We show that
\(
	\rho_2 \lesssim s^{-1/3}.
\)
The method used is analogous to that for the corresponding mode~in~\S\ref{sec-p1:cutoff:3}.

The method which we use to estimate this `marginal' probability is amenable to considering the `joint' law of $C_{i,j}$ for all $i,j \in [m]$ simultaneously, for any fixed $m \ge 2$.
Denote by $\rho_m$ the mode of $( C_{i,j} )_{i,j \in [m]}$, given $(W_+, W_-)$.
Then
\[
	\pr{ C_{i,j} = C'_{i,j} \: \forall \, i,j \in [m] \midb W = W' }
\le
	\rho_m.
\]
Via an extension of the analysis for $\rho_2$,
we show that
\(
	- \log \rho_m
\asymp
	m^2 \log \rho_2
\gtrsim
	m^2 \log s.
\)
This first relation would be an equality if the events were independent.
The events are not independent, but the logarithm of the probability is only affected by a constant factor. In some imprecise sense, this is `approximate independence'.
The method of extension used is extremely similar to that for the corresponding extension in \S\ref{sec-p1:cutoff:d}.
Given this, the proof then follows exactly as in the directed case.

\smallskip

We now briefly recall how we proved \cref{res-p1:cutoff:d:prod-bound:<:clm} for the DRW.
There we considered the relative order of the choices of generators $Z_i$ and $Z_j$.
We redacted one $Z_j$ from this partial order and then considered the longest interval without $Z_i$.
The fact that the map $t \mapsto C_{i,j}(t)$ is increasing was crucial for this.
It meant that there was a unique interval in which we could place the redacted generator to obtain a specific value of $C_{i,j}$.
Now this is not the case because $C_{i,j}$ can decrease, when a generator is picked with an inverse.
We control this difference via a certain local time.
\end{Proof}

\medskip

\newcommand{\xx}{x}
\newcommand{\yy}{y}
\newcommand{\XX}{X}
\newcommand{\YY}{Y}

\renewcommand{\MM}{M}


We now briefly sketch the ideas behind these adaptations.

\begin{Proof}[Sketch for $k \ll \log \abs \uab$]
Here the technical requirements of typicality are that a large proportion of the generators are picked within a (large) constant factor of the expected number of times, say at least $\eta s$ and at most $\eta^{-1} s$ for some small constant $\eta > 0$, where $s = t/k$;
cf \cref{eq-p1:cutoff:d:pick-gen,eq-p1:cutoff:d:typ:<}.

Consider the relative order of two generators $Z_1^{\pm1}$ and $\pm Z_2^{\pm1}$.
Abbreviate $\pm \xx \cq Z_1^{\pm1}$ and $\pm \yy \cq Z_2^{\pm1}$.
We redact the signs of the $\pm \xx$ and location of one of the $\yy$-s from this order.
Suppose the redacted order has length $N$.
The redacted $\yy$ is then placed back in one of $N+1$ slots. Consider these slots sequentially.
As the location of the redacted $\yy$ moves from left-to-right, when it passes over a $\pm \xx$ the value of $C_{i,j}$ changes by $\pm 1$:
	if it passes over a $+ \xx$, then it increases;
	if it passes over a $- \xx$, then it decreases.
(In the directed case, we only had $+ \xx$.)
The value of $C_{i,j}$ is thus doing a type of SRW as $\yy$ moves along the available~slots.

Given an arbitrary target $\gamma$ for $C_{1,2}$, we wish to know how many of the $N+1$ slots for the redacted $\yy$ give rise to $C_{1,2} = \gamma$.
The value does not change unless a $\pm \xx$ is passed.
Thus the slots can be partitioned into sequential intervals according to the locations of the $\pm \xx$. All slots inside a given interval give rise to the same value.
Write $\MM$ for the maximal length of an interval.

In the directed case, the value was a DRW, increasing each time a $+ \xx$ is passed, so there was at most one interval which would give rise to $\gamma$. We thus bounded the probability by $\MM/N$.
We had $N \asymp s \cq t/k$ and $\MM \lesssim \log s$ with sufficiently high probability.
We thus upper bounded the probability by $M/N \lesssim (\log s) / s \lesssim s^{-1/3}$.
This bound was sufficient in the proof of \cref{res-p1:cutoff:d:prod-bound:<:clm}.

Now, however, the value is a SRW so there may be many such intervals.
Let $R = (R_\ell)_{\ell=0}^L$ be the corresponding SRW: $R_\ell$ is the value of the sum in interval $\ell$; set $R_0 \cq 0$.
We now have to multiply the value $M/N$ by the number of intervals giving rise to the value $\gamma$. But this factor is at most the maximal local time of the SRW $R$.
It is standard that the maximal local time of a SRW run for time order $s$ is at most $\sqrt s \log s$ with very high probability.
(It is typically order $\sqrt s$.)
We thus end up with a final probability order at most
\(
	(\log s / s) \cdot (\sqrt s \log s)
=
	(\log s)^2 / \sqrt s
\lesssim
	s^{-1/3}.
\)
	%
\end{Proof}


\medskip

We now explain the required adaptations more carefully.
First we set up typicality.

\begin{Proof}[Typicality]
\qedtriangle
	%
We use the general adaptations described above.
We define the technical requirements of typicality, $\mcw'$, and describe how to use this to establish the cutoff result.
We restrict to generators chosen the correct order of times:
	for $(w_+, w_-) \in \mbz_+^k \times \mbz_+^k$
	write
	\[
		\mcc(w_+, w_-)
	\cq
		\brb{ i \in [k] \midb \eta s \le \min\bra{w_{+,i}, w_{-,i}} \le \max\bra{w_{+,i}, w_{-,i}} \le \eta^{-1} s },
	\]
	where $\eta > 0$ is a (small) constant.
We set
\[
	\mcw'
\cq
	\brb{
		(w_+, w_-) \in \mbz_+^k \times \mbz_+^k
	\midb
		\abs{\mcc(w_+, w_-)} \ge \tfrac45 k,
	\:	\maxt{i} \bra{w_{+,i}, w_{-,i}} < \tfrac12 p
	};
\]
cf \cref{eq-p1:cutoff:d:pick-gen,eq-p1:cutoff:d:typ:<}.
We can choose $\eta > 0$ small enough so that typicality holds whp.
We restrict attention to indices $i$ with
\(
	\eta s
\le
	\min\bra{W_{+,i},W_{-,i}}
\le
	\max\bra{W_{+,i},W_{-,i}}
\le
	\eta^{-1} s.
\)
\end{Proof}


\medskip

Next we analyse the `marginal' law $\bra{ C_{i,j} = C'_{i,j} }$ for a single pair $(i,j)$; without loss of generality we take $(i,j) = (1,2)$.
Afterwards we analyse the `joint' law $\bra{ C_{i,j} = C'_{i,j} \: \forall \, i,j \in [m] }$, for some $m$.
We do this by showing a form of `approximate independence' for the events $\bra{ C_{i,j} = C'_{i,j} }$.

To analyse the marginal law, we bound the mode of $C_{1,2}$, which we denote $\rho_2$, and use the trivial inequality
\(
	\pr{ C_{1,2} = C'_{1,2} }
\le
	\rho_2.
\)
Here we are suppressing the conditioning from the notation.

We first set-up some notation and a partitioning of the relative order which will be key.
We subsequently use this partitioning and a comparison with local times of SRW to estimate the mode.

\begin{Proof}[Partitioning]
\qedtriangle
	%
Write $\pm \xx \cq Z_1^{\pm1}$ and $\pm \yy \cq Z_2^{\pm1}$.
We calculate $C_{1,2}$ via~the~following~sum:
\begin{itemize}[noitemsep, topsep = \smallskipamount, label = \bcdot]
	\item 
	for each $+\yy$ that comes after a $+\xx$ (not necessarily directly), add $1$;
	
	\item 
	for each $-\yy$ that comes after a $-\xx$ (not necessarily directly), add $1$;
	
	\item 
	for each $-\yy$ that comes after a $+\xx$ (not necessarily directly), subtract $1$;
	
	\item 
	for each $+\yy$ that comes after a $-\xx$ (not necessarily directly), subtract $1$.
\end{itemize}
An inspection of \cref{res-p1:cutoff:d:matrix-product} shows that this does indeed calculate $C_{1,2}$ correctly.
We reveal the relative order of all generators, except for one appearances of $+\yy$; we term this the \textit{redacted $\yy$}.
The relative order is uniformly random and thus the location of the redacted $\yy$ is uniform.

We try to control the mode of $C_{1,2}$, as mentioned above.
Fix some arbitrary target $\gamma$ for $C_{1,2}$.
\begin{quote}
	Given what we see,
	namely the order without the redacted $\yy$,
	what is the probability that the redacted $\yy$ is at a location that makes the above sum equal to $\gamma$?
\end{quote}
Specifically, in the proof of \cref{res-p1:cutoff:d:prod-bound:<:clm} it is enough to show that the probability in the above question decays like an inverse power of $s$.
We show here a decay of at least $s^{-1/3}$.

Assume that we did not reveal the \emph{sign} of the $\pm \xx$-generators, only whether they are either $\pm \xx$, not $\pm \yy$.
Indicate this by writing $\XX$ for a generator which is revealed to be one of $\pm \xx$ and $\YY$ for $\pm \yy$ similarly.
We have redacted a single $+ \yy$.
We wish to control the sum upon replacement of this redacted $\yy$.
If there are $N$ letters in the redacted word, then there are $N+1$ slots for the~redacted~$\yy$.

We partition the redacted order into consecutive intervals in $\bra{0, 1, ..., N}$ according to the locations of the $\XX$:
	the first interval is $\bra{0}$ union the locations of any $\YY$-s before the first $\XX$;
	subsequent intervals start at the location of an $\XX$ and end one before the location of the following $\XX$.
For example,
	the word $\XX \YY \YY \XX \YY$
gives rise to
	the intervals $(\bra{0}, \bra{1,2,3}, \bra{4,5})$
and
	the word $\YY \XX \XX \YY \YY$
gives rise to
	the intervals $(\bra{0,1}, \bra{2}, \bra{3,4,5})$.
Placing the redacted $\yy$ in position $r \in \bra{1, ..., N-1}$ corresponds to placing it between the $r$-th and $(r+1)$-th letters in the (non-redacted) word; $r = 0$ corresponds to the left-most slot and $r = N$ to the right-most slot.

If there are $L$ appearances of $\XX$, then there are $L+1$ intervals.
Write $I = (I_\ell)_{\ell \in [L+1]}$ for the sequence of intervals---here $L$ is non-random after conditioning on the number of $\XX$-s and, by typicality, satisfies $L \in [2 \eta s, 2 \eta^{-1} s]$.
Reordering the $\YY$-s inside an interval has no effect on the value of the sum.
Thus one does not need to know the exactly location of the redacted $\yy$, but rather just the interval.
We control the number of intervals which give rise to a given value of the sum by comparison with the local times of a \textit{random bridge} on $\mbz$, which we define now.

Let $R = (R_\ell)_{\ell=0}^L$ be a discrete-time process on $\mbz$ with $R_0 = 0$.
We couple the process $R$ with the redacted order in the following way.
	Let $\Sigma$ be the vector of $\pm1$-s denoting the signs of the $\XX$-s at the start of the intervals; then $\Sigma \in \bra{\pm1}^L$.
	Let the $\ell$-th step of $R$ be given by $\Sigma_\ell \in \bra{\pm1}$.
Importantly, we are conditioning on how many $+$-steps and how many $-$-steps this takes, but not the order in which they are taken; this is equivalent to conditioning on the final value $R_L$.
Except for this conditioning, $R$ is the normal SRW in $\mbz$; hence it is a random bridge.
\end{Proof}

We now describe how to couple the redacted order with the random bridge $R$.
We then upper bound the mode $\rho_2$ in terms of the local time of $R$ and the maximal interval length.

\begin{Proof}[Comparison with Local Times]
\qedtriangle
We establish the following relation, with terms explained below:
\[
	\rho_2
\lesssim
	s^{-1} \sqrt{ \ex{\mcl(\Sigma,I)^2} \ex{\MM(I)^2} }.
\]
Here $\mcl$ is the maximal local time of $R$ defined below and $\MM(I)$ is the maximal interval length in~$I$.

\smallskip

Recall how we calculated the sum $C_{1,2}$ above.
Denote by $S_0$ the value of the sum $C_{1,2}$ without the redacted $\yy$.
The coupling between $R$ and the redacted order is key:
	if we place the redacted $\yy$ in the $\ell$-th interval,
	then $C_{1,2} = S_0 + R_\ell$.
Thus the number of intervals which give rise to a specific value of $C_{1,2}$ is at most the \textit{maximum local time} of $R$, defined precisely now:
	define $\mcl_r$ to be the \textit{local time} of $R$ at $r \in \mbz$,
	ie
	\(
		\mcl_r
	\cq
		\sumt[L]{\ell=0}
		\one{ R_\ell = r };
	\)
	define $\mcl \cq \maxt{r \in \mbz} \mcl_r$ to be the \textit{maximum local time}.
We make explicit the dependence on $\Sigma$ and $I$ by writing $\mcl(\Sigma, I)$.
The only dependence on $I$ is through $L$, the number of partitions.
(The value of $L$ is needed to define $\Sigma$.)

The location of the redacted $\yy$ is uniform.
The intervals are not necessarily of the same length and thus into which \emph{interval} the redacted $\yy$ falls is not uniform: it is proportional to the length of the interval.
Denote by $\abs{I_\ell}$ the length of interval $I_\ell$, for $\ell \in [L]$, and write
\(
	\MM
\cq
	\maxt{\ell \in [L]}
	\abs{I_\ell}
\)
for the maximum interval length.
We make explicit the dependence on $I$ by writing $\MM(I)$.

Combining the local time and maximal interval length,
we get
\(
	\rho_2
\le
	\ex{ \mcl(\Sigma, I) \MM(I) / L }.
\)
Again, we are hiding implicit typicality conditioning from the notation.
By typicality, $L \ge 2 \eta s$.
Hence
\[
	\rho_2
\le
	(2 \eta s)^{-1} \ex{ \mcl(\Sigma, I) \MM(I) }
\asymp
	s^{-1} \ex{ \mcl(\Sigma, I) \MM(I) }
\le
	s^{-1} \sqrt{ \ex{\mcl(\Sigma,I)^2} \ex{\MM(I)^2} }.
\qedhere
\]
\end{Proof}

Again, we implicitly suppressed conditioning for notational clarity.
An enthusiastic reader may add the appropriate conditioning to all the probabilities and expectations if she desires.
This does nothing but obfuscate the ideas behind the method, in our opinion, hence its omission below.

It remains to estimate these two expectations.
We consider the maximal interval length first.

\begin{Proof}[Expected Maximal Interval Length]
\qedtriangle
This interval length has the same distribution as in the directed case.
Recall that $L \le 2 \eta^{-1} s$, by typicality.
Thus applying \cref{res-p1:cutoff:3:na} gives
\[
	\ex{ \MM(I)^2 }
\lesssim
	\log(2 \eta^{-1} s)^2
\eqsim
	(\log s)^2.
\qedhere
\]
\end{Proof}

\begin{Proof}[Maximal Local Time]
\qedtriangle
It is an immediate consequence of
\cref{res-p1:cutoff:undir:local-time:bridge}
that
\[
	\maxt{r}
	\ex{ \mcl(\Sigma,I)^2 \mid R_L = r }
\lesssim
	(s^{3/4})^2.
\qedhere
\]
\end{Proof}

\begin{Proof}[Conclusion for Marginal Analysis]
We have shown that
\[
	\rho_2
\lesssim
	s^{-1} \ex{ \mcl(\Sigma, I) \MM(I) }
\le
	s^{-1} \sqrt{ \ex{\mcl(\Sigma,I)^2} \ex{\MM(I)^2} }
\lesssim
	(\log s)^2 / \sqrt s
\lesssim
	s^{-1/3}.
\qedhere
\]
\end{Proof}


\medskip

We have established the desired `marginal' estimate:
\[
	- \log \rho_2
\gtrsim
	\log s.
\]
We now extend this to the `joint' analysis.
Recall that we studied a single pair in \S\ref{sec-p1:cutoff:3} ($d = 3$) and extended it to a collection of pairs in \S\ref{sec-p1:cutoff:d} (general $d$), showing `approximate independence'.
We use exactly the same methodology here.
The only difference is that we have to include the local time in the analysis---other than that, the following extension is exactly the same as the extension in \S\ref{sec-p1:cutoff:d}.
The moments of this local time are sufficiently well behaved that they do not cause any major concerns.
Denote by $\rho_m$ the mode of $(C_{i,j})_{i,j \in [m]}$.
We establish the desired `joint' estimate:
\[
	- \log \rho_m
\gtrsim
	- m^2 \log \rho_2
\gtrsim
	m^2 \log s.
\]

\medskip

We consider pairs $(i,j) \in [m]^2$ such that $i \ne j$, even when we do not say so explicitly.
Write $[m]^2_* \cq \bra{(i,j) \in [m]^2 \mid i \ne j}$.
Let $x_{i,j} \in \mbn_0$ for all $(i,j) \in [m]_*^2$; set $\bm x_m \cq (x_{i,j})_{(i,j) \in [m]^2_*}$.
Take an arbitrary ordering of all $K \cq m(m-1)$ pairs $(i,j) \in [m]^2_*$.
Write the $\kappa$-th term in this ordering as $\ell_\kappa$.
We are interested in controlling $\pr{ \bm C_m = \bm x_m \mid W = w }$; cf \cref{eq-p1:cutoff:3:<:mode-i}.
Then $\rho_m$ is the mode of this distribution, with the maximum taken over choices of $\bm x_m$.
We do this by sequentially estimating the conditional probabilities that $C_{i,j} = x_{i,j}$.
For each $\kappa \in [K]$, let $\chi_\kappa \cq \one{C_{\ell_\kappa} = x_{\ell_\kappa}}$, recalling that $(\ell_1, ..., \ell_K)$ is the chosen ordering of all pairs $(i,j) \in [m]^2_*$.
Then we want to bound
\(
	\ex{ \chi_1 \cdots \chi_K }.
\)

\begin{Proof}[Joint Analysis]
We use the following general bound, which we already used for directed graphs:
	\[
		\ex{ f_1(V_1, ..., V_K) \cdots f_K(V_1, ..., V_K) }
	&
	\le
		\MAX{v_1, ..., v_K}
		\ex{ f_1(V_1, v_2, ..., v_K) \cdots f_K(v_1, ..., v_{K-1}, V_K) }
	\\&\hspace{-3em}
	=
		\MAX{v_1, ..., v_K}
		\exB{ \prodt[K]{\kappa=1}
			\ex{ f_\kappa(v_1, ..., v_{\kappa-1}, V_\kappa, v_{\kappa+1}, ..., v_K \midb V_1, ..., V_{\kappa-1} }
		}.
	\]
We are looking at $K = m(m-1)$ pairs simultaneously.
We did exactly this type of argument for directed graphs.
The only difference is that now instead of obtaining a product over longest intervals we also have the expected maximum local time.
We obtain an analogue of \cref{eq-p1:cutoff:d:chi-prod}:
\[
	\ex{ \chi_1 \cdots \chi_K }
&
\lesssim
	s^{-m(m-1)} \ex{ (\mcl(\Sigma, I) \MM(I))^{m(m-1)} }
\\&
\le
	s^{-m(m-1)} \sqrt{ \ex{ \mcl(\Sigma, I)^{2m(m-1)} } \ex{ \MM(I)^{2m(m-1)} } }.
\]
We control the length of the longest interval via \cref{res-p1:cutoff:3:na}.
\cref{res-p1:cutoff:undir:local-time:bridge} states that
\[
	\maxt{r}
	\ex{ \mcl(\Sigma,I)^C \mid R_L = r }
\lesssim
	C^C L^{3C/5}
\Quad{whenever}
	\text{$C \ge 1$},
\]
where the implicit constant is absolute.
Recall that $L \asymp s$.
We desire the relation
\[
	\maxt{r}
	\ex{ \mcl(\Sigma,I)^C \mid R_L = r }
\lesssim
	s^{3C/4}.
\]
Thus $C \lesssim s^{3/20}$ suffices.
Assume this for now.
We then deduce, as in the directed case,~that
\[
	- \log \ex{ \chi_1 \cdots \chi_k }
\asymp
	m^2 \log s^{3/5}
\asymp
	m^2 \log s,
\]
We have thus upper bounded the mode by $s^{-c m^2}$, for some constant $c > 0$.
The proof is completed exactly as before, with $m \cq \ceil{d/c} \asymp d$.

It remains to show that we need only consider $C \lesssim s^{3/20}$.
Above we are considering the $2m(m-1)$-th moment. Thus we take $C \cq 2m(m-1) \asymp m^2 \asymp d^2$.
Recall that $k \ll \log \uab / \log d$ by \cref{hyp-p1:cutoff} and $s \asymp \abs \uab^{2/k}$.
Thus
\(
	k \ll k \log s / \log d,
\)
ie $d = s^{\oh1}$.
Thus $C \asymp d^2 = s^{\oh1} \le s^{3/20}$.
\end{Proof}


\medskip

It remains to control the maximal local time.
The following lemma is deferred to \cite[\cref{res-p0:cutoff:undir:local-time:bridge}]{HOt:rcg:supp}.

\begin{lem}
\label{res-p1:cutoff:undir:local-time:bridge}
	Let $r \in \mbz$ and let $R = (R_\ell)_{\ell=0}^L$ be a random bridge from $0$ to $r$ of length $L$.
	Denote by $\mcl(r,L)$ the maximal local time of $R$.
	Let $C \ge 1$.
	Then
	\[
		\maxt{r}
		\ex{ \mcl(r,L)^C }
	\lesssim
		C^C L^{3C/5},
	\]
	where the implicit constant is absolute.
\end{lem}

We prove this via a reduction from a random bridge to a random walk.
We then use the following lemma on the tails of the local time at $0$ for SRW on $\mbz$, proved in \cite[\cref{res-p0:cutoff:undir:local-time:walk}]{HOt:rcg:supp}.

\begin{lem}
\label{res-p1:cutoff:undir:local-time:walk}
	Let $X = (X_\ell)_{\ell=0}^\infty$ be a SRW on $\mbz$ run for $L$ steps.
	Denote by $\mcl'(0,L)$ the local time of $X$ at $0$ in the first $L$ steps.
	Then, for all $k \ge 0$, we have
	\[
		\pr{ \mcl'(0,L) > k }
	\le
		\expb{ - \tfrac12 k^2 / L }.
	\]
\end{lem}

\section{Extensions to Cutoff Theorem}
\label{sec-p1:ext}

\renewcommand{\ugr}{\UU_{p,d}}
\renewcommand{\uab}{\UU_{p,d}^\ab}

In this section, we describe some extensions to the cutoff arguments.

\begin{itemize}[itemsep = 0pt, topsep = \smallskipamount, label = \bcdot]
	
	\item [\S\ref{sec-p1:ext:window}]
	We describe the limit profile in the regime $1 \ll k \ll \log \abs \uab$.
	
	\item [\S\ref{sec-p1:ext:comp}]
	We relax that condition that $p$ (in $\ugr$) is prime when $k \gtrsim \log \abs \ugr$.
	
	\item [\S\ref{sec-p1:ext:heis}]
	We adapt the argument to consider the $d$-dimensional Heisenberg group.
\end{itemize}

\subsection{Lifting the Primality Condition for $p$}
\label{sec-p1:ext:comp}

Thus far, we have always been assuming that $p$ is prime.
Here we describe how to remove this assumption under suitable conditions.
To emphasise the lack of primality, we consider $\UU_{m,d}$ with $m,d \in \mbn$---ie replace the letter $p$ by the letter $m$.

\renewcommand{\ugr}{\UU}
\renewcommand{\uab}{\UU^\ab}
\renewcommand{\uk}{\UU_k}


\begin{customthm}{\ref{res-p1:intro:ext:comp}}
	Let $\ugr \cq \UU_{m,d}$.
	Suppose that
		$1 \ll \log k \ll \log \abs \ugr$,
		$k \gg \max\bra{ \sqrt{\log m}, \: \log \div m }$
	and
		$d \ge 3$ is fixed or diverges sufficiently slowly.
	Then the RW on $\uk^\pm$ exhibits cutoff whp~at~$t_*^\pm(k, m, d)$.
\end{customthm}

\renewcommand{\ugr}{\UU_{m,d}}
\renewcommand{\uab}{\UU_{m,d}^\ab}
\renewcommand{\uk}{(\UU_{m,d})_k}

\begin{Proof}[Proof of Lower Bound]
First note that the lower bound holds easily: all it required was that the walk projected to the Abelianisation was not mixed. The lower bound for mixing on an Abelian group is valid for any Abelian group and any choice of generators.
\end{Proof}

We now turn to the upper bound.
Fix some $\eps > 0$ and let $t = t_* (1 + \eps)$;
herein we suppress $t$ (and $\eps$) from the notation.
Let $W'$ be an independent copy of $W$ and define $S'$ via $W'$; write $V \cq W - W'$.
We are interested in the probability that $S = S'$ when $W$ and $W'$ are~typical.

We use \emph{almost} exactly the same definition of typicality as for $m = p$ prime.
We perform some gcd-analysis in the proof of \cref{res-p1:ext:comp:S=S'|W=/=W'}. This requires strengthened bounds on $\max_i \abs{W_i}$.
The precise details play no role in the analysis below, other than being required in the references invoked in the proof of \cref{res-p1:ext:comp:S=S'|W=/=W'} to control a certain gcd.
Precise references for the adjustment to typicality are deferred to that~proof.
The strengthened requirements still hold whp.

\renewcommand{\ugr}{\UU}
\renewcommand{\uab}{A}

We now state a series of lemmas.
The proofs of the lemmas are deferred until after we explain how to combine the lemmas to deduce the theorem.
We abbreviate $\ugr \cq \UU_{m,d}$ and $\uab \cq \UU_{m,d}^\ab$.

\begin{lem}
\label{res-p1:ext:comp:S=S'|W=/=W'}
	We have
	\[
		\abs \ugr \, \pr{ S = S' \mid W \ne W', \: \typ }
	\le
		1 + \oh1
	\]
\end{lem}

The entropic analysis does not require $m = p$ to be prime.
Thus
\cref{res-p1:cutoff:prelim:W=W',res-p1:cutoff:prelim:typ,res-p1:cutoff:prelim:omega-choice,res-p1:cutoff:prelim:h0} still hold with $h_0$ and $\omega$ given by \cref{def-p1:cutoff:prelim:h0}.
In particular, $1 \ll \omega \ll \min\bra{k, \: \log \abs \uab}$ and
\begin{empheq}
[ left = {h_0 = \empheqlbrace} ]
{alignat = 3}
	&\log \abs \uab
		&\Qwhen&
	k \le (\log \abs \uab)^{1 + 2/(d-2)},
\nonumber
\\
	&(1 - 1/\rho) \log \abs \ugr
		&\Qwhen&
	k \ge (\log \abs \uab)^{1 + 2/(d-2)};
\nonumber
\end{empheq}
set $h \cq h_0 + \omega$.
We restate \cref{res-p1:cutoff:prelim:W=W'} here for convenience---we label it \cref{res-p1:ext:comp:W=W'} here.

\begin{lem}
\label{res-p1:ext:comp:W=W'}
	We have
	\[
		\pr{ W = W' \mid \typ }
	\le
		e^{-\omega} e^{-h_0} / \prt{\typ}.
	\]
\end{lem}

Recall that $\div m = \sumt{r \in [m]} \one{ r \wr m }$ is the number of divisors of $m$.

\begin{lem}
\label{res-p1:ext:comp:k-large:S=S'|W=W'}
	Suppose that $k \gtrsim \log \abs \uab$.
	We have
	\[
		\abs \ugr \, \pr{ S = S' \mid W = W', \: \typ }
	\lesssim
		e^{h_0} \cdot (\div m)^{d^2}.
	\]
\end{lem}

When considering all $m$,
we consider the regime $k \gtrsim \log \abs \uab$.
(When we restrict $m$ to a density-$1$ subset of $\mbn$, we require only $k \gg \sqrt{d^3 \log m}$.)
We can thus conclude in the case of all $m$.

\begin{Proof}[Proof of \cref{res-p1:intro:ext:comp}: $k \gtrsim \log \abs \uab$]
Combining \cref{res-p1:ext:comp:S=S'|W=/=W',res-p1:ext:comp:W=W',res-p1:ext:comp:k-large:S=S'|W=W'},
we deduce that
\[
	\abs \ugr \, \pr{ S = S' \mid \typ } - 1
=
	\oh1 + \Ohb{ e^{h_0} (\div m)^{d^2} \cdot e^{-\omega} e^{-h_0} }
=
	\oh1 + \Ohb{ e^{-\omega} (\div m)^{d^2} }.
\]
It remains to show that this final $\Oh{\cdot}$ term is $\oh1$ and apply the modified $L_2$ calculation.

We have $k \gtrsim \log \abs \uab$ and may take any $\omega$ with $1 \ll \omega \ll \min\bra{k, \: \log \abs \uab} \asymp \log \abs \uab \asymp d \log m$.
It is well-known that $\div r \le r^{\Oh{1/\log\log r}}$ for all $r \in \mbn$; see \cite[\S 18.1]{HW:number-theory}.
Thus
\(
	\div m = m^{\oh1}.
\)
Thus if $d$ does not diverge too quickly then we can choose $\omega \ll d \log m$ such that
\(
	e^{-\omega} (\div m)^{d^2}
=
	\oh1.
\)
\end{Proof}

For the regime $k \ll \log \abs \uab$, we need to restrict $m$ to a density-$1$ subset of $\mbn$.

\begin{lem}
\label{res-p1:ext:comp:k-small:S=S'|W=W'}
	Suppose that $k \ll \log \abs \uab$.
	We have
	\[
		\abs \ugr \, \pr{ S = S' \mid W = W', \: \typ }
	\lesssim
		e^{h_0} m^{5d^3/k} (\div m)^{d^2}.
	\]
\end{lem}


We now find a suitable density-$1$ subset of $\mbn$ to relax the conditions on $k$.

\begin{Proof}[Proof of \cref{res-p1:intro:ext:comp}: $k \ll \log \abs \uab$]
Combining \cref{res-p1:ext:comp:S=S'|W=/=W',res-p1:ext:comp:W=W',res-p1:ext:comp:k-small:S=S'|W=W'},
we deduce that
\[
	\abs \ugr \, \pr{ S = S' \mid \typ } - 1
=
	\oh1 + \Ohb{ e^{-\omega} m^{5d^3/k} (\div m)^{d^2} }.
\]
It remains to show that this final $\Oh{\cdot}$ term is $\oh1$ and apply the modified $L_2$ calculation.
We have
\[
	e^{-\omega} m^{5d^3/k} (\div m)^{d^2}
=
	\expb{ - \omega + 5 d^3 \log m / k + d^2 \log \div m }.
\]
We have $k \ll \log \abs \uab \asymp d \log m$ and may take any $\omega$ with $1 \ll \omega \ll \min\bra{k, \: \log \abs \uab} \asymp k$.
Since we require $d \asymp 1$ or $d \gg 1$ sufficiently slowly,
we thus require
\[
	\omega
\gg
	\max\brb{ \log m / k, \: \log \div m },
\Quad{or equivalently}
	k
\gg
	\max\brb{ \sqrt{\log m}, \: \log \div m }.
\qedhere
\]
\end{Proof}

It remains to give the deferred proofs of the above lemmas.

\begin{Proof}[Proof of \cref{res-p1:ext:comp:S=S'|W=/=W'}]
When $m = p$ is prime, for any $v \in \mbz_p^k \setminus \bra{0}$, each of $\sumt{i} A_i v_i$, $\sumt{i} B_i v_i$ and $\sumt{i} C_i v_i$ was an independent $\Unif(\mbz_p)$ random variable, since $\bra{A_i, B_i, C_i}_1^k$ is an independent collection of $\Unif(\mbz_p)$-s.
Hence
in the $(a,b)$-th coordinate there is a sum $K_{a,b} \cq \sumt{i} Z_i(a,b) V_i(t) \sim \Unif(\mbz_p)$, independent of the other terms in the coordinate.
We deduce that
\[
	\pr{ S = S' \mid W \ne W', \: \typ }
=
	1 / \abs{ \UU_{p,d} }
=
	1 / p^{d(d-1)/2}.
\]

This is not so when $m$ is composite
Instead, $K_{a,b} \sim \mfgcd_v \mbz_m$ where $\mfgcd_v \cq \gcd(v_1, ..., v_k, m)$ given $(W,W') = (w,w')$ with $v \cq w - w'$.
This is proved in \cite[\cref{res-p0:deferred:unif-gcd}]{HOt:rcg:supp}.
Exactly the same arguments used above,
noting that there are $\tfrac12 d(d-1)$ non-trivial coordinates in the matrix,
give
\[
	\pr{ S = S' \mid W = w, \: W' = w' }
\le
	\rbr{ \mfgcd_{w-w'}/m }^{d(d-1)/2}
\Quad{for any}
	w,w' \in \mbz^k.
\label{eq-p1:ext:comp:S=S'|WW'}
\nt
\]
We then need to take expectation of $(W, W')$ with $W \ne W'$.
We argue that
\[
	\ex{ \mfgcd_{W-W'}^{d(d-1)/2} \mid W \ne W', \: \typ }
=
	1 + \oh1.
\label{eq-p1:ext:comp:gcd=1+o1}
\nt
\]
Given this, the lemma follows immediately from the combination of
\cref{eq-p1:ext:comp:S=S'|WW',eq-p1:ext:comp:gcd=1+o1}.

We carry out exactly this type of gcd-expectation when analysing Abelian groups in companion papers \cite{HOt:rcg:abe:cutoff,HOt:rcg:abe:extra}.
The argument can be sketched as follows.
	The probability that the gcd of the non-zero coordinates of SRW on $\mbz^k$ is $\gamma$ decays like at least $\gamma^{-\ell}$ for $\gamma \ge 2$ given that there are $\ell$ non-zero coordinates.
	The contribution to the gcd from values larger than $1$ is then $\oh1$ when $\ell - d \gg 1$.
	Smaller $\ell$ is atypical and can be controlled.
The precise references are as follows:\par\smallskip
{\centering
	see \cite[\cref{res-p2:cutoff1:gcd-ex}]{HOt:rcg:abe:cutoff} for $k \ll \log \abs \uab$;\quad%
	see \cite[\cref{res-p5:profile:gcd-ex}]{HOt:rcg:abe:extra}  for $k \gtrsim \log \abs \uab$.%
\par\smallskip}\noindent%
We briefly describe the notation there:
	$G$ is an Abelian group and $d(G)$ is the minimal size of a generating subset of $G$;
	$\mfgcd = \gcd\rbr{ V_1(t), ..., V_k(t), \abs G }$ where $V \cq W - W'$ (which have the same definition as here).
There are hypotheses on $(k, G)$ in each reference; see \cite[\cref{hyp-p2:cutoff1}]{HOt:rcg:abe:cutoff} and \cite[\cref{hyp-p5:profile}]{HOt:rcg:abe:extra}.
One should apply these results \emph{to the Abelianisation}, ie set $G \cq \uab \cong \mbz_m^{d-1}$.
Indeed, the time considered there is the entropic time for $G$; here it is for $\uab$.
If $d$ grows sufficiently slowly (or is fixed), then these hypotheses are easily seen to be satisfied by $(k, \uab)$.

We note a couple of minor differences between the set-up there and the set-up here.
The exponent of $\mfgcd$ there is $d(G) = d(\uab) = d-1$; here it is $\tfrac12 d(d-1)$.
A simple inspection of the proof shows that the this makes no difference if $d$ grows sufficiently slowly (or is fixed).
Also, the references require a slightly strengthened definition of typicality, as noted above.
The only change is in imposing a stronger restriction on $\max_i \abs{W_i}$:
see
	\cite[\cref{def-p2:cutoff1:typ}]{HOt:rcg:abe:cutoff}
and
	\cite[\cref{def-p5:profile:typ}]{HOt:rcg:abe:extra}.
This new condition still holds whp:
see
	\cite[\cref{res-p2:cutoff1:typ}]{HOt:rcg:abe:cutoff}
and
	\cite[\cref{res-p5:profile:typ}]{HOt:rcg:abe:extra}.
\end{Proof}

\begin{Proof}[Proof of \cref{res-p1:ext:comp:k-large:S=S'|W=W'}]
Part of the typicality conditions for this regime were that a large number of generators are picked at most once; see \cref{eq-p1:cutoff:d:typ:>=}.
When calculating the probability that $C_{i,j} = C'_{i,j}$ for all $i,j \in [k]$, we restricted ourselves to just looking at those $(i,j)$ with $C_{i,j}, C'_{i,j} \in \bra{0,1}$; see \cref{eq-p1:cutoff:3:>:pick-once,eq-p1:cutoff:3:>:typ}.
There we considered only the directed case, so $C_{i,j} \ge 0$ necessarily.
Here we consider both the directed and undirected simultaneously; in the latter case, $C_{i,j}$ may be negative.
We thus instead consider $(i,j)$ with $C_{i,j}, C'_{i,j} \in \bra{0,\pm1}$, which is the case if $\max\bra{ \abs{W_i}, \abs{W_j}, \abs{W'_i}, \abs{W'_j} } \le 1$.
(The extension from directed to undirected was explained in detail earlier, in \S\ref{sec-p1:cutoff:undir}.)

We consider first the (simpler) case that $d = 3$.
By the previous paragraph,
the event
\[
	\mce \cq \bra{ \bm C = \bm C' }
\Qwhere
	\bm C \cq (C_{i,j})
\Quad{and}
	\bm C' \cq (C'_{i,j})
\]
from \cref{eq-p1:cutoff:3:E-def}
may be replaced with the event
\[
	\widetilde \mce
\cq
	\brb{ \exists \, (i',j') \st \abs{D_{i',j'}} = 1 }^c
\Qwhere
	D_{i,j} \cq C_{i,j} - C'_{i,j}
\Qfor
	i,j \in [k].
\]
Herein assume that there exists $(i',j')$ with $\abs{D_{i',j'}} = 1$.
By \cref{res-p1:cutoff:d:matrix-product},
we need to control $A B$ for $A, B \sim^\iid \Unif(\mbz_m)$,
representing the $Z_{\gamma_{i'}}(a,a+1)$ and $Z_{\gamma_{j'}}(a+1,b)$, respectively.
(When $m = p$ was prime, it did not matter whether or not $\abs{D_{i',j'}} = 1$, just that $D_{i',j'} \not\equiv 0$ mod $p$.)

The law of $A B$ conditional on $\gcd(B, m)$ is uniform on
\(
	\gcd(B, m) \mbz_m.
\)
Hence
\[
	\MAX{x \in \mbz_m} \, \pr{ A B \equiv x \MOD m }
=
	\MAX{x \in \mbz_m} \, \ex{ \pr{ A B \equiv x \mod m \mid \gcd(B, m) } }
\le
	\ex{ \gcd(B, m) }/m.
\]
It remains to bound this expectation.
For $\alpha, \beta \in \mbn$, write $\alpha \wr \beta$ to indicate that $\alpha$ divides $\beta$.
We~have
\[
	\ex{ \gcd(B, m) }
=
	\sumt{r} r \, \prt{ \gcd(B, m) = r }
\le
	\sumt{r} r \, \prt{ r \wr B } \one{ r \wr m }
=
	\sumt{r} \one{ r \wr m }
=
	\div m.
\label{eq-p1:ext:comp:3:ex-gcd}
\nt
\]
We deduce the analogue of \cref{eq-p1:cutoff:3:Ec} from this:
\[
	\MAX{w \in \mcw} \, \pr[w]{ S = S' \midb \exists \, (i',j') \st \abs{D_{i',j'}} = 1 }
\le
	(\div m)/m
\label{eq-p1:ext:comp:3:Sab=S'ab|good}
\nt
\]
where, as previously, we write
\(
	\pr[w]{\cdot}
\cq
	\pr{ \, \cdot \mid W = W' = w }
\)
and $\mcw$ is the `typicality' set.	

We bound
\(
	\pr[w]{ S = S' \mid \widetilde \mce } \le 1.
\)
This leads us to the analogue of \cref{eq-p1:cutoff:3:decomp}, still for $d = 3$:
\[
	\MAX{w \in \mcw} \, \pr[w]{ S = S' \midb \typ }
\le
	(\div m) / m + q
\Qwhere
	q
\cq
	\MAX{w \in \mcw} \, \prt[w]{ \widetilde \mce }.
\label{eq-p1:ext:comp:3:S=S'|W=W':k-large}
\nt
\]
The case $d = 3$ now follows immediately from \cref{res-p1:cutoff:3:prod-bound}.

Consider now general $d$.
We use the partitioning argument as given in \S\ref{sec-p1:cutoff:d} here.
There we controlled each partition separately, having an error probability $q_b$ corresponding to block $P_b$ for each $b = 3, ..., d$.
For each block $P_b$ (ie each $b \in \bra{3, ..., d}$), the event
\[
	\widetilde \mce_b
\cq
	\brb{ \exists \, (i',j') \in P_b^2 \st \abs{D_{i',j'}} = 1 }^c
\quad
	\text{holds with probability at most $q_b$}.
\]
When $\widetilde \mce_b$ does not hold, we use the same analysis as above, leading to \cref{eq-p1:ext:comp:3:Sab=S'ab|good}: we obtain
\[
	\MAX{w \in \mcw} \, \pr[w]{ S_{a,b} = S'_{a,b} \midb \exists \, (i',j') \st \abs{D_{i',j'}} = 1 }
\le
	(\div m)/m.
\label{eq-p1:ext:comp:d:Sab=S'ab|good}
\]
The analogue of \cref{eq-p1:ext:comp:3:S=S'|W=W':k-large} for general $d$ is obtained in exactly the same way as \cref{eq-p1:cutoff:d:decomp}; it is
\[
	\MAX{w \in \mcw} \, \pr[w]{ S = S' \midb \typ }
\le
	\prodt[d]{3}
	\rbb{ \rbb{ (\div m) / p }{}^{b-2} + q_b }
\Quad{where}
	q
\cq
	\MAX{w \in \mcw} \, \prodt[d]{3} \pr[w]{ \widetilde \mce_b }.
\label{eq-p1:ext:comp:k-large:S=S'|W=W':d}
\nt
\]
The lemma now follows immediately from \cref{res-p1:cutoff:d:prod-bound}, noting that $\div m \ge 2$.
\end{Proof}

\begin{Proof}[Proof of \cref{res-p1:ext:comp:k-small:S=S'|W=W'}]
The high-level method applied here is similar to that used for $k \gtrsim \log \abs \uab$.

We start with $d = 3$ for the sake of exposition; we then describe the adaptations to general (including diverging) $d$.
We cannot assume that there is some $(i',j')$ with $\abs{D_{i',j'}} = 1$.
Instead, as in \S\ref{sec-p1:cutoff:3}, we restrict consideration to those generators which are picked order $s \cq t/k$ times;
cf \cref{eq-p1:cutoff:3:<:pick-gen:def,eq-p1:cutoff:3:<:typ}.
It then costs an error probability $\oh{1/p}$ to assume that there exists $(i',j')$, each chosen order $s$ times, such that $D_{i',j'} \not\equiv 0$ mod $m$; cf \cref{eq-p1:cutoff:3:<:restrict-i,eq-p1:cutoff:3:<:q<<1/p}.
Previously $D_{i',j'} \not\equiv 0$ was invertible (as $m = p$ was prime), thus we needed only control $A B$ with $A, B \sim^\iid \Unif(\mbz_m)$ representing the $Z_{\gamma_i}(a,a+1)$ and $Z_{\gamma_j}(a+1,b)$, respectively.
Now we must control $D' A B$ where $D' \cq D_{i',j'}$.

The law of $D' A B$ conditional on $\gcd(D' B, m)$ is uniform on $\gcd(D' B, m) \mbz_m$.
Further,
\[
	\gcd(D' B, m)
\le
	\gcd(D', m) \gcd(B, m)
\Quad{and}
	\mfgcd \cq \gcd(D', m) \ \text{and} \ \gcd(B, m) \ \text{are independent}.
\]
Also, $\gcd(B, m)$ is independent of $(W,W')$
and
\(
	\ex{ \gcd(B, m) }
\le
	\div m,
\)
as in \cref{eq-p1:ext:comp:3:ex-gcd}.
Thus
\begin{subequations}
	\label{eq-p1:ext:comp:k-small:S=S'|W=W':pf}
\[
	\abs \ugr \, \pr{ S = S' \mid W = W', \: \typ }
\lesssim
	e^{h_0} \cdot \ex{ \mfgcd \mid \typ } \cdot (\div m),
\label{eq-p1:ext:comp:k-small:S=S'|W=W':pf:3}
\nt
\]
analogously to \cref{res-p1:ext:comp:k-large:S=S'|W=W'} above.
It remains to control $\ex{ \mfgcd \mid \typ }$, recalling that $\mfgcd = \gcd(D', m)$ where the generator indices $i'$ and $j'$ are both chosen order $s$ times and $D' \not\equiv 0$ mod $m$.

We now turn to general $d$.
Again, we use the same partitioning argument.
The underlying structure of this argument is the same as for $k \gtrsim \log \abs \uab$.
We obtain a pair $(i'_b, j'_b)$ for each block $P_b$ (ie each $b \in \bra{3, ..., d}$) such that $D_{i'_b, j'_b} \not\equiv 0$ mod $m$ and each of $i'_b$ and $j'_b$ is picked order $s$ times.

Set $\mfgcd_b \cq \gcd(D'_b, m)$ where $D'_b \cq D_{i'_b,j'_b}$ for each $b \in \bra{3, ..., d}$.
We obtain
\[
	\abs \ugr \, \pr{ S = S', \: W = W' \mid \typ }
\lesssim
	e^{h_0}
\cdot
	\ex{ \prodt[d]{3} \mfgcd_b^{b-2} \midb \typ }
\cdot
	(\div m)^{d^2},
\label{eq-p1:ext:comp:k-small:S=S'|W=W':pf:d}
\nt
\]
\end{subequations}
via arguments analogous to those used in the proof of \cref{res-p1:ext:comp:k-large:S=S'|W=W'}, particularly \cref{eq-p1:ext:comp:k-large:S=S'|W=W':d}, and \cref{eq-p1:ext:comp:k-small:S=S'|W=W':pf:3}.
\cref{res-p1:ext:comp:k-small:S=S'|W=W'} now follows from \cref{eq-p1:ext:comp:k-small:S=S'|W=W':pf:d} and the fact that the $\mfgcd_b$ are independent. Indeed, $i'_b, j'_b \in P_b$ for each $b$ and $(P_3, ..., P_d)$ is a partition of $[k]$. Thus there is no overlap in indices.
Independence of the gcds is immediate from the independent evolution of $W_i(\cdot)$ for different $i$.
Hence
\[
	\ex{ \prodt[d]{3} \mfgcd_b^{b-2} \midb \typ }
=
	\prodt[d]{3} \ex{ \mfgcd_b^{b-2} \midb \typ }.
\]
It remains to control this product of gcds.
An upper bound of $m^{5 d^3 / k}$ is given by \cref{res-p1:ext:comp:k-small:gcd}.
\end{Proof}

\begin{lem}
\label{res-p1:ext:comp:k-small:gcd}
	In the notation of the above proof,
	\(
		\prodt[d]{3} \ex{ \mfgcd_b^{b-2} \midb \typ }
	\le
		m^{5 d^3 / k}.
	\)
\end{lem}

\begin{Proof}
We use the trivial bound
\(
	\mfgcd
=
	\gcd(D'_b, m)
\le
	\abs{D'_b}.
\)
Since $i'_b$ and $j'_b$ are both chosen order $s$ times, $\abs{D'_b} \lesssim s^2$.
Here
\(
	s
\asymp
	m^{2(d-1)/k}
\gg
	1.
\)
Hence
\(
	\mfgcd_b
\le
	\abs{D'_b}
\le
	m^{5d/k}.
\)
Thus
\[
	\prodt[d]{3} \ex{ \mfgcd_b^{b-2} \midb \typ }
\le
	\prodt[d]{3} (m^{5d/k})^{b-2}
\le
	m^{5 d^3 / k}.
\qedhere
\]
\end{Proof}

\subsection{Cutoff for the Heisenberg Group}
\label{sec-p1:ext:heis}

\renewcommand{\ugr}{\UU_{p,3}}
\renewcommand{\uab}{\UU_{p,3}^\ab}

In this section we analyse the Heisenberg group $\hgr$:
	an element of this group is represented as $(x,y,z)$ with $x,y \in \mbz_p^{d-2}$ and $z \in \mbz_p$;
	multiplication is defined by
	\[
		(x,y,z) \circ (x',y',z')
	\cq
		(x+x', y+y', z+z' + x \cdot y'),
	\]
	where, here and below, $x \cdot y'$ is the usual dot product for vectors in $\mbz_p^{d-2}$.
Observe that $\UU_{p,3} = \HH_{p,3}$.

\renewcommand{\hgr}{\HH}
\renewcommand{\hab}{\HH^\ab}
\renewcommand{\hk}{\HH_k}

\begin{customthm}{\ref{res-p1:intro:ext:heis}}
	Let $\hgr \cq \HH_{m,d}$.
	Suppose
		that $1 \ll \log k \ll \log \abs \hgr$
	and
		$d \ge 3$ is fixed or diverges sufficiently slowly.
	If $m$ is not prime, then require
	\(
		k \gg \max\bra{ \sqrt{\log m}, \: \log \div m }.
	\)
	Then the RW on $\hk^\pm$ exhibits cutoff whp at
	\[
		\tilde t^\pm_*(k,m,d)
	\cq
		\max\brb{t^\pm_0(k, \abs \hab), \: \log_k \abs \hgr}
	\eqsim
	\begin{cases}
		t^\pm_0(k, \abs \hab)
			&\text{when}\quad
		k \le \rbr{ \log \abs \hab }^{2d-3},
	\\
		\log_k \abs \hgr
			&\text{when}\quad
		k \ge \rbr{ \log \abs \hab }^{2d-3}.
	\end{cases}
	\]
\end{customthm}

\renewcommand{\hgr}{\HH_{m,d}}
\renewcommand{\hab}{\HH_{m,d}^\ab}
\renewcommand{\hk}{(\HH_{m,d})_k}

As stated in the introduction,
the case of diverging $d$ is established in \cite[\cref{res-p2:intro:comp:cor:heis}]{HOt:rcg:abe:cutoff}.
In fact, rather more can be shown when $m = p$ is prime; see \cite[\cref{res-p2:intro:comp:cor:special}]{HOt:rcg:abe:cutoff}.

\renewcommand{\hgr}{\HH}
\renewcommand{\hab}{A}

We initially assume that $m = p$ is prime and describe the adaptations required to handle non-prime $m$ at the end.
To ease notation,
e write
\(
	\hgr \cq \HH_{p,d}
\)
and
\(
	\hab \cq \HH_{p,d}^\ab.
\)
Observe that $\hab \cong \mbz_p^{2(d-2)}$.

\begin{Proof}[Lower Bound on Mixing]
The lower bound follows from exactly the same argument as in \S\ref{sec-p1:cutoff:lower}.
\end{Proof}

It remains to consider the upper bound on the mixing time.

\begin{Proof}[Preliminaries]
\qedtriangle
To analyse $\hgr$, we adjust slightly our analysis of $\UU_{p,3}$.
Recall from the start of \S\ref{sec-p1:cutoff:3} that we wrote the $i$-th generator $Z_i$ as $(A_i, B_i, C_i)$ with $(A_i, B_i, C_i)_{i \in [k]}$ a collection of independent $\Unif(\mbz_p)$-s.
We do the same here except that now $A_i, B_i \sim^\iid \Unif(\mbz_p^{d-2})$; still $C_i \sim^\iid \Unif(\mbz_p)$.

The first observation is that now $S(t)$ has exactly the same representation as for $d = 3$:
\begin{gather*}
	S(t)
=
	\rbb{ \sumt[k]1 A_i W_i(t), \: \sumt[k]1 B_i W_i(t), \: \sumt[k]1 C_i W_i(t) + f\rbr{ \bm\alpha, \bm\beta } }
\\
	\text{where}
\quad
	f\rbb{ (\alpha_j)_1^t, \: (\beta_j)_1^t }
\cq
	\sumt[t]{s=1} \beta_s \cdot \sumt[s-1]{r=1} \alpha_r
\end{gather*}
and $(\alpha_1, \beta_1, \gamma_1), ..., (\alpha_N, \beta_N, \gamma_N)$ are the steps taken by the walk $S$ up to time $t$; see \cref{eq-p1:cutoff:3:product,eq-p1:cutoff:3:S}.

Recall that for $\ugr$ there was a `phase transition' in the mixing time $t_*$ at $k = (\log \abs \uab)^3$.
We considered
the phases
separately.
Similarly, for $\hgr$ there is a `phase transition' in
\[
	\tilde t_*(k,p,d)
\cq
	\max\brb{ t_0(k, p^{2d-4}), \: \log_k(p^{2d-3}) }.
\]
Now, however, the location of the transition is different.
Recall that
\(
	t_0(k,N) \eqsim \tfrac{\rho}{\rho-1} \log_k N
\)
when $k = (\log N)^\rho \gg \log N$.
Plugging in $N = \abs \hab = p^{2d-4}$,
we deduce that
the transition occurs at
\[
	k = (\log \abs \hab)^\rho
\Quad{with}
	\tfrac{\rho}{\rho-1} \cdot (2d-4) = 2d-3,
\Quad{ie}
	\rho = 2d-3.
\]
The regime $k \le (\log \abs \uab)^{2d-3}$ in the current set-up thus corresponds to $k \le (\log \abs \uab)^3$ in \S\ref{sec-p1:cutoff:3}.

We require estimates on the behaviour of the entropy of $W$ around the proposed mixing time.
These were provided by \cref{res-p1:cutoff:prelim:h0,res-p1:cutoff:prelim:omega-choice,def-p1:cutoff:prelim:h0} in the set-up of $\ugr$.
We now describe their analogues in the current set-up of $\hgr$.
	Fix $\xi > 0$ and choose $\omega$ with $1 \ll \omega \ll \min\bra{k, \log \abs \hab}$.
	Define $\rho$ such that $k = (\log \abs \hab)^\rho$.
	The analogue of \cref{res-p1:cutoff:prelim:h0}, for which we briefly defer rigorous justification, is the following:
	the entropy at $(1 + \xi) \tilde t_*$ is at least $h_0 + 2 \omega$ where
	\begin{empheq}
	[ left = {h_0 \cq \empheqlbrace} ]
	{alignat = 3}
		&\log \abs \hab
			&\Qwhen&
		\rho \le 2d-3,
	\nonumber
	\\
		&(1 - \tfrac1\rho) \log \abs \hgr
			&\Qwhen&
		\rho \ge 2d-3.
	\nonumber
	\end{empheq}
With these definitions, \cref{res-p1:cutoff:prelim:typ,res-p1:cutoff:prelim:W=W'} follow.
The analogue of \cref{eq-p1:cutoff:3:W=W'} follows from these:
	\[
		\pr{ W = W' \mid \typ }
	\le
		e^{-h} / \pr{\typ}
	=
		e^{-\omega} e^{-h_0} / \pr{\typ},
	\label{eq-p1:ext:heis:W=W'}
	\nt
	\]
where we use exactly the same definition of typicality as for the upper triangular group $\ugr$.

We now justify the analogue of \cref{res-p1:cutoff:prelim:h0}.
The claim for $\rho \le 2d-3$ follows immediately from the generic \cref{res-p1:cutoff:prelim:omega-choice}, as did the corresponding claim for $\rho \ge 1 + 2/(d-2)$ in the $\UU_{p,d}$ set-up.
We now turn to $\rho \ge 2d-3 \ge 3$.
Important in our calculation of the entropy of $W$ at \text{$\tdiam = \log_k \abs \hgr$} is whether $\tdiam/k \ll 1$. Indeed, inspecting \cref{res-p1:ent:t0}---in particular \cref{eq-p1:ent:t0:>} vs \cref{eq-p1:ent:t0:<,eq-p1:ent:t0:=}---one sees that the entropy can take significantly different forms.
We have $k \ge (\log \abs \hab)^3$ and $\log \abs \hab \asymp \log \abs \hgr$ here. Thus $\tdiam/k = \log_k \abs \hgr / k \ll 1$.
In the set-up of $\UU_{p,d}$, our conditions on $d$ enforced this relation.
(More generally, for $\UU_{p,d}$ it is possible to choose $d$ and $k$ with $k \ge (\log \abs{\UU_{p,d}^\ab})^{1 + 2/(d-2)}$
such that $\tdiam / k \gg 1$.
The conditions which we imposed on $d$ forbade this.)
The claim for $\rho \ge 2d-3$ then follows exactly analogously to the claim for $\rho \ge 1 + 2/(d-2)$ in the $\UU_{p,d}$ set-up.
\end{Proof}

\begin{Proof}[Analysis when $W(t) \ne W'(t)$]
\qedtriangle
An argument completely analogous to that used for \cref{eq-p1:cutoff:3:W=/W'} gives
\[
	\pr{ S = S' \mid W \ne W', \: \typ }
=
	1/p^{d-2} \cdot 1/p^{d-2} \cdot 1/p
=
	1/p^{2d-3}
=
	1/\abs \hgr.
\label{eq-p1:ext:heis:S=S'_W=/W'}
\nt
\]
Indeed, the argument there shows that $S' S^{-1} \sim \Unif(H)$ on $\bra{ W(t) \ne W'(t) } \cap \typ$.
\end{Proof}

\begin{Proof}[Analysis when $W(t) = W'(t)$]
\qedtriangle
Here we let
\(
	\mce
\cq
	\bra{ C_{i,j} = C'_{i,j} \: \forall \, i,j \in [k] }.
\)
We argued that
\[
	\pr{ S = S' \mid W = W' = w, \: \mce^c }
\le
	2/p,
\]
\emph{in the set-up of $\UU_{p,3}$},
see \cref{eq-p1:cutoff:3:Ec}.
We now argue that the same inequality holds here, in the set-up for $\hgr$.
Write $D_{i,j} \cq C_{i,j} - C'_{i,j}$.
On the event $\mce^c$, there exist $i',j' \in [k]$ with $D_{i',j'} \not\equiv 0$.
Then
\[
	f(\bm\alpha,\bm\beta) - f(\bm\alpha',\bm\beta')
=
	A_{i'} \cdot \rbb{ D_{i',j'} B_{j'} + \sumt{j \ne j'} D_{i,j} B_j }
+	\sumt{i \ne i'} A_i \cdot \sumt{j} D_{i,j} B_j;
\]
cf \cref{eq-p1:cutoff:3:f-f'}.
Recall that $A_{i'}, B_{j'} \sim \Unif(\mbz_p^{d-2})$ independently.
As in \cref{eq-p1:cutoff:3:U(V+X)+Y:def}, we write this as
\[
	U \cdot (V + X) + Y.
\]
To analyse this we use arguments analogous to those used for \cref{eq-p1:cutoff:3:U(V+X)+Y:pr}.
Crucial is that $V + X \sim \Unif(\mbz_p^{d-2})$ independently of $U \sim \Unif(\mbz_p^{d-2})$.
If $V + X \equiv 0$ (as a vector in $\mbz_p^{d-2}$), then $U \cdot (V + X) \equiv 0$.
On the other hand, conditional on $V + X$ and that there is some $\ell \in [d-2]$ so that $(V + X)_\ell \ne 0$, we have $U \cdot (V + X) \sim \Unif(\mbz_p)$.
(Since $p$ is prime, any non-zero elements of $\mbz_p$ is invertible.)
Thus
\[
	\pr{ U \cdot (V + X) + Y \equiv 0 }
&
\le
	\pr{ U \cdot (V + X) + Y \equiv 0, \: V + X \not\equiv 0 }
+	\pr{ V + X \equiv 0 }
\\&
\le
	1/p + 1/p^{d-2}
\le
	2/p.
\]
This is analogous to \cref{eq-p1:cutoff:3:U(V+X)+Y:pr}.
From this and \cref{eq-p1:ext:heis:W=W'},
we deduce exactly the same inequality as \cref{eq-p1:cutoff:3:final}:
\[
\begin{gathered}
	\pr{ S = S', \: W = W' \mid \typ }
\le
	2 e^{-h} \rbb{ 1/p + q(t) } / \pr{\typ}
\quad
	\text{where}
\quad
	q(t)
\cq
	\MAX{w \in \mcw} \,
	\pr{ \mce \midb W = W' = w }.
\end{gathered}
\]
It remains to analyse $q(t)$.
We establish an analogue of \cref{res-p1:cutoff:3:prod-bound}, namely show that
\[
	\pr{ S = S', \: W = W' \mid \typ }
\le
	2 e^{-h} \rbb{ 1/p + q(t) } / \pr{\typ}
=
	\oh1.
\label{eq-p1:ext:heis:S=S'_W=W'}
\nt
\]
We perform this analysis below.
For the moment, assume that the above relation holds.
	%
\end{Proof}

\begin{Proof}[Conclusion Given Analysis of $q(t)$]
\qedtriangle
Immediate from \cref{eq-p1:ext:heis:S=S'_W=/W',eq-p1:ext:heis:S=S'_W=W'}, we deduce that
\[
	\abs \hgr \,
	\pr{ S = S' \mid \typ }
-	1
=
	\abs \hgr \, \pr{ S = S', \: W \ne W' \mid \typ }
-	1
+	\abs \hgr \, \pr{ S = S', \: W = W' \mid \typ }
=
	\oh1.
\]
Also, $\pr{\typ} = 1 - \oh1$.
The upper bound on mixing follows from the modified $L_2$ calculation.
\end{Proof}


\renewcommand{\hgr}{\HH}
\renewcommand{\hab}{\HH^\ab}
\renewcommand{\ugr}{\UU}
\renewcommand{\uab}{\UU^\ab}

\begin{Proof}[Analysis of $q(t)$]
	%
Write $q \cq q(t)$.
In analogy with \S\ref{sec-p1:cutoff:3},
we argue that
\begin{subequations}
	\label{eq-p1:ext:heis:q}
\begin{empheq}
[ left = {1/p + q \le \empheqlbrace} ]
{alignat = 3}
	&2/p
		&\Qwhen&
	\rho \le 2d-3,
\label{eq-p1:ext:heis:q:<}
\\
	&2 / \abs \hgr^{1/\rho}
		&\Qwhen&
	\rho \ge 2d-3.
\label{eq-p1:ext:heis:q:>}
\end{empheq}
\end{subequations}
Given this, the same manipulations as used in \cref{res-p1:cutoff:3:prod-bound} establish \cref{eq-p1:ext:heis:S=S'_W=W'}.
These details are routine and exactly analogous to those used in \cref{res-p1:cutoff:3:prod-bound}; they are omitted.
We adjust the proof of \cref{res-p1:cutoff:3:prod-bound} to establish \cref{eq-p1:ext:heis:q}; in particular, the regimes and entropy must be adjusted appropriately.
We do not give a full proof, but rather point out one place where the analysis differs.

For the regime $k \ge (\log \abs \hab)^{2d-3}$,
which corresponds to the previous $k \ge (\log \abs \uab)^3$,
using $t_* \ge \tdiam$,
we upper bound
\(
	q \le \abs \ggr^{-1/\rho}
\)
and compare this with $1/p$.
This involves comparing either $\abs \ugr^{-1/\rho} = p^{-3/\rho}$ or $\abs \hgr^{-1/\rho} = p^{-(2d-3)/\rho}$ with $1/p$.
For the regime $k \le (\log \abs \hab)^{2d-3}$,
which corresponds to the previous $k \le (\log \abs \uab)^3$,
using $t \ge t_0(k, \gab)$,
we upper bound
\(
	q \le \abs \gab^{1/(\rho-1)}.
\)
and compare this with $1/p$.
This involves comparing either $\abs \uab^{-1/(\rho-1)} = p^{-2/(\rho-1)}$ or $\abs \hab^{-1/(\rho-1)} = p^{-(2d-4)/(\rho-1)}$.
In all cases, these critical values are precisely the point at which the `phase transition' occurs.
The inequalities in \cref{eq-p1:ext:heis:q} now follow immediately.
	%
\end{Proof}

\begin{Proof}[Adaptations for Composite $m$]
	%
We can adapt the above arguments to obtain results for composite $m$.
Exactly the same adaptations as for $\UU_{p,3}$ apply and the same conditions appear.
\end{Proof}


\renewcommand{\hgr}{\HH_{m,d}}
\renewcommand{\hab}{\HH_{m,d}^\ab}


\subsection{Cutoff Window}
\label{sec-p1:ext:window}

We determine the limit profile in the regime $k \ll \log \abs \uab$.
Here the mixing time is $t_0(k, \abs \uab)$.

\begin{customthm}{\ref{res-p1:intro:ext:profile}}
	Let $p$ be prime and $d \ge 3$.
	Abbreviate $\ugr \cq \UU_{p,d}$.
	Suppose that
		$1 \ll k \ll \log \abs \uab$
	and
		$d \ge 3$ is fixed or diverges sufficiently slowly.
	Then there exist times $(t_\alpha)_{\alpha \in \mbr}$ satisfying
	\[
		t_0 \eqsim k \cdot \tfrac1{2 \pi e} \abs \uab^{2/k},
	\quad
		t_\alpha - t_0 \eqsim \alpha \sqrt 2 t_0 / \sqrt k = \oh{t_0}
	\Qand
		d^\pm_{\uk}(t_\alpha) \eqsim \Psi(\alpha) \ \whp,
	\]
	where $\Psi$ is the standard Gaussian tail, ie $\Psi(\alpha) \cq (2 \pi)^{-1/2} \intt[\infty]{\alpha} e^{-x^2/2} dx$ for $\alpha \in \mbr$.
\end{customthm}

\begin{Proof}
Recall that $t_0$ is determined by the entropy of $W$, which was the expectation of the random variable $Q$; see \cref{def-p1:ent:t0}.
For $\alpha \in \mbr$, we now define $t_\alpha(k, N)$ according to the variations of $Q$:
\[
	\ex{ Q_1(t_\alpha(k, N)) }
=
	\rbb{ \log N + \alpha \sqrt{vk} }/k
\Qwhere
	v \cq \Varb{ Q_1(t_0(k,N)) }.
\]
Analogously to before, we consider the entropic time $t_\alpha \cq t_\alpha(k, \abs \uab)$.
We show in the supplementary material,
in \cite[\S\ref{sec-p0:se}]{HOt:rcg:supp},
that,
for all $\alpha \in \mbr$,
if $1 \ll \omega \ll \sqrt{vk}$ and $1 \ll k \ll \log \abs \uab$,~then%
\[
	t_0 \eqsim k \abs \uab^{2/k},
\quad
	t_\alpha - t_0 \eqsim \alpha \sqrt 2 t_0 / \sqrt k = \oh{t_0}
\Qand
	\pr{ Q(t_\alpha) \le \log \abs \uab \pm \omega } \eqsim \Psi(\alpha),
\]
where $\Psi$ is the standard Gaussian tail;
see \cite[\cref{res-p0:se:t0a,res-p0:se:CLT}]{HOt:rcg:supp}.

Recall from the analysis of the total variation distance \emph{given typicality}, from \S\ref{sec-p1:cutoff:d}, in the regime $k \ll \abs \uab$, that the particular value of $t_0$ is unimportant---changing it by a constant would not affect the proof or the result.
The above says that $t_\alpha \eqsim t_0$ for all $\alpha \in \mbr$.
Hence this contribution is $\oh1$ when $t_0$ is replaced by $t_\alpha$, regardless of $\alpha \in \mbr$.
All that changes is \cref{res-p1:cutoff:prelim:typ}: now
\[
	\pr{ \mu_{t_\alpha}\rbb{W(t_\alpha)} \le e^{-h} }
=
	\pr{ Q(t_\alpha) \ge \log \abs \uab \pm \omega }
\eqsim
	1 - \Psi(\alpha).
\]
Hence exactly the same argument gives
\(
	d_{\uk}(t_\alpha) \eqsim \Psi(\alpha)
\
	\whp.
\)
\end{Proof}

\section{Typical Distance and Diameter}
\label{sec-p1:dist}

\renewcommand{\ugr}{\UU_{p,d}}
\renewcommand{\uab}{\UU_{p,d}^\ab}
\renewcommand{\ucom}{\UU_{p,d}^\com}

This section focusses on distances from a fixed point in the directed random Cayley graph of the upper triangular matrix group $\ugr$, with $p$ prime and $d \ge 3$.
Recall the definition of \textit{typical distance}:
when $G \cq \ugr$ and there are $k$ generators,
for $R \ge 0$ and $\beta \in (0,1)$,
write
\[
	\mcb_k(R) \cq \brb{ x \in G \midb \dist_{G_k}(\id,x) \le R }
\Qand
	\mcd_{k}(\beta) \cq \min\brb{ R \ge 0 \midb \abs{ \mcb_k(R) } \ge \beta \abs G },
\]
emphasising explicitly the dependence on $d$ for the latter statistic.


\subsection{Precise Statement and Remarks}

In this section, we state the more refined version of \cref{res-p1:intro:typdist}.
Again, there are some simple conditions that the parameters must satisfy.

We now state the main result of this section; it is in essence a restatement of \cref{res-p1:intro:typdist}.

\begin{subtheorem}{thm}
	\label{res-p1:dist:res:typdist}

\begin{thm}[Typical Distance for $k \ll \log \abs \uab$]
\label{res-p1:dist:res:typdist:<}
	Let $k, d, p \in \mbn$ with $p$ prime and $d \ge 3$.~%
	Write
	\[
		M^+_k \cq k \abs \uab^{1/k}/e
	\Qand
		M^-_k \cq k \abs \uab^{1/k} / (2e);
	\Quad{recall that}
		\abs \uab = p^{d-1}.
	\]
	Suppose that $1 \ll k \ll \log \abs \uab$ and that $d$ is either a fixed constant or diverges sufficiently~slowly.
%
	For all constants $\beta \in (0,1)$,
	we have
	\[
		\mcd^\pm_{\uk}(\beta) / M^\pm_k \to^\mbp 1
	\Quad{(in probability)}
		\asinf N.
	\]
	Moreover,
		the implicit lower bound holds deterministically, ie for all choices of generators,
	and
		for all upper triangular groups, requiring only $1 \ll k \ll \log \abs \uab$.
\end{thm}

\begin{thm}[Typical Distance for $k \gtrsim \log \abs \uab$]
\label{res-p1:dist:res:typdist:>}
	Let $k, d, p \in \mbn$ with $p$ prime and $d \ge 3$.
	Write
	\[
		M_k
	\cq
		\log \uab) / \log( k / \log \uab )
	=
		\tfrac{\rho}{\rho-1} \log_k \abs \uab 
	\Qwhere
		\rho \cq \log k / \log \log \abs \uab.
	\]
		Suppose
	that $k \gtrsim \log \abs \uab$ and $\log k \ll \log \abs \ugr$
	and
	that $d$ is either a fixed constant or diverges sufficiently~slowly.
	For all $\lambda \in (0,\infty)$, there exists a constant $\alpha^\pm_\lambda \in (0,\infty)$ so that,
	for all constants $\beta \in (0,1)$,
	the following convergences in probability hold:
	\begin{alignat*}{3}
		\mcd^\pm_{\uk}(\beta) / \rbb{ \alpha^\pm_\lambda k} &\to^\mbp 1
	&&\Quad{if}&
		k &\eqsim \lambda \log \abs \uab;
	\\
		\mcd^\pm_{\uk}(\beta) / \max\brb{ M_k, \: \log_k \abs \ugr } &\to^\mbp 1
	&&\Quad{if}&
		k &\gg \log \abs \uab.
	\end{alignat*}
	Moreover,
		the implicit lower bound holds deterministically, ie for all choices of generators,
	and
		for all upper triangular groups, requiring only $k \gtrsim \log \abs \uab$ and $\log k \ll \log \abs \ugr$.
	Note that
	\[
		\max\brb{ M_k, \: \log_k \abs \ugr }
	=
		\max\brb{ \tfrac{\rho}{\rho-1}, \: \tfrac12 d } \log \abs \uab.
	\]
\end{thm}

\end{subtheorem}

We also state the result on the diameter; it is a restatement of \cref{res-p1:intro:diam}.

\begin{thm}[Diameter for $k \gtrsim \log \abs \uab$]
\label{res-p1:dist:res:diam}
	Let $k, d, p \in \mbn$ with $p$ prime and $d \ge 3$.
	Write
	\[
		M_k
	\cq
		\log \uab) / \log( k / \log \uab )
	=
		\tfrac{\rho}{\rho-1} \log_k \abs \uab 
	\Qwhere
		\rho \cq \log k / \log \log \abs \uab.
	\]
	Suppose
	that $k \gtrsim \log \abs \uab$ and $\log k \ll \log \abs \ugr$
	and
	that $d$ is either a fixed constant or diverges sufficiently~slowly.
%
	For all $\lambda \in (0,\infty)$,
	letting $\alpha^\pm_\lambda$ be the constant from \cref{res-p1:dist:res:typdist:>},
	the following convergences in probability hold:
	\begin{alignat*}{3}
		\rbb{ \diam{} \uk } / \rbb{ \alpha^\pm_\lambda k} &\to^\mbp 1
	&&\Quad{if}&
		k &\eqsim \lambda \log \abs \uab;
	\\
		\rbb{ \diam{} \uk } / \max\brb{ M_k, \: \log_k \abs \ugr } &\to^\mbp 1
	&&\Quad{if}&
		k &\gg \log \abs \uab.
	\end{alignat*}
	Moreover,
		the implicit lower bound holds deterministically, ie for all choices of generators,
	and
		for all upper triangular groups, requiring only $k \gtrsim \log \abs \uab$ and $\log k \ll \log \abs \ugr$.
	Note that
	\[
		\max\brb{ M_k, \: \log_k \abs \ugr }
	=
		\max\brb{ \tfrac{\rho}{\rho-1}, \: \tfrac12 d } \log_k \abs \uab.
	\]
\end{thm}

As a proxy for the size of the ($L_1$) balls in the (directed) Cayley graph with $k$ generators, denoted $\mcb_k(\cdot)$, we use the size of discrete, directed $L_1$ balls in dimension $k$, denoted $B_k(\cdot)$:
for $R \ge 0$, define
\(
	B_k(R) \cq \bra{ x \in \mbz_+^k \mid \dist_{\mbz_+^k}(0,x) \le R }.
\)
This is done in \cref{res-p1:dist:size} below.

Were the underlying group Abelian, we would have the easy inequality
\(
	\abs{ \mcb_k(R) } \le \abs{ B_k(R) }.
\)
For the upper triangular group $\ugr$, we develop a similar inequality; roughly, we use the inequality for the Abelianisation
and upper bound the number of elements which can be seen by the other vertices by the maximum amount, ie $\abs \ugr / \abs \uab$.
In \cite[\S\ref{sec-p3:typdist1}]{HOt:rcg:abe:cutoff}, we studied typical distance for general Abelian groups, using the same (overall) method; there, the radius $R$ of the balls in question was chosen so that $\abs{B_k(R)} \approx \abs G$.
Here our candidate radius $M^*_{k}$ satisfies $\abs{B_k(M^*_{k})} \gg \abs \uab$.

\begin{defn}
\label{def-p1:dist:M}
	Set
	\(
		\omega
	\cq
		\max\bra{ (\log k)^2, \: k / \abs \uab^{1/(2k)} }
	=
		\max\bra{ (\log k)^2, \: k / p^{(d-1)/(2k)} }.
	\)
	Choose $M^*_{k}$ to be the minimal integer satisfying $\abs{B_k(M^*_{k})} \ge e^\omega \abs \uab$.
\end{defn}

For the sake or presentation, we first analyse in \S\ref{sec-p1:dist:outline}--\S\ref{sec-p1:dist:upper} typical distance for directed graphs with $1 \ll k \ll \log \abs \uab$.
We then describe in \S\ref{sec-p1:dist:ext} the requisite extensions to the argument to handle undirected graphs and $k \gtrsim \log \abs \uab$ for both typical distance and diameter.

\subsection{Outline of Proof}
\label{sec-p1:dist:outline}

As remarked after the summarised statement (in \S\ref{sec-p1:intro:res}), when considering the mixing time on a graph, geometric properties of the graph are often derived and used.
In a reversal of this, we use knowledge about the mixing properties of the random walk to derive a geometric result; the style of proof is similar enough that we even quote lemmas from the mixing section.

The main difference between the proofs is the following:
	previously, $W(\cdot)$ was a DRW on $\mbz_+^k$;
	we replace this $W(t)$ by $A$ which is uniformly distributed on a $\mbz_+^k$-ball of radius $R$, for $R$ defined later;
	this $A$ tells us how many times each generator is used;
	we apply the sequence of generators, with multiplicities, in an order chosen uniformly at random;
	call the resulting element $S$.

We choose $M$ so that this ball has size slightly larger than $\abs \uab$---recall that this size was used for the entropic time $t_0(k, \abs \uab)$ in the mixing.
For a constant $\xi > 0$, if $R \cq M(1-\xi)$, then we use a counting argument to show that the ball cannot cover more than a proportion $\oh1$ of the vertices of the graph; hence this gives a \emph{deterministic} lower bound, valid \forallZ.
For a constant $\xi > 0$, if $R \cq M(1+\xi)$, then we show that not only does the ball cover (almost) all the graph, but the random variable $S$ is well-mixed whp, in the sense that it is very close to the uniform distribution. From this we deduce that, for a proportion $1-\oh1$ of the vertices, there is a non-zero probability that $S$ is at that vertex, and hence a path to it must exist; furthermore, by choice of $A$, the path must have length at most $R = M(1+\xi)$.
To prove this, we even use an analogous $L_2$ calculation to that used for the mixing, namely \cref{res-p1:cutoff:prelim:D-equiv}.

\subsection{Size of Ball Estimates and Lower Bound}
\label{sec-p1:dist:size}

\begin{lem}
\label{res-p1:dist:size}
	For all $R \ge 0$,
	we have
	\[
		\absb{ B_k(R) }
	=
		\binomt{\floor{R}+k}{k}.
	\]
\end{lem}

\begin{Proof}
	Assume that $R \in \mbn$.
	It is a standard combinatorial identity that
	\[
		\absb{ B_k(R) }
	=
		\absb{ \brb{ \alpha \in \mbz_+^k \midb \sumt[k]{1} \alpha_i \le R } }
	=
		\binomt{R+k}{k}.
	\qedhere
	\]
\end{Proof}

Recall that $M_k = M^+_k = k \abs \uab^{1/k}/e$.
The next lemma shows that the difference between $M^*_{k}$ and $M^*_{k}$ is only in subleading order terms, and so can be absorbed into the error terms.

\begin{lem}
\label{res-p1:dist:M}
	For $k \ll \log \abs \uab$,
	for all constants $\xi \in (0,1)$,
	we have
	\[
		M^*_{k} \le \ceilb{ M_{k} (1 + \xi) }
	\Qand
		\absb{ B_k\rbb{ M_{k} (1 - \xi) } } \ll \abs \uab.
	\]
\end{lem}

\begin{Proof}
\emph{Upper bound.}
Set $M \cq e^\xi k \abs \uab^{1/k} / e$.
By Stirling's approximation, we have
\[
	\binomt{M+k}{k}
\ge
	M^k/k!
\gtrsim
	k^{-1/2} (eM/k)^k
=
	k^{-1/2} e^{\xi k} \abs \uab.
\]
Since $\omega \ll k$, we have $k^{-1/2} e^{\xi k} \gg e^\omega$ and $\binom{M+k}{k} \ge e^\omega \abs \uab$.

\smallskip

\emph{Lower bound.}
Set $M \cq e^{-\xi} k \abs \uab^{1/k} / e$.
Using the inequality $\binom Nk \le (eN/k)^k$, we have
\[
	\binomt{M+k}{k}
\le
	\rbb{ e(M+k)/k }^k
\le
	(eM/k)^k \expb{k^2/M}
\le
	e^{-\xi k} \abs \uab \expb{ e^{1+\xi} k / \abs \uab^{1/k} }.
\]
Since $k \ll \log \abs \uab$,
we have
\(
	k/\abs \uab^{1/k} \ll \xi k
\)
and $\binom{M+k}{k} \ll \abs \uab$.
\end{Proof}

From these, it is straightforward to deduce the lower bound (for all $Z$) in \cref{res-p1:dist:res:typdist:<}.

\begin{Proof}[Proof of Lower Bound in \cref{res-p1:dist:res:typdist:<}]
Were the underlying group Abelian, we would be able to upper bound
\(
	\abs{ \mcb_k(M) }
\le
	\abs{ B_k(M) }.
\)
However, this does not hold for general groups.

Recall that the Abelianisation of $\ugr$ corresponds to `modding out all but the super-diagonal':
\[
	\uab
=
	\ugr/[\ugr,\ugr]
\cong
	\mbz_p^{d-1}
\Quad{and}
	\abs \ucom
=
	\abs{[\ugr,\ugr]}
=
	p^{(d-1)(d-2)/2}.
\]
for a given number of steps, the number of different elements
that can be seen is at most $\abs \ucom$ times the number that can be seen in the Abelianisation
$\uab$;
that is,
\[
	\absb{ \mcb_k(M) }
\le
	L \cdot \abs{ \ucom }
\Qwhere
	L
\cq
	\absb{ \brb{ g \ucom \midb g \in \mcb_k(M) } }.
\]
Since $L \le \abs{ B_k(M) }$,
we have
\(
	\abs{ \mcb_k(M) } \le \abs{ B_k(M) } \abs \ucom.
\)
We choose $M$ so that
\(
	\abs{ B_k(M) } \approx \abs \uab.
\)
More precisely,
for any constant $\xi > 0$,
by \cref{res-p1:dist:M},
we have
\[
	\abs{ B_k\rbb{M^*_{k} (1 - \xi )} }
\ll
	\abs \uab \cdot \abs \ucom
=
	\abs \ugr.
\]
Thus, for any constant $\beta \in (0,1)$, we have $\mcd_{k}(\beta) \ge M^*_{k} (1 - \xi)$, asymptotically.
\end{Proof}

\begin{rmkt}
	This proof generalises further.
	Instead of looking at just upper triangular groups, we can take \emph{any} group $G$.
	We then obtain a lower bound analogously, but where now $M^*_{k}$ is defined so that
	\(
		\abs{ B_k(M^*_{k}) }
	\ge
		e^\omega \abs \gab,
	\)
	for some suitable $\omega \gg 1$.
	(For $\ugr$, this is $e^\omega p^{d-1}$.)
\end{rmkt}

\begin{rmkt}
	The statements are for \emph{directed} lattice balls, in $\mbz_+^k$. Changing to \emph{undirected} lattice balls, in $\mbz^k$, increases the size by a factor at most $2^k$. Since $k \ll \log \abs \uab$ and we are looking at sizes of the order $\abs \uab$, analogous statements can easily be proved for directed balls.
\end{rmkt}

\subsection{Mixing-Type Results and Upper Bound}
\label{sec-p1:dist:upper}

As stated in the outline (\S\ref{sec-p1:dist:outline}), we replace the auxiliary $W$ with $A \sim \Unif(B_k(M^*_{k}))$, and then apply the generators in a uniformly chosen order.
More precisely, we have the following algorithm.

\begin{defn}
\label{def-p1:dist:S-def}
	Define $S$ via the following random algorithm.
	\begin{itemize}[itemsep = 0pt, topsep = \smallskipamount, label = $\bcdot$]
		\item 
		First draw $A \sim \Unif(B_k(M^*_{k}))$;
		this tells us how many times we use each of the $k$ generators.
		Define the vector $g$ by $g_1 = \cdots = g_{A_1} = 1$, $g_{A_1+1} = \cdots = g_{A_1+A_2} = 2$ and so on.
		
		\item 
		To decide in which order we apply the generators, label the steps $1, ..., N$, so $N \cq \sum_1^k A_i$, and then draw a uniform permutation $\sigma$ on $[N] = \bra{1, ..., N}$;
		this will tell us in which order we the generators:
		\(
			S \cq Z_{g_{\sigma(1)}} \cdots Z_{g_{\sigma(N)}}.
		\)
	\end{itemize}
\end{defn}

\noindent
In words, we choose how many times each generator is going to be used by $A$, and then apply them in a uniformly chosen order.
In particular, we can define $(C_{i,j})_{i,j \in [k]}$ as before; see \cref{res-p1:cutoff:d:matrix-product}.

\smallskip

We now present our `mixing-type' result, showing that $S$ is close to uniform.


\begin{prop}
\label{res-p1:dist:L2mixing}
	Assume that the conditions of \cref{res-p1:dist:res:typdist:<} hold.
	Then
	\[
		\ex{ \tvb{ \pr[G_k]{ S \in \cdot } - \pi_G } } = \oh1.
	\]
\end{prop}

\begin{Proof}
For notational ease, write $M \cq M^*_{k}$.
Let $S'$, $A'$ and $\sigma'$ be independent copies of $S$, $A$ and $\sigma$, respectively.
For a set $\mca$ (to be defined), the modified $L_2$ calculation used in \cref{res-p1:cutoff:prelim:D-equiv}~gives
\[
	\ex{ \tvb{ \pr[G_k]{ S = \cdot \mid Z } - \pi_G } }
\le
	n \, \pr{ S = S' \mid \typ } - 1 + \pr{ A \notin \mca },
\]
where $\typ \cq \bra{ A,A' \in \mca }$.
Write $n \cq \abs \ugr = p^{d(d-1)/2}$.
Similarly to the mixing case, we separate according to whether or not $A = A'$.
If $A = A'$, then we do an analysis similar to that of $W = W'$ from \S\ref{sec-p1:cutoff}.
Using \cref{res-p1:a.X-unif} in an analogous way as was used to obtain \cref{eq-p1:cutoff:d:W=/W'}, we obtain
\[
	\pr{S = S' \mid A \ne A', \: \typ }
=
	1/n
=
	1/\abs \ugr
=
	1/ p^{d(d-1)/2}.
\label{eq-p1:dist:A=/A'}
\nt
\]
recall that the coefficients in \cref{res-p1:a.X-unif} (corresponding to the entries of $A - A'$ here) are deterministic, and hence \cref{eq-p1:dist:A=/A'} holds \emph{regardless} of the choice of $\mca$.
This also uses the fact that $M \ll p$, and so $\abs{A_i - A'_i} \ll p$ for all $i$; this follows from manipulating the conditions of \cref{res-p1:dist:res:typdist:<} and using $M \asymp k \abs \uab^{1/k}$.
Using the definition of $M = M^*_{k}$, it is easy to calculate
\[
	\pr{ A = A' \mid \typ }
\le
	\absb{ B_k(M) }^{-1} / \pr{ \typ }
\le
	\abs \uab^{-1} e^{-\omega} / \pr{ \typ };
\label{eq-p1:dist:A=A'}
\nt
\]
this replaces the entropic calculation \cref{eq-p1:cutoff:d:W=W'}.
Combining \cref{eq-p1:dist:A=/A',eq-p1:dist:A=A'} establishes
\[
	n \, \pr{ S = S' \mid \typ } - 1
\le
	n \abs \uab^{-1} \, \pr{ S = S' \mid A = A', \: \typ } / \pr{ \typ }
\]

As stated above, the analysis of $\pr{ S = S' \mid A = A', \: \typ }$ is analogous to the $W = W'$ case from \S\ref{sec-p1:cutoff}.
There we stated that it was not important that $W$ was a DRW, and that we would apply the same proof here (\S\ref{sec-p1:dist}) for a ``different $W$''---the $A$ just defined is this ``different $W$''.
Recall that we separated the generators using the partition $\bra{P_3, ..., P_d}$; we do the same here.~%
Define~$\mce_b$~as~in~\cref{eq-p1:cutoff:d:CbEb}:
\[
	\bm C_b \cq \rbr{C_{i,j}}_{i,j \in P_b},
\quad
	\bm C'_b \cq \rbr{C'_{i,j}}_{i,j \in P_b}
\Qand
	\mce_b
\cq
	\bra{ \bm C_b = \bm C'_b }.
\]
So far, we have in essence been following \textit{Proof of \cref{res-p1:cutoff:res} Given Lemmas \ref{res-p1:cutoff:d:prod-bound} and \ref{res-p1:cutoff:prelim:omega-choice}} from the start of \S\ref{sec-p1:cutoff:d}, but with $W(t)$ replaced by $A$ and $\mcw$ replaced by $\mca$;
it is not until \cref{res-p1:cutoff:d:prod-bound}, which upper bounds the analogue of $\pr{ S = S' \mid A = A', \: \typ }$, that the choice of typicality ($\mcw$ there; $\mca$ here) is made.
Hence \cref{eq-p1:cutoff:d:decomp} still holds here:
write $\widebar\mbp_a\rbr{\cdot} \cq \pr{\cdot \mid A = A' = a, \: \typ}$;
we have
\[
	\widebar\mbp_a\rbb{ S = S' }
\le
	2^{d^2/2} \prodt[d]{3} \rbb{ 1/p^{d-2} + q_b }
\Qwhere
	q_b \cq \maxt{a \in \mca} \prodt[d]{3} \widebar\mbp_a\rbr{\mce_b}.
\label{eq-p1:dist:S=S':int}
\nt
\]

In the mixing context, \cref{res-p1:cutoff:d:prod-bound} upper bounded this probability by $2^{d^2} n^{-1} e^{h_0}$, where $h_0$ was the entropy; for $k \ll \log \abs \uab$, we chose $h_0 = \log \abs \uab$, so this upper bound became $2^{d^2/2} n^{-1} \abs \uab$.
Combined with the entropic calculation \cref{eq-p1:cutoff:d:W=W'}, of which \cref{eq-p1:dist:A=A'} is the analogue, established \cref{eq-p1:cutoff:d:L2:final}:
\(
	n \, \pr{S = S' \mid \typ} - 1
\le
	2 \cdot e^{-\omega} 2^{d^2}.
\)
Conditions on $d$ ensured that we could choose $\omega$ to~make~this~$\oh1$.

We claim that we can copy the proof of \cref{res-p1:cutoff:d:prod-bound} to show that
\(
	q_b \le 1/p^{b-2}.
\)
The proof of this claim is deferred to the end of the subsection (\S\ref{sec-p1:dist:upper}).
From this claim, we deduce that
\[
	\widebar\mbp_a\rbb{ S = S' }
\le
	2^{d^2} \prodt[d]{3} 1/p^{b-2}
=
	2^{d^2} p^{-(d-1)(d-2)/2}
=
	2^{d^2} \abs \uab / n.
\label{eq-p1:dist:S=S':final}
\nt
\]
Combining \cref{eq-p1:dist:S=S':int,eq-p1:dist:S=S':final}, we obtain
\[
	n \, \pr{S = S' \mid \typ } - 1
\le
	e^{-\omega} 2^{d^2} / \pr{\typ},
\]
where we shall choose $\typ$ so that $\pr{\typ} = 1 - \oh1$; this is analogous to \cref{eq-p1:cutoff:d:L2:final}.
We now check that our conditions on $d$ allow us to choose $\omega$ so that $\omega \gg d^2$:
	recall from \cref{def-p1:dist:M} that
	\(
		\omega
	=
		\max\bra{(\log k)^2, k / \abs \uab^{1/(2k)}};
	\)
	the condition $d^2 \ll \omega$ is included in the conditions of \cref{res-p1:dist:res:typdist:<}.

\smallskip

It remains to prove our claim that we can copy of the proof of \cref{res-p1:cutoff:d:prod-bound} to prove that
\(
	q_b \le 1/p^{b-2}.
\)
In said proof, we were particularly interested in the (expected) number of times that an individual generator was picked; this was $s \cq t/k$, and, in the regime $k \ll \log \abs \uab$, satisfied $s \asymp \abs \uab^{2/k}$.
At the start of \textit{Proof of \cref{res-p1:cutoff:d:prod-bound}} for $k \ll \log \abs \uab$, we emphasised that the proof did not rely heavily on the distribution of $W$, nor did it need $s \asymp p^{2(d-1)/k}$; we apply the same arguments with $W$ replaced by $A$, and in this case the expected number of times that an individual generator is picked, which we still denote $s$, satisfies $s \asymp \abs \uab^{1/k}$ since $A \sim \Unif(B_k(M^*_{k}))$ and, by \cref{res-p1:dist:M}, $M^*_{k} \asymp k \abs \uab^{1/k}$.
We elaborate further on how to adapt the proof to this context.

Let $\eta \in (0,1)$ be a (small) constant.
For $a \in \mbz_+^k$, writing
\[
	\mcc(a) \cq \brb{ i \in [k] \midb \eta s \le a_i \le \eta^{-1} s },
\Quad{we have}
	\pr{ \abs{\mcc(A)}/k \ge \tfrac45 } = 1 - \oh1,
\]
if $\eta$ is sufficiently small; this is analogous to \cref{eq-p1:cutoff:d:pick-gen}.
We use this to define typicality, analogously to \cref{eq-p1:cutoff:d:typ:<} except recalling that we no longer require the entropic part:
\[
	\mca
\cq
	\brb{ a \in \mbz_+^k \midb \abs{\mcc(a)} \ge \tfrac45 k, \: \maxt{i} a_i < p }.
\]

We use exactly the same decomposition of generators; we look at $m$-tuples, and require $m \ll k/d^2$. We take $m \asymp d$, and so need $d^3 \ll k$;
this is implied by the conditions of \cref{res-p1:dist:res:typdist:<}.
Fixing some $a \in \mca$, consider the mode $\rho_m$ of the vector $\bm C_m \cq (C_{i,j})_{i,j \in [m]}$, conditional on $A = a$; here $C_{i,j}$ is defined as in \cref{eq-p1:cutoff:3:Cij-def}, but now with $S$ defined using $A$ instead of $W$.
Since $\log s = \log \abs \uab / k \asymp d \log p/k$, as it did in \S\ref{sec-p1:cutoff:d}, we may apply \cref{res-p1:cutoff:d:prod-bound:<:clm}, provided $m \ll \abs \uab^{1/k}$.
This follows from the condition $k \le \tfrac12 \log \abs \uab / \log d$ in \cref{res-p1:dist:res:typdist:<}.
Applying \cref{res-p1:cutoff:d:prod-bound} gives
\[
	q_b \le p^{-c(b-2) m / d},
\]
for some constant $c > 0$,
exactly as in \cref{eq-p1:cutoff:d:qb:final}.
Setting $m \cq d/c$ gives $q-b \le 1/p^{b-2}$.
\end{Proof}

\section{Extensions to Typical Distance Theorem}
\label{sec-p1:dist:ext}

In this section we describe two extensions to the typical distance arguments.

\begin{itemize}[itemsep = 0pt, topsep = \smallskipamount, label = \bcdot]
	\item [\S\ref{sec-p1:dist:ext:k-log-ab}]
	We consider typical distance $k \gtrsim \log \abs \uab$.
	
	\item [\S\ref{sec-p1:dist:ext:diam}]
	We consider diameter for $k \gtrsim \log \abs \uab$.
	
	\item [\S\ref{sec-p1:dist:ext:undir}]
	We consider undirected Cayley graphs for $1 \ll k \lesssim \log \abs \uab$.
\end{itemize}

\subsection{Extending Typical Distance to $k \gtrsim \log \abs \uab$}
\label{sec-p1:dist:ext:k-log-ab}

In this subsection we extend the argument to $k \gtrsim \log \abs \uab$.
We consider first $k \asymp \log \abs \uab$.

\begin{Proof}[Typical Distance for $k \asymp \log \abs \uab$]
Key to the typical distance analysis was adapting the mixing time analysis from \S\ref{sec-p1:cutoff:d} to the case where the auxiliary process $W$, which is a RW, is replaced with $A$, which is uniform on a certain ball.
To define the element $S$ of $G$, we use the generators chosen by $A$ (ie $Z_i$ $A_i$ times), applied in a uniformly chosen order; see \cref{def-p1:dist:S-def}.
We then adapted the mixing analysis of \S\ref{sec-p1:cutoff} to show, for $k \ll d \log p$, that $S$ is well-mixed.

The method for $k \asymp \log \abs \uab$ is exactly the same, except that now we use mixing analysis for $k \asymp \log \abs \uab$.
We define $\mca$ analogously to $\mcw$, given in \cref{eq-p1:cutoff:d:typ:>=}, except without the entropic part:
\[
	\mca
\cq
	\brb{ a \in \mbz_+^k \midb \abs{ J(a) - R e^{-R/k} } \le \tfrac12 \eps R e^{-R/k}, \: \maxt{i} a_i < \sqrt p }
\Qwhere
	J(a)
\cq
	\sumt[k]{1} \one{a_i = 1}.
\]
Write
\(
	\typ \cq \bra{ A, A' \in \mca }.
\)
It then suffices to show that
\[
	\ex{ \normb{ \pr[G_k]{ S \in \cdot \mid A \in \mca } - \pi_G }_2^2 }
=
	n \, \pr{ S = S' \mid A = A', \: \typ } - 1
=
	\oh1.
\]
(Indeed, if the $L_2$ norm is $\oh1$ in expectation, then it is $\oh1$ whp, and in particular the support of $S$ must be a proportion $1 - \oh1$ of the vertices of $G$ whp.)
As in \cref{eq-p1:dist:S=S':int}, we have
\[
	\widebar\mbp_a\rbb{ S = S' }
\le
	2^{d^2/2} \prodt[d]{3} \rbb{ 1/p^{d-2} + q_b }
\Qwhere
	q_b \cq \maxt{a \in \mca} \prodt[d]{3} \widebar\mbp_a\rbr{\mce_b},
\]
recalling that we wrote $\widebar\mbp_a\rbr{\cdot} \cq \pr{\cdot \mid A = A' = a, \: \typ}$.
Now, however, to bound $\widebar\mbp_a\rbr{\mce_b}$, we use the method from \cref{res-p1:cutoff:d:prod-bound} with $k \asymp \log \abs \uab$, rather than $1 \ll k \ll \log \abs \uab$.
There, key to the analysis was to consider only generators which are chosen either once or not at all---we needed such generators to have the same relative order in $S$ as in $S'$.
For more details, see (\ref{eq-p1:cutoff:3:1/J!}, \ref{eq-p1:cutoff:3:1/t!}, \ref{eq-p1:cutoff:d:1/t!}, \ref{eq-p1:cutoff:d:q:>=:=}).
The remainder of the mixing-type proof follows analogously to the regime $1 \ll k \ll \log \abs \uab$.

\smallskip

It remains to describe the adaptations in calculating the minimal radius of a ball of cardinality at least $\abs \uab = p^{d-1}$.
When $k \asymp \log \abs \uab$, it is not difficult to see that we need a radius order $k$; so typically each of the $k$ coordinates gets displaced by an order 1 amount.
We wish to choose $R \asymp k \asymp \log \abs \uab$ so that
\(
	\binomt{R+k}{k} \approx p^{d-1}.
\)
A more refined approximation to $\binom{R+k}{k}$ is needed, with $R \asymp k$:
an application of Stirling's approximation gives
\[
	\log \binomt{R+k}{k}
=
	(R+k) h\rbb{ k/(R+k) } \cdot \rbb{1 + \oh1},
\]
where $h : [0,1] \to [0,1] : q \mapsto - q \log q - (1-q) \log(1-q)$ is the entropy (in nats) of $\Bern(q)$.
Writing $R = \alpha k$,
if $k \eqsim \lambda \log \abs \uab$, then,
motivated by the above display,
we choose $\alpha$ to satisfy
\[
	(\alpha + 1) h\rbb{ 1/(\alpha+1) }
=
	\lambda
\]
in the directed case.
In the undirected case, the same analysis holds, but the lattice balls have slightly different sizes; cf \S\ref{sec-p1:dist:ext:undir}.
The desired radius $R$ still satisfies $R \asymp k$, but now with a different implicit constant; cf $\alpha^\pm$ in \cite[\cref{res-p3:typdist2:res}]{HOt:rcg:abe:geom}.
Further details can be found in \cite[\S\ref{sec-p3:typdist2:balls}]{HOt:rcg:abe:geom} and \cite[\S\ref{sec-p0:balls}]{HOt:rcg:supp}.
(In \cite[\S\ref{sec-p3:typdist2}]{HOt:rcg:supp}, the groups $G$ are Abelian, so $\abs \gab = n$.)
\end{Proof}

The analysis for $k \gg \log \abs \uab$ is almost identical.

\begin{Proof}[Typical Distance for $k \gg \log \abs \uab$]
The mixing time of the RW gives an upper bound on the typical distance.
As noted in \cref{rmk-p1:intro:typdist:interp}, the mixing time agrees with the desired upper bound on the typical distance.
This completes the proof of the upper bound.

We turn to the lower bound.
The lower bound of $\log_k \abs \ugr$ is trivial.
For the other part of the maximum, we project to the Abelianisation $\uab$, as in \S\ref{sec-p1:cutoff:lower}.
As in \cref{res-p1:dist:size},
we have
\(
	\abs{B_k(R)} \le 2^R \binom{R+k}{R}.
\)
It is easy to check that this is $\oh{\abs \uab}$ when $R = (1 - \xi) \cdot \tfrac{\rho}{\rho-1} \log_k \abs \uab$ and $k = \rbr{ \log \abs \uab }^\rho$.
This is done carefully in \cite[\S\ref{sec-p3:typdist3}]{HOt:rcg:abe:geom}. (There the groups are Abelian, so $G = \gab$.)
\end{Proof}

\subsection{Diameter for $k \gtrsim \log \abs \uab$}
\label{sec-p1:dist:ext:diam}

In this subsection, we outline the diameter argument.
Here $d \asymp 1$, so $\log \abs \ugr \asymp \log \abs \uab$.
In \cite[\S\ref{sec-p3:typdist2}]{HOt:rcg:abe:geom} we consider typical distance for Abelian groups in the regime $k \asymp \log \abs G$. (For Abelian groups, $G = \gab$.)
In \cite[\S\ref{sec-p3:diam}]{HOt:rcg:abe:geom} we show carefully how to adapt the argument to the diameter; here we give a sketch of the argument there, adapted slightly to upper triangular groups.

\begin{Proof}[Proof of \cref{res-p1:dist:res:diam}]
Clearly typical distance is a lower bound on the diameter.

\smallskip

We first assume that $k \eqsim \lambda \log \abs \uab$ with $\lambda \in (0,\infty)$.
Our aim is to show that the diameter is, up to smaller order terms, the same as the typical distance; by \cref{res-p1:intro:typdist}, this is $\alpha k$ for some constant $\alpha \cq \alpha^\pm_\lambda$.
We thus let $\xi > 0$ and set $R \cq \alpha k (1 + \xi)$; we show that $\diam{} \uk \le R + 1$~\whp.

Split the generators into two sets:
	$A \cq [Z_1, ..., Z_{(1-\eps)k}]$
and
	$B \cq [Z_{(1-\eps)k+1}, ..., Z_k]$,
with $\eps = \oh1$ to be determined.
Roughly, one first uses the generators in $A$ to connect the identity to the set $S \subseteq \ugr$ of elements of $G$ which can be reached by paths of length at most $R$.
By \cref{res-p1:intro:typdist}, we have $\abs{\ugr \setminus S} \ll \abs \ugr$ whp; assume this herein.
Given an arbitrary $g \in G$, if $g \notin S$ then one uses the remaining generators from $B$ to connect $g$ to $S$ directly, and by extension to the identity.
More precisely, we try to connect $g \in G \setminus S$ to $S$ via $hz = g$ for some $z = Z_i$ for some $i > (1-\eps)k$ and $h \in S$.
The probability that this fails for a given $g$ for all such $Z_i$ and all $h \in S$ is at most $\rbr{ \abs{\ugr \setminus S} / \abs \ugr}^{\eps k}$.
Since $\abs{\ugr \setminus S} \ll \abs \ugr$ and $k \asymp \log \abs \uab \asymp \log \abs \ugr$, we can choose $\eps \to 0$ sufficiently slowly so that this latter probability is $\oh{1/\abs \ugr}$.
By the union bound, the probability that this fails for some such $g$ is $\oh1$.
When this does not fail, $\diam{} \uk \le R + 1$.

Finally, when $k \asymp \log \abs \uab$ replacing $k$ with $(1 - \eps)k$ changes the typical distance by a factor $1 + o_{\tozero \eps}(1)$.
This completes the proof for this regime.

\smallskip

The same argument holds for $k \gg \log \abs \uab$, using the typical distance result of that regime.
\end{Proof}

\subsection{Undirected Cayley Graphs}
\label{sec-p1:dist:ext:undir}

Here we describe the required adaptations to the proof to allow undirected graphs.
The only major difference is that the size of the discrete lattice balls changes:
previously we considered $a \in \mbz_+^k$ with $\sumt[k]{1} a_i \le R$,
while now we consider $a \in \mbz^k$ with $\sumt[k]{1} \abs{a_i} \le R$.
Indeed, besides estimates on the sizes of balls, the only tool required was an adaptation of the mixing proof. Since this proof works for both the directed and undirected cases, the same holds true here.

When $1 \ll k \ll \log \abs \uab$, the desired radius $R$ satisfies $R \gg k$.
Comparing \cite[\cref{res-p0:balls:R:1}]{HOt:rcg:supp} with \cref{res-p1:dist:size} shows that the radius of the desired ball changes by approximately a factor 2.
This is shown carefully in \cite[\S\ref{sec-p3:typdist1:balls}]{HOt:rcg:abe:geom}; there the groups $G$ are Abelian, so $\abs \gab = \abs G$.

When $k \asymp \log \abs \uab$, the desired radius $R$ satisfies $R \asymp k$.
In general there is not an easy closed form for the implicit constant in either the directed or undirected cases; cf \cite[\cref{res-p3:typdist2:res}]{HOt:rcg:abe:geom}.

\section{Concluding Remarks and Open Questions}
\label{sec-p1:conc-rmks}


\begin{itemize}[itemsep = 0pt, topsep = \smallskipamount, label = $\bcdot$]
	\item [\S\ref{sec-p1:conc-rmks:const-k}]
	We discuss some statistics in the regime where $k$ is a fixed constant.
	
	\item [\S\ref{sec-p1:conc-rmks:nil}]
	We discuss very briefly how our methods can be extended to more general nilpotent groups.
	
	\item [\S\ref{sec-p1:conc-rmks:open-conj}]
	To conclude, we discuss some questions which remain open and gives some conjectures.
\end{itemize}

\subsection{Lack of Cutoff when $k$ Is Constant}
\label{sec-p1:conc-rmks:const-k}

%

Throughout the paper we have always been assuming that \toinf k \asinf{\abs G}.
It is natural to ask what happens when $k$ does not diverge.
This case has actually already been covered by \textcite{DSc:growth-rw}, using their concept of \textit{moderate growth}.
There is no cutoff.

Recall that a group $G$ is called \textit{nilpotent of step at most $L$} if its lower central series terminates in the trivial group after at most $L$ steps:
	$G_0 \cq G$ and $G_\ell \cq [G_{\ell-1}, G]$ for $\ell \in \mbn$ with $G_L = \bra{\id}$.

\smallskip

For a Cayley graph $G(Z)$, use the following notation.
	Write $\Delta \cq \diam G(Z)$ for its diameter.
	For the lazy simple random walk on $G(Z)$,
	write
		$\trel \cq \trel\rbr{ G(Z) }$ for the relaxation time (ie inverse of the spectral gap)
	and
		$\tmix \cq \tmix\rbr{ \eps; G(Z) }$ for the (TV) $\eps$-mixing time, for $\eps \in (0,1)$.
When considering sequences $(G_N(Z_{(N)}))_\Ninn$, add an $N$-sub/superscript.

We phrase the result of \textcite{DSc:growth-rw} in our language.

\begin{thm}[{cf \cite[Corollary~5.3]{DSc:growth-rw}}]
\label{res-p1:conc-rmks:const-k}
	Let $(G_N)_\Ninn$ be a sequence of finite, nilpotent groups.
	For each $\Ninn$, let $Z_{(N)}$ be a symmetric generating set for $G_N$ and write $L_N$ for the step of $G_N$.
	Suppose that $\sup_N \abs{Z_{(N)}} < \infty$ and $\sup_N L_N < \infty$.
	Then
	\(
		\tmix^N / k_N
	\lesssim
		\Delta_N^2
	\lesssim
		\trel^N
	\lesssim
		\tmix^N
	\)
	\asinf N;
	in particular, $(\tmix^N)_\Ninn$ does not exhibit the cutoff phenomenon
\end{thm}

We give a very brief exposition of the results of \textcite{DSc:growth-rw}, including the definition of moderate growth, leading to this conclusion in \cite[\S\ref{sec-p5:const-k}]{HOt:rcg:abe:extra}.

\subsection{Extending Arguments from Upper Triangular to Other Nilpotent Groups}
\label{sec-p1:conc-rmks:nil}

Our analysis has focussed on upper triangular matrix groups; these are a canonical class of nilpotent groups.
In the introduction, in \cref{rmk-p1:intro:cutoff:heis-to-general}, we claimed that some of our analysis extends from upper triangular groups to more general nilpotent groups.
This extension is based primarily on observations made by P\'eter Varj\'u during discussions of our work with him.

Recall that we wrote $S$ for the location of the walk and $W$ for its auxiliary variable; let $W'$ be an independent copy of $W$, and define $S'$ correspondingly.
Recall the definition of $C_{i,j}$ from \cref{eq-p1:cutoff:3:Cij-def}.
The difference $C_{i,j} - C'_{i,j}$ has a natural group theoretic interpretation.
Indeed, for a step-2 nilpotent group, one can write
\(
	S' S^{-1}
\)
in the canonical form
\(
	\prodt[k]{1}
	Z_i^{W_i-W'_i}
	Z_i^{V_i}
\cdot
	\prodt{i < j}
	[Z_i, Z_j]^{C_{i,j}-C'_{i,j}};
\)
for general groups, $S' S^{-1}$ can be written like this up to multiplication on the right by a term in~$[[G, G], G]$.
When $W = W'$, we are left just with the product of commutator terms in the expression for $S' S^{-1}$.

For a $p$-group $G$, when $W \ne W'$ mod $p$, one can show that $S' S^{-1}$ is uniformly distributed on $G$.
(In the current article, we used the specific structure of $\ugr$ to reach this conclusion.)

We give a fairly detailed discussion, elaborating on the above points, in \cite[\S\ref{sec-p5:heis-nil}]{HOt:rcg:abe:extra}.

\subsection{Regime with $p$ Small or $k$ Close to $d(G^\ab)$}
\label{sec-p1:conc-rmks:small}

\renewcommand{\ugr}{\UU_{p,d}}
\renewcommand{\uab}{\UU_{p,d}^\ab}
\renewcommand{\uk}{(\UU_{p,d})_k}
\renewcommand{\hgr}{\HH_{p,d}}
\renewcommand{\hab}{\HH_{p,d}^\ab}
\renewcommand{\hk}{(\HH_{p,d})_k}

The discussion in this subsection applies to both the upper triangular group $\ugr$ and the Heisenberg group $\hgr$.
In order for the Cayley graph to generate the group, we trivially need to generate the Abelianisation. This is in fact sufficient, since the group is nilpotent; recall \cref{rmk-p1:intro:typdist:nil-gen}.
For an Abelian group $A$, we wrote $d(A)$ for the minimal size of a generating set of $A$.
In this notation, we need $k \ge d(\gab)$.
We have
\(
	d(\uab) = d-1
\)
and
\(
	d(\hab) = 2d-4.
\)

We always assumed that $d$ is sufficiently smaller than $k$:
	for $\UU$, it was of smaller order;
	for $\HH$, we had $9d < k$.
Additionally---and, perhaps somewhat surprisingly, relatedly---we assumed that $p$ grows sufficiently quickly.
Here we briefly discuss some difficulties in relaxing either of these conditions, ie allowing $p$ small or $k$ closer to $d(\gab)$.
For a fuller description, going into more detail and with potential solutions, we defer the reader to \cite[\S\ref{sec-p5:upd:p-small-k-close-to-d}]{HOt:rcg:abe:extra}.

Since it is an Abelian group, the walk projected to $\gab$ at time $t$ depends only on the final count $W(t)$, not the order in which the generators are applied.
Every element of $\gab$ is of order precisely $p$.
Thus the entropic lower bound also holds if we replace the RW on $\mbz^k$ (ie $W$) with the RW on $\mbz_p^k$---denote this $\widebar W$.
Even when $k - d(\gab) \asymp k$, if $p \asymp 1$ then the entropic time for $\widebar W$ is a constant factor larger than that for $W$.
At least for certain regimes, one needs to consider $\widebar W$.

Another, potentially more substantial, obstacle regards $\abs{C_{i,j}}$.
Crucial in controlling the error $q(t)$, a typicality requirement was that $\abs{C_{i,j}} < p$, at least for many pairs $(i,j)$.
For $k \gtrsim \log \abs \gab$, this is not an issue because an order 1 proportion of generators are picked at most once. (We used this in \S\ref{sec-p1:ext:comp} to relax the requirement that $p$ be prime.)
For $1 \ll k \ll \log \abs \gab$, however, we have
\(
	\ex{\abs{W_1}} \asymp p^{2r/k}
\)
where $\gab \cong \mbz_p^r$ if $k - d(\gab) \asymp k$. Thus if $p \gg 1$ and, eg, $r \ge \tfrac12 k$, then the assumption on $\abs{C_{i,j}}$ will not hold.
(If $p \asymp 1$ then necessarily $k \ge d(\gab) \asymp d \asymp d \log p \asymp \log \abs \gab$.)

For $\ugr$,
we partitioned the generators according to the columns, so as to get some desired independence:
	see the analysis leading up to \cref{eq-p1:cutoff:d:decomp};
	key in this is the factor $2$ in \cref{eq-p1:cutoff:d:Sab=S'ab}.
Taking a product over the partitions means that this factor $2$ becomes $2^{d^2}$.
When $d$ is very large, this causes a serious problem---even when $p \asymp 1$.
(In fact, the smaller $p$ is the harder it is to sufficiently reduce this factor.)
One likely needs to do something more creative than simply partition the generators.

When $k - d(\gab) \ll k$, these problems are exasperated.
The entropic time for $\widebar W$ can then even be of larger order than that for $W$.
More substantial is that, at least in many cases, it seems unavoidable that $\abs{C_{i,j}} > p$.
One needs to study the `wrap around' effect of taking $C_{i,j}$ mod $p$.
If one tries to replace prime $p$ by general $m$, the gcd analysis may get even more complicated.

\smallskip

We conjecture that there is cutoff whp for the RW on $\gk$ at $\max\bra{t_p(k, \abs \gab), \: \log_k \abs \gab}$, where $t_p(k, \abs \gab)$ is the time at which the entropy of RW on $\mbz_p^k$ reaches $\log \abs \gab$.
See Conjecture~\ref{oq-p1:cutoff:p-small} below.

\subsection{Open Questions and Conjectures}
\label{sec-p1:conc-rmks:open-conj}

We close the paper with some questions which are left open.

\newcounter{oq}

\refstepcounter{oq}
\label{oq-p1:cutoff:general-groups}
\subsubsection*{\theoq: Sufficient Conditions for Cutoff for Nilpotent and General Groups}

We have established cutoff for a certain family of non-Abelian groups at a time depending on the algebraic structure of the group.
The main part of the Aldous--Diaconis conjecture remains~open.

\begin{conj-ind}\theoq
	For all groups $G$,
	for $k \gg \log \abs G$ with $\log k \ll \log \abs G$,
	the random walk on $G_k$ exhibits cutoff whp.
\end{conj-ind}

It is natural to ask at which time this cutoff occurs.

\begin{openproblem-ind}\theoq
	Find an expression for the cutoff time in Conjecture~\ref{oq-p1:cutoff:general-groups}.
	
	Find conditions under which this time can be in terms of a few statistics of the group, eg
		the size of the Abelianisation,
		the number of low dimensional irreducible representations
	or
		the size of the largest Abelian subgroup.
\end{openproblem-ind}

In our companion article \cite[\cref{res-p2:intro:tv}]{HOt:rcg:abe:cutoff},
we establish cutoff for all Abelian groups in the regime $1 \ll k \lesssim \log \abs G$.
There are some necessary conditions to generate the group whp.
Using the techniques from there, one can look to extend Open Problem~\ref{oq-p1:cutoff:general-groups} into the regime $1 \ll k \lesssim \log \abs \gab$.

\medskip

Since $t_0(k, \abs{G^\ab})$ is a lower bound on the mixing time, clearly the size of the Abelianisation plays a role.
Relatedly, the size of the largest Abelian subgroup appears to play a role.
Indeed, consider the dihedral group $D_{2n}$ of order $2n$. This has an Abelian subgroup congruent to $\mbz_n$ (corresponding to rotations).
Some basic calculations suggest that the mixing time likely should be the same as that of $\mbz_n$.
This is perhaps related to irreducible representations (\textit{irreps}):
	$D_{2n}$ has at most 4 irreps of dimension 1 (and hence $\abs{D_{2n}^\ab} \le 4$)
and
	all the remaining irreps are of dimension 2.
Instead of just considering $\abs{G^\ab}$, which is the number of irreps of dimension 1, more generally the number of low dimensional irreps (in some precise sense)
is likely to affect the mixing time.

As a starting point, one should perhaps study nilpotent groups.
Finally, we mention work by \textcite{G:quasirandom-groups}, on \textit{quasirandom groups}.
He looks for groups whose (non-trivial) irreps all have high dimension; such groups he describes as being `very far from Abelian'.
The farther a group is from Abelian, the faster we expect its mixing time to be.
Perhaps a similar criterion would be useful for this question of comparing the mixing of the Abelianisation with that of the full group.

\refstepcounter{oq}
\label{oq-p1:cutoff:variant}
\subsubsection*{\theoq: Cutoff for a Variant of the Upper Triangular Model}

There is a well-known variant of the upper triangular model:
	the super-diagonal coordinates are taken modulo $p$ (ie in $\mbz_p$)
and
	the corner coordinate is taken modulo $p^2$ (ie in $\mbz_{p^2}$);
	we denote this $\widetilde \UU_{p,3}$.
It can be extended to $d \times d$ matrices in one of two natural ways:
	coordinates at a distance $\ell$ from the diagonal are taken modulo
	either	$p^\ell$
	or		$p^{2^{\ell-1}}$.
(It is easy to check that both form groups for general $d \ge 3$.)
In some senses, eg growth of balls, this variant is more similar to the (infinite) upper triangular group with entries in $\mbz$.
Note that
	$\abs{\widetilde \UU_{p,3}} = p^4$,
	$\widetilde \UU_{p,3}^\ab \cong \mbz_p^2$
and
	$\widetilde \UU_{p,3}^\com \cong \mbz_{p^2}$.
Suppressing the $p$-dependence, write $\widetilde \UU_k \cq (\widetilde \UU_{p,3})_k$ for the random Cayley graph of $\widetilde \UU_{p,3}$.

\smallskip

Consider the mixing time of the RW on $\widetilde \UU_k$.
We \emph{strongly believe} that our proof---with no additional ideas, just some algebraic checks---implies that $\max\bra{ t_0(k, p^3), \: \log_k(p^4) }$ gives an upper bound whp.
Also, $\log_k(p^4)$ always gives a lower bound.
Thus when $k$ is large enough so that this maximum is asymptotically $\log_k(p^4)$, eg $k \ge (\log p)^4$, we would have cutoff.
What about smaller~$k$?

Our proof from \S\ref{sec-p1:cutoff:lower} gives a lower bound of $t_0(k, p^2)$ whp as $\abs{\widetilde \UU_{p,3}^\ab} = p^2$.
We can thus bound the mixing time between the two entropic times $t_0(k, p^2)$ and $t_0(k, p^3)$.
Which is correct, if either?

\begin{question-ind}\theoq
	Consider the RW on $(\widetilde \UU_{p,3})_k$.
	We conjecture that there is cutoff whenever $1 \ll \log k \ll \log p$.
	Assuming that this is the case, at what time does this cutoff occur?
\end{question-ind}

\refstepcounter{oq}
\label{oq-p1:cutoff:p-small}
\subsubsection*{\theoq: Cutoff for RW on $G_k$ for $G \in \{U_{p,d}, H_{p,d}\}$ with $p$ Small or $k$ Close to $d(G^\ab)$}

\renewcommand{\ugr}{\UU_{p,d}}
\renewcommand{\uab}{\UU_{p,d}^\ab}
\renewcommand{\uk}{(\UU_{p,d})_k}
\renewcommand{\hgr}{\HH_{p,d}}
\renewcommand{\hab}{\HH_{p,d}^\ab}
\renewcommand{\hk}{(\HH_{p,d})_k}

We discussed the case where $p$ is small or $k$ is close to $d$ in \S\ref{sec-p1:conc-rmks:small}.
We consider both $\ggr = \ugr$ and $\ggr = \hgr$.
Recall that $\uab \cong \mbz_p^{d-1}$ and $\hab \cong \mbz_p^{2d-4}$.
Let $t_p(k, N)$ denote the time at which the entropy of RW on $\mbz_p^k$ reaches $\log N$.
Based on the discussion, we make the following conjecture.

\begin{conj-ind}\theoq
	Let $p$ be prime and $d \ge 3$.
	Let $\ggr \in \bra{\ugr, \hgr}$.
	Suppose that $1 \ll \log k \ll \log \abs \ggr$ and $\rbr{ k - d(\gab) } p \gg 1$.
	Then the RW on $\ggr_k$ exhibits cutoff whp at
	\[
		\max\brb{ t_p(k, \abs \gab), \: \log_k \abs \ggr }.
	\]
\end{conj-ind}

We analyse cutoff for the group $\mbz_p^r$ in \cite[\cref{res-p5:intro:p}]{HOt:rcg:abe:extra}.
This analysis will be required, at least for establishing the entropic--Abelianisation lower bound.
We gave details for the upper bound in \S\ref{sec-p1:conc-rmks:small}, including the main obstacles and potential solutions; see \cite[\S\ref{sec-p5:upd:p-small-k-close-to-d}]{HOt:rcg:abe:extra} for a fuller description.

\refstepcounter{oq}
\label{oq-p1:cutoff:p-not-prime}
\subsubsection*{\theoq: Cutoff for Upper Triangular Group $\ugr$ with $p$ Not Prime}

Throughout this paper, we primarily considered $\ugr$ with $p$ prime.
In \cref{res-p1:intro:ext:comp} we relieved this condition, but only in the regime $k \gtrsim \log \abs \uab$.
We conjecture that the analogous behaviour holds for $p$ not prime in the regime $1 \ll k \ll \log \abs \uab$.
This is work in progress.

\begin{conj-ind}\theoq
	Subject to potentially stronger conditions, analogous results hold when $p$ is not prime with the same cutoff time.
\end{conj-ind}


When we studied Abelian groups in our companion article \cite{HOt:rcg:abe:cutoff}, we did not assume that the analogue of $p$ was prime; we did the gcd analysis.
It is not unreasonable to imagine that similar techniques applied there (see \cite[\S\ref{sec-p2:cutoff1:upper}/\S\ref{sec-p2:cutoff2:upper}]{HOt:rcg:abe:cutoff}) may well be applicable here too.

\refstepcounter{oq}
\label{oq-p1:expander}
\subsubsection*{\theoq: Spectral Gap for Upper Triangular Group $\ugr$ with $k \gtrsim \log \abs \uab$}

We studied typical distance for $k \lesssim \log \abs \uab$.
In the boundary regime $k \asymp \log \abs \uab$, the typical distance is order $k$.
We study typical distance for Abelian groups in \cite[\cref{res-p3:intro:typdist}]{HOt:rcg:abe:geom}, obtaining analogous results---for an Abelian group, $G = \gab$. The regime $k \asymp \log \abs G$ is the point at which the Cayley graph of an Abelian group becomes an expander; see \cite[\cref{res-p3:intro:gap}]{HOt:rcg:abe:geom}.
It is natural to conjecture an analogue for upper triangular groups.
It is known that the $G_k$ is an expander whp for \textit{any} group is when, eg, $k \ge 2 \log_2 {}\abs G$; this is the celebrated Alon--Roichman theorem~\cite{AR:cayley-expanders}.

\begin{conj-ind}\theoq
	If $k \gtrsim \log \abs \uab$ with $k - d \asymp k$,
	then $\uk$ is an expander \whp.
\end{conj-ind}

If this were proved for some diverging $d$, then it would provide the first example of a group with the property that its $k$-uniform Cayley graph is an expander whp for some $k$ with $k \ll \log \abs G$.

\refstepcounter{oq}
\label{oq-p1:diameter:k>>1}
\subsubsection*{\theoq: Diameter for Upper Triangular Group $\ugr$ for Diverging $k$}

We have shown concentration of typical distance, but never considered the diameter.
It is trivial that the typical distance is a lower bound on the diameter, and that twice the typical distance is an upper bound.
Further, in the regime $k \asymp \log \abs \uab$, we argued that the diameter and typical distance are asymptotically equivalent; see \S\ref{sec-p1:dist:ext:k-log-ab}.
Can more be determined?

\begin{conj-ind}\theoq
	For $G = \ugr$ and $Z_1, ..., Z_k \sim^\iid \Unif(\ugr)$, write $\Delta_Z$ for the diameter of the Cayley graph with generators $Z$.
	Assume that $k$ diverges, sufficiently rapidly in terms of $d$.
	Then the law of $\Delta_Z$ concentrates.
\end{conj-ind}

\subsubsection*{Questions for Typical Distance}

Questions for typical distance can be asked analogous to those detailed in Questions \ref{oq-p1:cutoff:p-small} and~\ref{oq-p1:cutoff:p-not-prime}.

\subsubsection*{Replacing the Upper Triangular Group with a Nilpotent Group}

Questions~\ref{oq-p1:expander}--\ref{oq-p1:diameter:k>>1} for the upper triangular group can all be extended by replacing $G = \ugr$ with a general nilpotent group $G$; replacing the Abelianisation $\uab$ with $\gab$.

\renewcommand{\bibfont}{\sffamily}
\renewcommand{\bibfont}{\sffamily\small}
\printbibliography[heading=bibintoc]

\end{document}